\documentclass[11pt]{amsart}

\usepackage{amsmath, amscd,hyperref,amssymb}
\usepackage{amsfonts}
\usepackage{amsthm, upref}
\usepackage{graphicx}
\usepackage{enumerate}
\usepackage[dvipsnames]{xcolor}
\usepackage{tikz}
\usepackage{pgf}
\usetikzlibrary{matrix,arrows, decorations.markings}
\usepackage{accents}
\usepackage{subcaption}

\input xy
\xyoption{all}

    \numberwithin{equation}{section}

\usepackage{bm}
        
\bmdefine\alphab{\mathbf{\alpha}}
\bmdefine\betab{\mathbf{\beta}}
\bmdefine\pib{\mathbf{\pi}}
\bmdefine\xib{\mathbf{\xi}}
\bmdefine\sigmab{\mathbf{\sigma}}

\newcommand{\comment}[1]{}
\newcommand{\eq}{\begin{equation}}
\newcommand{\en}{\end{equation}}

\newcommand{\norm}[1]{\left\lVert #1 \right\rVert}

\newcommand{\mcal}[1]{\mathcal{#1}}

\newcommand{\ev}{\mathbb E}

\begin{document}

\theoremstyle{plain}
\newtheorem{thm}{Theorem}[section]
\newtheorem{rem}{Remark}[section]
\newtheorem{lemma}[thm]{Lemma}
\newtheorem{prop}[thm]{Proposition}
\newtheorem{cor}[thm]{Corollary}
\newtheorem{defn}{Definition}[section]

\theoremstyle{definition}
\newtheorem{cond}{Condition}
\newtheorem{asmp}{Assumption}
\newtheorem{notn}{Notation}
\newtheorem{prb}{Problem}

\theoremstyle{remark}
\newtheorem{rmk}{Remark}[section]
\newtheorem{exm}{Example}[section]
\newtheorem{clm}{Claim}

\numberwithin{equation}{section}

\newcommand{\ms}[1]{\color{RedOrange} {#1}}

\title[Dynamical phase transition]{Dynamical phase transition in Brownian motion via sharp large deviations}

\author{Benjamin Budway, Makiko Sasada and Mykhaylo Shkolnikov}
\address{ORFE Department \\ Princeton University \\ Princeton, NJ 08544}
\email{bb2584@princeton.edu}
\address{Graduate School of Mathematical Sciences, The University of Tokyo, 3-8-1 Komaba, Meguro-ku, Tokyo 153-8914, Japan}
\email{sasada@ms.u-tokyo.ac.jp}
\address{Department of Mathematical Sciences and Center for Nonlinear Analysis\\Carnegie Mellon University \\  Pittsburgh, PA 15213}
\email{mshkolni@gmail.com}



\thanks{M.~Shkolnikov is partially supported by the National Science Foundation grant DMS-2342349.}

\date{\today}

\begin{abstract}
We prove a dynamical phase transition in Brownian motion with drift, a phenomenon recently predicted in the physics literature. Specifically, when a Brownian motion with drift on a long time interval is conditioned on an upper large deviation event for the (normalized) occupation measure of the interval $[-1,1]$, in the spirit of the Gibbs conditioning principle, a phase transition occurs. If one conditions on a large deviation event in the far upper tail, the conditioned process follows the Langevin diffusion suggested by the Donsker-Varadhan variational formula for additive functionals. In contrast, if one conditions on a large deviation event in the near upper tail, the conditioned process follows a different Langevin diffusion on a fraction of the time interval and subsequently behaves like the unconditioned process. We establish this phase transition via a new sharp large deviation result for the occupation measure, which is of independent interest. We also provide further corollaries of this result and discuss generalizations to other Markov processes and additive functionals. 
\end{abstract}

\maketitle

\section{Introduction} \label{sec1}

\begin{figure}[h] 
	\centering
	\begin{subfigure}[b]{0.48\textwidth}
		\centering
		\includegraphics[width=\textwidth]{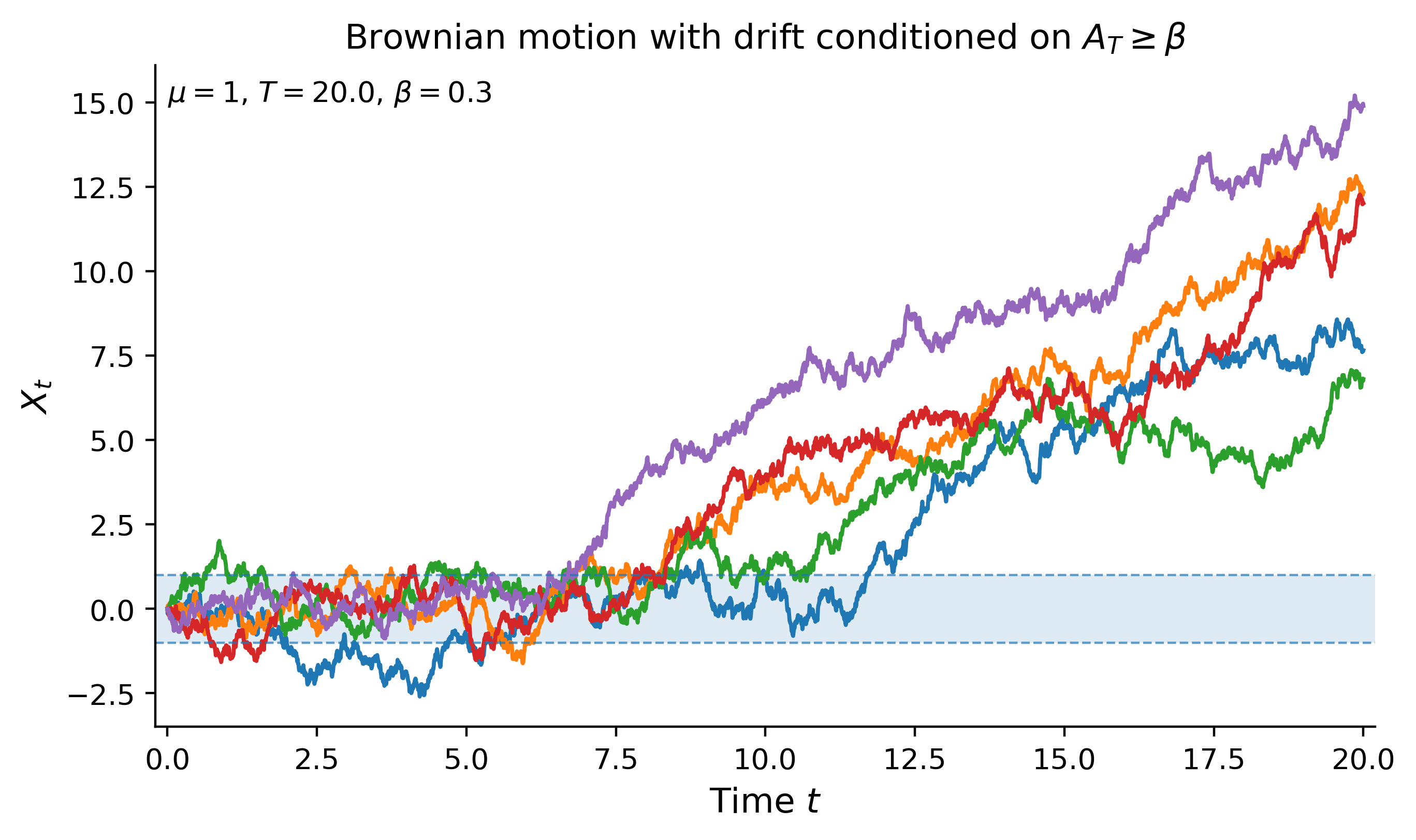}
		\caption{$\beta < \beta_c$}
		\label{fig:path2}
	\end{subfigure}
		\hfill
	\begin{subfigure}[b]{0.48\textwidth}
		\centering
		\includegraphics[width=\textwidth]{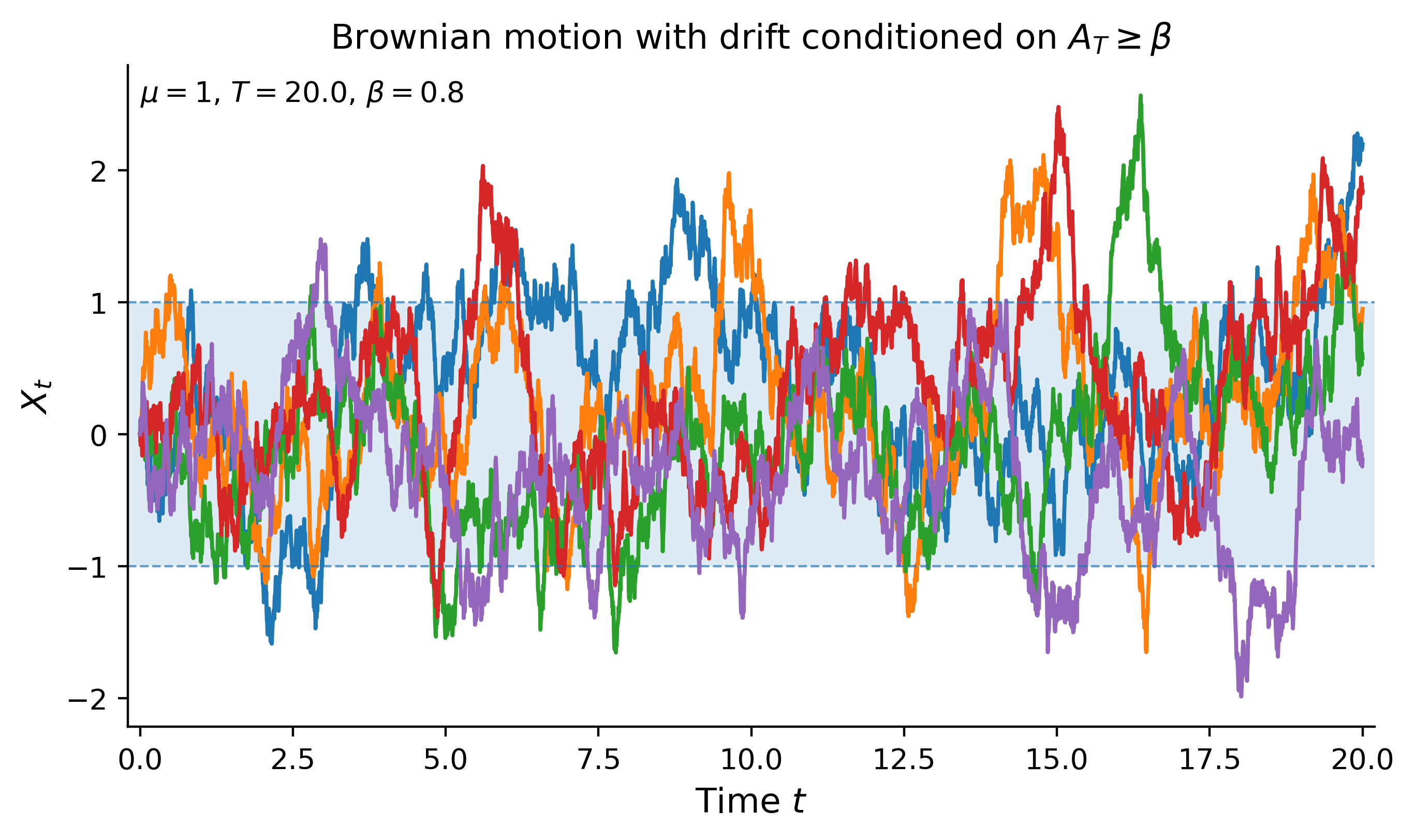}
		\caption{$\beta \geq \beta_c$}
		\label{fig:path1}
	\end{subfigure}
	\caption{Brownian paths conditioned on the occupation constraint.}
	\label{fig:two_paths}
\end{figure}

In the past decade, physicists have studied a phenomenon they call a \textit{dynamical phase transition} in Markov processes, most notably in Brownian motion (see \cite{Tsobgni_Nyawo2016-vw}, \cite{Nyawo2018-ri}, \cite{Mukherjee2024-xt}, \cite{Kanazawa2025-pu}, \cite{Yerrababu2026-zf}, and the references therein). More specifically, consider a Brownian motion with drift $X:=(X_t)_{t\ge0}:=(x+\mu t+B_t)_{t\ge0}$, where $x\in{\mathbb R}$, $\mu\ge0$, and $(B_t)_{t \geq 0}$ is a standard Brownian motion, and let $A_T:=\frac{1}{T}\int_0^T \mathbf{1}_{[-1,1]}(X_t)\,\mathrm{d}t$ denote the fraction of time $X$ spends in the interval $[-1,1]$, for $T\in(0,\infty)$. If $\beta>0$ is less than a certain critical value $\beta_c\ge0$, it has been predicted (see \cite[Section 3]{Nyawo2018-ri}) that, conditional on $\{A_T\ge\beta\}$, the process $(X_t)_{t\in[0,T]}$ acts like a particular Langevin diffusion for an initial fraction of the interval $[0,T]$ and like a Brownian motion with drift $\mu$ thereafter. On the other hand, if $\beta\ge\beta_c$, the conditioned process behaves like a Langevin diffusion depending on $\beta$ on all of $[0,T]$. The simulations in Figure \ref{fig:two_paths} illustrate the dynamical phase transition in Brownian motion with drift. 


\medskip

The primary aim of this paper is to provide a comprehensive proof of the dynamical phase transition in Brownian motion with drift outlined above. To set the stage for our main results, recall that $X$ is reversible with respect to the measure $e^{2\mu x}\,\mathrm{d}x$. Hence, upon conditioning on $\{A_T\ge\beta\}$, one naturally expects the (normalized) occupation measure $\rho_T(\cdot):=\frac{1}{T}\int_0^T \mathbf{1}_{\cdot}(X_t)\,\mathrm{d}t$ to concentrate around a probability measure $\pi$ for which $f:= \big(\frac{\mathrm{d}\pi}{e^{2\mu x}\,\mathrm{d}x}\big)^{1/2}$ is the minimizer in the Donsker-Varadhan variational formula (cf.~\cite{DV1}):
\begin{equation}
	\begin{split}
		&\,\qquad\inf \bigg\{\frac{1}{2}\int_{-\infty}^\infty |f'(x)|^2\, e^{2\mu x}\,\mathrm{d}x:\; f \ge 0,\, \int_{-\infty}^\infty f(x)^2\, e^{2\mu x}\,\mathrm{d}x=1,\, \int_{-1}^1 f(x)^2\, e^{2\mu x}\,\mathrm{d}x = \beta  \bigg\} \\
		&\qquad= \frac{\mu^2}{2} + \inf \bigg\{\frac{1}{2}\int_{-\infty}^\infty |\phi'(x)|^2\,\mathrm{d}x:\; \phi \ge 0,\, \int_{-\infty}^\infty \phi(x)^2\,\mathrm{d}x=1, \int_{-1}^1 \phi(x)^2\,\mathrm{d}x = \beta  \bigg\}.
	\end{split}
\end{equation}
Here, the equality stems from the change of variables $\phi(x) = f(x)\,e^{\mu x}$. The Euler-Lagrange equation corresponding to the right-hand side suggests that the minimizer $\phi$ is an eigenfunction of the Schrödinger operator $\frac{1}{2}\partial_{xx} + \lambda \,\mathbf{1}_{[-1,1]}$ for some Lagrange multiplier $\lambda\ge0$. 

For reasons that will become apparent later, we select the eigenfunction
\begin{equation}\label{Eigenvector}
	\phi_{\lambda}(x) := \begin{cases}
		AMe^{-\sqrt{2\Lambda(\lambda)}|x|} \hspace{39.3pt}\text{ if } |x|\geq 1, \\
		A\cos\big(\sqrt{2(\lambda - \Lambda(\lambda))}\,x\big)\text{ if } |x|\leq 1,
	\end{cases}
\end{equation}
where $M =e^{\sqrt{2\Lambda(\lambda)}}\cos\sqrt{2(\lambda-\Lambda(\lambda))}$, $A$ is fixed by $\int_{-\infty}^\infty \phi_{\lambda}(x)^2\,\mathrm{d}x=1$, and the  associated eigenvalue $\Lambda(\lambda)$ solves the transcendental equation
\begin{equation}\label{trans1}
\Lambda(\lambda)= (\lambda - \Lambda(\lambda))\,\big(\!\tan\sqrt{2(\lambda - \Lambda(\lambda))}\,\big)^2.
\end{equation}
It turns out (Lemma \ref{translemma}) that with the conventions $\Lambda(0):=0$ and that the argument of the tangent function in \eqref{trans1} belongs to $\big(0,\frac{\pi}{2}\big)$ for $\lambda>0$, the equation \eqref{trans1} defines a bijective continuous map $\Lambda:\,[0,\infty)\to[0,\infty)$ satisfying $\Lambda\in C^2((0,\infty))$, $\Lambda'((0,\infty))=(0,1)$, and $\Lambda''>0$ on $(0,\infty)$. Let $\lambda_c \geq 0$ be the unique solution of $\Lambda(\lambda_c) = \frac{\mu^2}{2}$, and let $\beta_c := \Lambda'(\lambda_c)$ (if $\mu =0$, then $\lambda_c= 0$, and we define $\beta_c = 0$). 

\medskip

Finally, we let $\pi_\lambda(\mathrm{d}x):=\phi_\lambda(x)^2\,\mathrm{d}x$ (recall \eqref{Eigenvector}), and, for $\lambda>0$, define $Y^\lambda$ as the Langevin diffusion corresponding to $\pi_\lambda$:
\begin{equation}\label{newsde}
	\mathrm{d}Y^\lambda_t = \partial_x \log \phi_\lambda (Y^\lambda_t)\,\mathrm{d}t + \mathrm{d}B_t
\end{equation}
(see Proposition \ref{SDEfacts} for the well-posedness of \eqref{newsde}). If $Y^\lambda_0$ has law $\pi$ on ${\mathbb R}$, we write $\mathbb{Q}_\lambda^\pi$ for the law of $Y^\lambda$; if $\pi=\delta_x$ for some $x\in{\mathbb R}$, we use $\mathbb{Q}_\lambda^x$ instead. We are now ready to state the dynamical phase transition in Brownian motion with drift $X$ as a theorem. 

\begin{thm}\label{MainTheorem}
Let $(X^{\beta,T}_t)_{0 \leq t \leq T}$ denote the process $(X_t)_{0 \leq t \leq T}$ conditioned on the event $\{A_T\ge\beta\}$. Fix any $\alpha \in [0,1)$ and $S>0$. 
	\begin{enumerate}[(a)]
		\item  If $\beta \in (\beta_c,1)$, define $\lambda$ by $\Lambda'(\lambda)=\beta$. Then, as $T\to\infty$, $(X_{\alpha T  + t}^{\beta,T})_{0 \leq t \leq S}$ converges in distribution to $\mathbb{Q}_\lambda^x$ if $\alpha =0$ and to $\mathbb{Q}_\lambda^{\pi_\lambda}$ if $\alpha > 0$.
		\item If $\beta \in (0,\beta_c]$, then the latter convergences in distribution hold with $\lambda=\lambda_c$ if $\alpha < \frac{\beta}{\beta_c}$. For $\alpha > \frac{\beta}{\beta_c}$, the process $(X_{\alpha T + t}^{\beta,T}-X_{\alpha T}^{\beta,T})_{0 \leq t \leq S}$ converges in distribution to $(\mu t+B_t)_{0\le t\le S}$.
	\end{enumerate}
\end{thm}

\begin{rem} \label{rem: stronger}
\begin{enumerate}[(a)]
\item We prove a stronger result for $\beta\in(0,\beta_c)$: For any $0 \le \alpha_1 \!<\! \frac{\beta}{\beta_c} \!<\! \alpha_2 \!<\! 1$ and $S_1,S_2>0$, the processes $(X^{\beta,T}_{t})_{\alpha_1 T \leq t \leq \alpha_1 T + S_1}$ and $(X_t^{\beta,T}-X_{\alpha_2 T}^{\beta,T})_{\alpha_2 T \leq t \leq \alpha_2 T + S_2}$ converge jointly in distribution to a product measure with the aforementioned marginals.
\item Theorem \ref{MainTheorem} concerns the non-degenerate cases $\beta \in (0,1)$. The degenerate case $\beta=1$, when $(X_t)_{0 \leq t \leq T}$ is conditioned to remain in $[-1,1]$, is partially addressed by the main theorem of \cite{Pinsky}.
\end{enumerate}
\end{rem}

Theorem \ref{MainTheorem} can be regarded as a \textit{Gibbs Conditioning Principle} for Brownian motion with drift. Informally, a Gibbs Conditioning Principle (see, e.g., \cite[Sections 3.3, 3.8, 7.3, 7.4]{Dembo2009-yr}) describes the distribution of a small part in a large random system when the entire system is conditioned on a large deviation event. A classical instance asserts that, under appropriate assumptions, if one conditions a sequence $(X_i)_{i=1}^N$ of i.i.d.\ random variables on its sample average $\frac{1}{N}\sum_{i=1}^N X_i$ lying in a set away from ${\mathbb E}[X_1]$, then for any fixed $k$, the subsequence $(X_i)_{i=1}^k$ converges in distribution to the $k$-fold product of an exponentially tilted distribution as $N\to\infty$. Standard proofs in the i.i.d.\ setting rely on exchangeability, which is lacking in Markov chains/processes. To our knowledge, the only general theorem for the latter is \cite[Theorem 2]{Csiszar1987-pm}, which establishes the principle for Markov chains on finite state spaces via the method of types, a technique inapplicable in our continuous state-space setting. Some related findings for Markov chains on general state spaces are available in \cite{Schroeder} and \cite{Meda1999}. However, the main results in these papers concern either a symmetrized quantity of the Markov chain (see \cite[Theorem 21.3.1]{Meda1999}), or prove an asymptotic quasi-Markov property (see \cite[Equation (1.10), Corollary 5.3]{Schroeder}) which, due to the lack of exchangeability, do not imply a Gibbs Conditioning Principle (see, e.g., \cite[proof of Theorem 7.3.21]{Dembo2009-yr} for the essential role of exchangeability). In the mathematical physics literature, the paper \cite{Chetrite2015-xx} establishes, under suitable assumptions, an equivalence at the large deviation scale between a Markov process conditioned on a large deviation event and a corresponding unconditioned Markov process. This equivalence does not imply the types of convergence in Theorem \ref{MainTheorem}.

\medskip

Our second main result shows the dynamical phase transition in Brownian motion with drift on the level of the (normalized) occupation measures $\rho_t(\cdot)=\frac{1}{t}\int_0^t \mathbf{1}_{\cdot}(X_s)\,\mathrm{d}s$, $t\in(0,\infty)$. 

\begin{thm}\label{Macroscopic} 
For $t\in(0,T]$, let $\rho_t^{\beta,T}$ denote the measure $\rho_t$ conditioned on  $\{A_T\ge\beta\}$.  
\begin{enumerate}[(a)]
		\item If $\beta \in (\beta_c,1)$, define $\lambda$ by $\Lambda'(\lambda) = \beta$. Then, for any $\alpha \in (0,1]$, it holds $\lim_{T\to\infty} \rho^{\beta,T}_{\alpha T}=\pi_\lambda$ weakly in probability.
		\item If $\beta \in (0, \beta_c]$, then it holds $\lim_{T\to\infty} \rho^{\beta,T}_{\alpha T}=\pi_\lambda$  weakly in probability with $\lambda = \lambda_c$ as long as $\alpha \in \big(0,\frac{\beta}{\beta_c}\big]$. If $\alpha \in\big(\frac{\beta}{\beta_c},1\big]$, 
		one has instead $\lim_{T\to\infty} \rho^{\beta,T}_{\alpha T}=\frac{\beta}{\alpha \beta_c}\pi_\lambda$ vaguely in probability with $\lambda=\lambda_c$.
\end{enumerate}
\end{thm}

Thus,  if $\beta \in (\beta_c,1)$, the conditioned occupation measures $\rho^{\beta,T}_{\alpha T}$ behave according to the Donsker-Varadhan theory (cf.~\cite{DV1}). In contrast, if $\beta\in(0,\beta_c)$, then the conditioned occupation measures $\rho^{\beta,T}_{\alpha T}$ are only tight for sufficiently small $\alpha$. For sufficiently small $\alpha$, the conditioned occupation measures do not converge to the measure suggested by the Donsker-Varadhan theory, but always to $\pi_{\lambda_c}$. 

\begin{rem} \label{extra1} 
For $\beta\in(\beta_c,1)$, $\alpha\in (0,1)$ and $\beta\in(0,\beta_c]$, $\alpha\in\big(0,\frac{\beta}{\beta_c}\big)$, the proofs of Theorems \ref{MainTheorem}, \ref{Macroscopic} also reveal that for any events $E_T$ belonging to the $\sigma$-algebras at times $\alpha T$, respectively, and satisfying $\lim_{T\to\infty} \mathbb{Q}_\lambda^x(E_T)=1$, it holds 
$\lim_{T\to\infty} \mathbb{P}(X^{\beta,T}\in E_T)=1$.
\end{rem}

One of the main challenges to proving conditional limit theorems in the context of large deviations is the overwhelmingly small probability of the event we are conditioning on. This means that expanding a sequence of conditional probabilities $\mathbb{P}(A|B_T)$ as $\frac{\mathbb{P}(A\cap B_T)}{\mathbb{P}(B_T)}$ can lead to major difficulties because $\mathbb{P}(B_T)$ is vanishing rapidly; in order to get the correct limit, any estimates used must be remarkably tight. Therefore, to prove Theorems \ref{MainTheorem} and \ref{Macroscopic}, we first obtain \textit{sharp} asymptotics of rare event probabilities in the spirit of the classical results by Bahadur and Ranga Rao \cite{R_r1960-sv}. We now present the sharp asymptotics as our third main result. Henceforth, when we say that $a_T(c) \sim b_T(c)$ uniformly over $c\in C_T$, we mean that $\lim_{T\to\infty}\sup_{c \in C_T} \big|\frac{a_T(c)}{b_T(c)}-1\big|=0$.

\begin{thm}\label{SharpAsymptotics} 
Fix $\mu \geq 0$, a compact $K\subset\mathbb{R}$, and $M\in(0,\infty)$.
	\begin{enumerate}[(a)]
		\item If $\beta \in (\beta_c,1)$, let $\lambda$ be given by $\Lambda'(\lambda)=\beta$. Then, uniformly over $x\in K$ and $|y|\leq M\sqrt{T}$,
		\begin{equation}\label{superasymptotics}
			\mathbb{P}\Big(A_T \ge \beta  + \frac{y}{T}\Big) 
			\sim e^{-\mu x}\phi_{\lambda}(x)e^{-T(\lambda \beta - \Lambda(\lambda)+\frac{\mu^2}{2})}\frac{e^{-\lambda y}}{\Lambda(\lambda)-\frac{\mu^2}{2}}\frac{ e^{-\frac{y^2}{2\Lambda''(\lambda) T}}}{\sqrt{2\pi T \Lambda''(\lambda)}}\int_{-1}^1 \phi_{\lambda}(r)e^{\mu r}\mathrm{d}r.
		\end{equation}
		\item If $\beta = \beta_c \in (0,1)$, then, uniformly over $x\in K$ and $|y|\leq M\sqrt{T}$,
		\begin{equation}\label{critasymptotics}
			\mathbb{P}\Big(A_T \ge \beta_c  + \frac{y}{T}\Big) \sim e^{-\mu x}\phi_{\lambda_c}(x)e^{-T\lambda_c \beta_c  }\frac{e^{-\lambda_c y}}{\beta_c}\Phi\bigg(\!-\frac{y}{\sqrt{T\Lambda''(\lambda_c)}}\bigg)
			\int_{-1}^1 \phi_{\lambda_c}(r)e^{\mu r}\mathrm{d}r,
		\end{equation}
		where $\Phi$ is the standard Gaussian cumulative distribution function. 
		\item If $\beta \in (0,\beta_c)$, then, uniformly over $x\in K$ and $|y|\leq M\sqrt{T}$,
		\begin{equation}\label{subasymptotics}
			\mathbb{P}\Big(A_T \ge \beta  + \frac{y}{T}\Big) 
			\sim e^{-\mu x}\phi_{\lambda_c}(x)e^{-T\lambda_c \beta  }\frac{e^{-\lambda_c y}}{\beta_c}\int_{-1}^1 \phi_{\lambda_c}(r)e^{\mu r}\mathrm{d}r.
		\end{equation}
	\end{enumerate}
\end{thm}

\begin{rem}\label{sharprmk}
	\begin{enumerate}[(a)]
		\item Our proof of Theorem \ref{SharpAsymptotics} can also be used when $A_T$ is replaced by $\frac{1}{T}\int_0^T V(X_t)\,\mathrm{d}t$, where $V:\,{\mathbb R}\to[0,\infty)$ is a bounded compactly supported function that is not identically zero.
		This extension is discussed in the appendix. 
        \item Theorem \ref{SharpAsymptotics} implies a large deviations principle for $A_T$ when $T\to\infty$, and the good rate function is given by 
        the Fenchel-Legendre transform of $J_{\mu}(\lambda) := \big(\Lambda(\lambda)-\frac{\mu^2}{2}\big)^+$. It is worth mentioning that, due to the non-differentiability of $J_{\mu}$ for $\mu > 0$, such a large deviations principle cannot be a direct consequence of the Gärtner-Ellis Theorem (see, e.g., \cite[Theorem 2.3.6]{Dembo2009-yr}). We provide the short proof in Section \ref{corollaries}.  
		\item Another interesting corollary of Theorem \ref{SharpAsymptotics} is that the law of $T(A_T-\beta)$ given $\{A_T\!\ge\!\beta\}$ converges to an exponential distribution. We give a formal statement in Section \ref{corollaries}. 
		\item Yet another interesting corollary of Theorem \ref{SharpAsymptotics} is concerned with the stopping time $\tau_T := \inf\{t \ge 0 : tA_t \ge \beta T\}$. Specifically, we show in Section \ref{corollaries} that the normalized law of $\tau_T$ given $\{A_T\ge\beta\}$ tends to a (one-sided if $\beta=\beta_c$) Gaussian distribution when $\beta \in (0,\beta_c]$. When $\beta \in (\beta_c,1)$, the law of $T- \tau_T$ is shown to converge to an exponential distribution.
	\end{enumerate}
\end{rem}

Let us briefly compare Theorem \ref{SharpAsymptotics} with the classical result of Bahadur and Ranga Rao. Their result (\cite[Theorem 1]{R_r1960-sv}, see also \cite[Theorem 3.7.4]{Dembo2009-yr}) states that for integer $T\in(0,\infty)$ and i.i.d.\ non-lattice random variables $(X_i)_{i=1}^T$ with finite logarithmic moment generating function $\Lambda_1(\lambda) := \log \mathbb{E}[e^{\lambda X_1}]$, it holds
\begin{equation}\label{BRR}
	\mathbb{P}\Big(\frac{1}{T}\sum_{i=1}^T X_i \geq \beta \Big) \sim  \frac{e^{-T(\lambda \beta -\Lambda_1(\lambda)) }}{\lambda \sqrt{2\pi T \Lambda_1''(\lambda)}},
\end{equation}
where $\lambda>0$ is chosen such that $\Lambda_1'(\lambda) = \beta$. If one were to shift from $\beta$ to $\beta+\frac{y}{T}$ in \eqref{BRR}, the exponential factor would yield the same correction terms in $y$ as in \eqref{superasymptotics}, and these corrections could be guessed via a Taylor expansion of the large deviations rate function. On the other hand, the mean parameter $\mu$ enters \eqref{superasymptotics} differently than it would \eqref{BRR}, and there is no obvious analogue of the initial condition $x$ in \eqref{BRR}. Another feature of the asymptotics \eqref{critasymptotics} and \eqref{subasymptotics} is that the dependence both on $T$ and on $y$ differs from the classical setting. The results of Bahadur and Ranga Rao were generalized to arbitrary sequences of random variables in \cite{Chaganty1993-fe} under assumptions on the associated characteristic functions. Sharp large deviation results analogous to \eqref{superasymptotics} for additive functionals of \textit{ergodic} Markov processes have appeared previously, for example, in \cite{Kontoyiannis2003-bl} and \cite{Fernando2021-dg}. To the best of our knowledge, however, no analogues of \eqref{critasymptotics} and \eqref{subasymptotics} have appeared in the literature on large deviations for additive functionals of Markov processes. Beyond Markov processes, sharp large deviations have been used in other contexts, such as asymptotic convex geometry \cite{Liao2024-jz}.

\medskip

The main technical contribution of our paper is the proof of Theorem \ref{SharpAsymptotics}, whose key ideas we now outline. Taking expectations in the elementary formula $\mathbf{1}_{(0,\infty)}(a)+\mathbf{1}_{\{0\}}(a)/2
=\frac{1}{2\pi} \int_{\mathbb R} \frac{e^{(\lambda+i\xi)a}}{\lambda+i\xi}\,\mathrm{d}\xi$, $\lambda>0$ with $a:=TA_T-\beta T-y$, and assuming that ${\mathbb P}(TA_T-\beta T-y=0)=0$ and that Fubini's Theorem is applicable, we anticipate the identity
\begin{equation}\label{Laplace inversion}
\mathbb{P}\Big(A_T \ge \beta  + \frac{y}{T}\Big) 
=\frac{1}{2\pi} \int_{\mathbb R} \frac{1}{\lambda+i\xi}\,\ev\Big[e^{(\lambda+i\xi)\int_0^T \mathbf{1}_{[-1,1]}(X_t)\,\mathrm{d}t}\Big]\, e^{-(\lambda+i\xi)(\beta T+y)}\,\mathrm{d}\xi.
\end{equation}
Here, we define $\lambda>0$ by $\Lambda'(\lambda)=\beta$ for $\beta\in(\beta_c,1)$ and by $\lambda=\lambda_c$ for $\beta\in(0,\beta_c]$. To avoid the technical issues associated with the right-hand side of \eqref{Laplace inversion}, we replace $\frac{1}{\sqrt{2\pi}} \frac{1}{\lambda+i\xi}$, the Fourier transform of $e^{-\lambda a}\,\mathbf{1}_{(0,\infty)}(a)$, in \eqref{Laplace inversion} by a generic Fourier transform $\widehat{\psi}\in C_c({\mathbb R})$. This leads us to study the $T\to\infty$ behavior of the operators $S^{\mu}_{\lambda + i \xi}(T)$ on $L^2(\mathbb{R},e^{2\mu x}\,\mathrm{d}x)$:
\begin{equation}
S^{\mu}_{\lambda+i\xi}(T)f := \ev\Big[e^{(\lambda+i\xi)\int_0^T \mathbf{1}_{[-1,1]}(X_t)\,\mathrm{d}t}\,f(X_T)\Big].
\end{equation}
For small $|\xi|$, which are responsible for the leading order contribution to the right-hand side of \eqref{Laplace inversion}, we apply Kato's perturbation theory \cite{Kato1995-kx}, specifically the spectral theory for holomorphic families of operators of type (A). We conclude by bounding the contribution of intermediate $|\xi|$ using the language of analytic semigroups (\cite[Definition II.4.5]{Engel1999-tt}, which we recall below), and by approximating $e^{-\lambda a}\,\mathbf{1}_{(0,\infty)}(a)$ with functions $\psi$ satisfying $\widehat{\psi}\in C_c^{\infty}(\mathbb R)$. We remark that perturbation theory also underlies the proofs of the sharp asymptotics in \cite[Theorem 6.3]{Kontoyiannis2003-bl} and \cite[Theorem 2.3]{Fernando2021-dg}. 

\medskip

We briefly discuss why the proof strategy of Pinsky \cite{Pinsky} fails in our setting. In his work, diffusion processes are conditioned to stay in a compact set, and he proves a conditional limit theorem similar to Theorem \ref{MainTheorem} using sharp asymptotics. This corresponds to the case $\beta=1$ in our notation, with the compact set being the interval $[-1,1]$. If the diffusion is reversible, $u(t,x)=\mathbb{P}(A_t \geq 1\,|\, X_0=x)$ satisfies the backward Kolmogorov equation $u_t = Lu$, where $L$ is a self-adjoint operator with compact resolvent. The sharp asymptotics then follow directly from the spectral decomposition of $L$. For $\beta \in (0,1)$, however, the natural backward Kolmogorov equation involves an augmented state variable recording the accumulated occupation time; the resulting operator contains a transport term in that variable, and is not naturally self-adjoint, even when the original diffusion is reversible.

\medskip

While the focus of our paper is on the foundational setting of Brownian motion with drift, our results should extend naturally to a broader class of Markov processes. Let $\mcal{L}$ be the generator of a Markov semigroup on an abstract $L^2(E,m)$ space. Our proofs rely heavily on the self-adjointness of $\mcal{L}$. Additionally, certain regularity assumptions on the semigroup are required (e.g., to ensure that pointwise evaluation is well-defined). Regarding the dynamical phase transition, our arguments indicate that recurrent Markov processes -- such as one-dimensional standard Brownian motion -- do not exhibit one. On the other hand, if the spectrum of $\mcal{L}$ is contained in $(-\infty,-\varepsilon]$ for some $\varepsilon>0$, as is the case for Brownian motion with drift, our arguments suggest that a dynamical phase transition occurs. In Lemma \ref{Generalspectrum} of the appendix, we prove that for Markov processes satisfying this condition and $\lambda \geq \lambda_c$, the corresponding Schrödinger operators have an isolated simple eigenvalue at the top of their spectrum, which equals $0$ when $\lambda=\lambda_c$. This result provides the primary spectral input required for the perturbation theory. However, we do not expect our arguments to address other cases of interest, such as standard Brownian motion in dimension $d \geq 5$, where a dynamical phase transition is also anticipated (see the physics paper \cite{Kanazawa2025-pu}).

\medskip

The rest of the paper is structured as follows. In Section \ref{sec:prelim}, we develop the spectral theory framework for our proof of Theorem \ref{SharpAsymptotics} and analyze the stochastic differential equation (SDE) \eqref{newsde}, as well as the transcendental equation \eqref{trans1}. We then prove Theorem \ref{SharpAsymptotics}(a) in Section \ref{sec: super}, and Theorem \ref{SharpAsymptotics}(b), (c) in Section \ref{sec: sub}. In Section \ref{sec: main}, we deduce Theorems \ref{MainTheorem} and \ref{Macroscopic} from Theorem \ref{SharpAsymptotics}. Section \ref{corollaries} concludes the paper with further interesting corollaries of Theorem \ref{SharpAsymptotics}. Finally, the appendix discusses generalizations to other Markov processes and additive functionals. 

\section{Preliminaries} \label{sec:prelim}

This section is divided into three subsections: Subsection \ref{subsec: basic} recalls useful definitions and results from spectral theory of unbounded operators; Subsection \ref{subsec: spectral} presents the spectral theory framework for our proof of Theorem \ref{SharpAsymptotics}; and Subsection \ref{subsec: SDE trans} examines the SDE \eqref{newsde} and the transcendental equation \eqref{trans1}.
 
\subsection{Basic definitions and results} \label{subsec: basic}

Let $(H,\langle \cdot, \cdot \rangle)$ be a complex Hilbert space and let $A$ be a (possibly unbounded) closed, densely defined linear operator with domain $\mcal{D}(A)$. We sometimes view $\mcal{D}(A)$ as a Hilbert space equipped with the graph inner product $\langle u , v \rangle_{A} := \langle u , v \rangle + \langle Au, Av \rangle$. The resolvent set of $A$, denoted by $\rho(A)$, is defined as the set of $\xi \in \mathbb{C}$ such that $(\xi-A):\,\mcal{D}(A) \to H$ is invertible. For $\xi \in \rho(A)$, the Closed Graph Theorem (recall that $A$ is closed) implies that $(\xi-A)$ has a bounded inverse $(\xi - A)^{-1}:\, H \rightarrow \mcal{D}(A)$, called the resolvent (note that we adopt a sign convention opposite to \cite{Kato1995-kx}). The \textit{spectrum} $\sigma(A)$ of $A$ is defined by $\sigma(A) = \mathbb{C}\setminus \rho(A)$. 

An operator $T$ is semi-Fredholm if its range is closed, and either the dimension of the kernel of $T$ is finite or the codimension of the range of $T$ is finite. The semi-Fredholm domain $\Delta(A)\supseteq \rho(A)$ of $A$ is defined as the set of $\xi \in \mathbb{C}$ such that $(\xi-A):\,\mcal{D}(A) \to H$ is semi-Fredholm. The \textit{essential spectrum} of $A$ is then $\sigma_{\text{ess}}(A) := \mathbb{C}\setminus\Delta(A)\subseteq \sigma(A)$. 

We call $\alpha \in \mathbb{C}$ an \textit{eigenvalue} if $A \phi = \alpha \phi$ for some $\phi \in \mcal{D}(A)\setminus\{0\}$. In the case of an isolated eigenvalue, the following definition makes sense. 

\begin{defn}
Let $A$ be a (possibly unbounded) closed, densely defined linear operator and let $\alpha \in \sigma(A)$ be an isolated (topologically from the rest of the spectrum) eigenvalue. Then, the {\em Riesz projection} is defined by 
	\begin{equation}\label{RieszDef}
		P = \frac{1}{2\pi i}\int_{\Gamma}(\xi-A)^{-1}\,\mathrm{d}\xi,
	\end{equation}
where $\Gamma$ is a (positively oriented) contour in $\mathbb{C}$ that separates $\alpha$ from $\sigma(A)\setminus\{\alpha\}$. The {\em algebraic eigenspace} of $\alpha$ is defined as the range of $P$, and the {\em algebraic multiplicity} as its dimension. An eigenvalue is called {\em algebraically simple} if its algebraic multiplicity is $1$.
\end{defn}

The Riesz projection can be defined more generally whenever the spectrum of an operator can be partitioned by a contour, but we will only need the specific form above (see \cite[Theorem III.6.17 and its proof]{Kato1995-kx}).

\medskip

We also consider strongly continuous semigroups. A particularly nice class of strongly continuous semigroups is that of \textit{analytic semigroups}. For $\delta \in \big(0,\frac\pi2\big]$, define the sector $\Theta_{\delta}= \big\{\xi \in \mathbb{C} :\, |\arg(\xi)| < \delta \big\}\setminus\{0\}$. For convenience, we recall the definition of an analytic semigroup from \cite[Definition II.4.5]{Engel1999-tt}.

\begin{defn}
	A family of bounded operators $(T_z)_{z \in \Theta_{\delta}\cup\{0\}}$ is called an {\em analytic semigroup} (of angle $\delta \in \big(0,\frac\pi 2\big]$) if 
	\begin{enumerate}[(i)]
		\item $T_0=I$ and $T_{z_1+z_2}=T_{z_1}T_{z_2}$ for all $z_1,z_2 \in \Theta_{\delta}$.
		\item The map $z \mapsto T_z$ is analytic in $\Theta_\delta$.
		\item For all $u \in H$ and all $0<\delta'<\delta$, we have that $\lim_{\Theta_{\delta'}\ni z \rightarrow 0}T_z u=u$.     
	\end{enumerate}
	If we also have that $\norm{T_z}$ is bounded in $\Theta_{\delta'}$ for every $0 < \delta' < \delta$, then $(T_z)_{z \in \Theta_{\delta}\cup\{0\}}$ is called a {\em bounded analytic semigroup}.
\end{defn}

The following definition, which plays a central role in our perturbative analysis, is adapted from \cite[Sections VII.2, VII.3]{Kato1995-kx}.

\begin{defn} \label{def: type A}
Given an open domain $D_0\subset{\mathbb C}$, we call a family $\{T_z\}_{z \in D_0}$ of (possibly unbounded) closed, densely defined linear operators on $H$ a {\em holomorphic family of type (A)} if $\mcal{D}(T_z)=\mcal{D}$ is independent of $z$, and for every $u \in \mcal{D}$, the map $z \mapsto T_z u$ is holomorphic on $D_0$. The family is referred to as {\em self-adjoint holomorphic of type (A)} if additionally, for every $z\in D_0$ with $\overline{z}\in D_0$, the adjoint of $T_z$ is given by $T_{\overline{z}}$. 
\end{defn}

If $A$ is a (possibly unbounded) closed, densely defined, self-adjoint linear operator on $H$, and $V$ is a bounded, self-adjoint linear operator on $H$, it is easy to see from Definition \ref{def: type A} that $T_z = A+zV$ defines a self-adjoint holomorphic family of type (A) on ${\mathbb C}$. 

\medskip

Next, we present three general results that we will use repeatedly. The first is a (very specialized) version of Weyl's Theorem.

\begin{prop}\label{WeylTheorem} 
Let $A$ be a (possibly unbounded) closed, densely defined linear operator on $H$. Let $B$ be a bounded linear operator on $H$ that is $A$-compact, i.e., $B:\, (\mcal{D}(A),\norm{\cdot}_{A}) \rightarrow H$ is a compact operator. Suppose that $\Delta(A)$ is connected and that $\rho(A+B)$ is non-empty. Then, $\sigma_{\text{ess}}(A+B)=\sigma_{\text{ess}}(A)$, and any $\alpha \in \sigma(A+B)\setminus \sigma_{\text{ess}}(A)$ is an isolated eigenvalue of $A+B$ of finite algebraic multiplicity.
\end{prop}

\begin{proof}
By \cite[Theorem IV.5.35]{Kato1995-kx}, $\sigma_{\text{ess}}(A+B)=\sigma_{\text{ess}}(A)$. Per the discussion on \cite[p.~243]{Kato1995-kx}, since the connected set $\Delta(A+B)$ contains a resolvent point, any $\alpha \in \Delta(A+B)\setminus{\rho(A+B)}$ is an isolated eigenvalue of $A+B$ of finite algebraic multiplicity.
\end{proof}

The following result is a well-known property of the resolvent of a self-adjoint operator.

\begin{lemma} \label{specdist}
Let $A$ be a (possibly unbounded) closed, densely defined, self-adjoint linear operator. Then, for all $\xi \in \rho(A)$, we have
\begin{equation} \label{eq: norm}
\norm{(\xi-A)^{-1}}=\frac{1}{\mathrm{dist}(\xi,\sigma(A))}.
\end{equation}
\end{lemma}

\begin{proof}
The spectral theorem in the form of \cite[Theorem I.4.9]{Engel1999-tt} implies that $A$ is unitarily equivalent to multiplication by a real-valued function $q$ on some $L^2(\Omega,\mu)$-space, and so it suffices to prove \eqref{eq: norm} for multiplication operators. Let us write $M_f$ for the multiplication operator associated with a measurable function $f$. By \cite[Proposition I.4.10]{Engel1999-tt}, the spectrum of $M_q$ is equal to the essential range of $q$. Moreover, for any $\xi \in \mathbb{C}$, $\xi - M_q = M_{\xi - q}$. Finally, by \cite[Proposition I.4.10]{Engel1999-tt}, $M_{\xi-q}$ has a bounded inverse if and only if $\xi \notin \sigma(M_q)$, and in this case, $\norm{M_{\xi-q}^{-1}} = \sup_{\alpha \in \sigma(M_q)} \frac{1}{|\xi-\alpha|}$, which proves the lemma.
\end{proof}

The last result of this subsection provides a locally uniform bound on the resolvents of holomorphic families of type (A).

\begin{lemma}\label{CompactResolvent}
Let $(T_z)_{z \in D_0}$ be a holomorphic family of type (A). Fix some $z_0 \in D_0$, and suppose that there exists a compact $K \subset \rho(T_{z_0})$. Set $R(\xi;z) = (\xi - T_z)^{-1}$ for $\xi\in\rho(T_z)$. Then, there exists a neighborhood $U\subset D_0$ of $z_0$ such that $K \subset \rho(T_z)$ for all $z \in U$, and
	\begin{equation}
	\sup_{z \in U} \,\sup_{\xi \in K}\, \norm{R(\xi;z)} < \infty.
	\end{equation}
	Furthermore, the statement remains true if we replace the operator norm of $R(\xi;z)$ by its norm as a map from the ambient Hilbert space $H$ to the domain of $T_z$.
\end{lemma}

\begin{proof}
By \cite[Theorem VII.1.3]{Kato1995-kx}, for each $\xi_0\in K$, there exists a neighborhood $D_{\xi_0}\times U_{\xi_0}$ of $(\xi_0,z_0)$ on which $(\xi - T_z)^{-1}$ is jointly bounded-holomorphic. The open cover $\{D_\xi\}_{\xi \in K}$ of the compact $K$ admits a finite subcover $(D_{\xi_i})_{i=1}^n$. Setting $D = \bigcup_{i=1}^n D_{\xi_i}$ and $U = \bigcap_{i=1}^n U_{\xi_i}$, we find that $(\xi - T_z)^{-1}$ is jointly holomorphic on $D \times U$. Therefore, shrinking $U$ if necessary, we arrive at
	\begin{equation} \label{lemma first part}
	\sup_{z \in U} \,\sup_{\xi \in K}\, \norm{(\xi - T_z)^{-1}}<\infty.
	\end{equation}
	
To extend \eqref{lemma first part} to the stronger norm $\norm{\cdot}_{\mcal{L}(H,\mcal{D}(T_z))}$, we note that if $u = (\xi-T_z)^{-1} v$, then $T_z u = \xi u - v$. Hence,
	\begin{equation}
	\norm{T_z u} \leq |\xi|\norm{(\xi - T_z)^{-1}v} + \norm{v},
	\end{equation}
which can be estimated in terms of $\norm{(\xi-T_z)^{-1}}$ and then controlled by \eqref{lemma first part}.
\end{proof}

\subsection{Spectral theory framework} \label{subsec: spectral}

The goal of this subsection is to prove Theorems \ref{DriftGeneration} and \ref{DriftSpectrum}, which study the Schrödinger operators $T^\mu_z:=\frac{1}{2}\partial_{xx}+\mu\partial_x+zV$ for $\mu\ge0$, $z\in{\mathbb C}$, and $V:=\mathbf{1}_{[-1,1]}$. Theorem \ref{DriftSpectrum} is the main spectral input needed for our proof of Theorem \ref{SharpAsymptotics}, whereas Theorem \ref{DriftGeneration} establishes important regularity properties of the semigroups generated by $T^\mu_z$. We let $L^2_{\mu} := L^2(\mathbb{R}\to\mathbb{C},e^{2 \mu x}\mathrm{d}x)$ and write $H^k_{\mu}$ for the space of functions $f \in L^2_{\mu}$ whose distributional derivatives up to order $k$ lie in $L^2_{\mu}$. When $\mu=0$, we omit the subscript. Moreover, we denote by $X$ the canonical stochastic process, by $\mathbb{P}_x^{\mu}$ a probability measure under which $X_t = x + \mu t + B_t$, $t\ge0$ (where $B$ is a standard Brownian motion), and by $\mathbb{E}_x^{\mu}$ the corresponding expectation.

\begin{thm} \label{DriftGeneration}
	\begin{enumerate}[(a)]
		\item For every $z\in\mathbb{C}$, the operator $T_z^{\mu}$ on $L^2_{\mu}$ with domain $\mcal{D}(T_z^\mu) = H^2_{\mu}$ is closed and densely defined. It generates an analytic semigroup on $L^2_{\mu}$.
		\item For every $f \in L^2_{\mu}$, there exists a $C^{1,1}$ and piecewise $C^{1,2}$ function $G$ on $(0,\infty)\times\mathbb{R}$ such that $e^{t T_z^{\mu}}f(x) = G(t,x)$ for almost every $(t,x)\in(0,\infty)\times\mathbb{R}$. We henceforth use $e^{tT_z^\mu}f(x)$ for this continuous representative.
		\item For all $f \in L^2_{\mu}$, $t > 0$, and $x\in\mathbb{R}$, we have the Feynman-Kac formula
		\begin{equation} \label{FK repr}
		G(t,x) = \mathbb{E}^{\mu}_x\bigg[f(X_t)\,\exp\bigg(z \int_0^t V(X_s)\,\mathrm{d}s \bigg)\bigg].
		\end{equation}
	\end{enumerate}
\end{thm}

\begin{thm} \label{DriftSpectrum}
For $z\in\mathbb{R}$, the operator $T^\mu_z$ is self-adjoint; thus, $\sigma(T_z^\mu)\subset\mathbb{R}$. For $\lambda>0$,
\begin{equation} \label{eq:spectral bound}
s(T^\mu_\lambda) := \sup_{\alpha\in\sigma(T^\mu_\lambda)} \mathrm{Re}(\alpha)
= \Lambda(\lambda) - \frac{\mu^2}{2} =: \Lambda_{\mu}(\lambda)
\end{equation}
is an algebraically simple eigenvalue of $T_{\lambda}^{\mu}$ corresponding to the strictly positive eigenvector $e^{- \mu x}\phi_\lambda$ (see equation \eqref{Eigenvector}). Furthermore, there is a neighborhood $U\subset\mathbb C$ of $\lambda$ (independent of $\mu$) such that: 
	\begin{enumerate}[(a)]
		\item For $z\in U$, the operator $T_z^{\mu}$ has an algebraically simple eigenvalue with real part $s(T_z^{\mu})$, which we denote $\Lambda_{\mu}(z)$. We have the relationship $\Lambda_{\mu}(z) = \Lambda_0(z) - \frac{\mu^2}{2}$. Furthermore, the map $z \mapsto \Lambda_{\mu}(z)$ is holomorphic on $U$.
		\item Let $\{P^{\mu}_z\}_{z \in U}:\, L^2_{\mu} \rightarrow H^2_{\mu}$ denote the Riesz projections for $\Lambda_{\mu}(z)$ and define $Q^{\mu}_z = (I-P^{\mu}_z):\, L^2_{\mu} \rightarrow L^2_{\mu}$. Then, for all $t\ge0$, $z\in U$, and $f \in L^2_{\mu}$,
		\begin{equation} \label{eq:decomp}
		e^{tT^{\mu}_z}f = e^{t \Lambda_{\mu}(z)}P^{\mu}_z f + e^{t T^{\mu}_z}Q^{\mu}_zf.
		\end{equation}
		\item For all $f \in L^2_{\mu}$, one can take $P^{\mu}_zf(x)$ as a continuous function in $x$ and a holomorphic function in $z \in U$ satisfying $P^{\mu}_zf(x) \rightarrow \langle e^{\mu \cdot}f, \phi_\lambda\rangle_{L^2}\, e^{-\mu x}\phi_\lambda(x)$ locally uniformly in $x$ as $z \rightarrow \lambda$.
		\item For any compact $K\subset \mathbb{R}$, there exist $\varepsilon > 0$ and $C=C(K)<\infty$ such that for all $t\ge1$, $z\in U$, $f \in L^2_{\mu}$, and $x \in K$,
		\begin{equation}
			\big|e^{t T^{\mu}_z}Q^{\mu}_zf(x)\big|
			\leq C\norm{f}_{L^2_{\mu}}e^{t(\Lambda_{\mu}(\lambda)-\varepsilon)}.
		\end{equation}
        Note that due to Theorem \ref{DriftGeneration}(b) and equation \eqref{eq:decomp}, we may also take $e^{t T^{\mu}_z}Q^{\mu}_zf(x)$ to be continuous in $(t,x)$.
		\item For any compact $K\subset \mathbb{R}$ and $0 < m < M <\infty$, there are $\varepsilon \!>\! 0$ and $C \!=\! C(K,m,M) \!<\! \infty$ such that for all $t\ge1$, $f \in L^2_{\mu}$, and $x\in K$, 
		\begin{equation*}
			\sup_{y\in\mathbb{R}:\, m \leq |y| \leq M} \big|e^{t T_{\lambda+iy}^{\mu}}f(x)\big| 
			\leq C \norm{f}_{L^2_{\mu}}e^{t(\Lambda_{\mu}(\lambda)-\varepsilon)}.
		\end{equation*}
	\end{enumerate}
\end{thm}

We will see that Theorems \ref{DriftGeneration} and \ref{DriftSpectrum} for $\mu>0$ can be readily deduced from their $\mu=0$ counterparts. Therefore, we first focus on the $\mu=0$ case, beginning with the proofs of Theorem \ref{DriftGeneration} and the first two statements of Theorem \ref{DriftSpectrum}.

\begin{prop} \label{MainGeneration}
Theorem \ref{DriftGeneration} holds for $\mu=0$. For $z \in \mathbb{R}$, the operator $T^0_z$ is self-adjoint; thus, $\sigma(T^0_z)\subset\mathbb{R}$. For $\lambda>0$, the quantity $s(T^0_\lambda) = \Lambda(\lambda)$ is an algebraically simple eigenvalue of $T^0_\lambda$ corresponding to the strictly positive eigenvector $\phi_\lambda$ (see equation \eqref{Eigenvector}).
\end{prop}

\begin{proof}
	\textit{Step 1.} First, note that the operator $T^0_0$ is unitarily equivalent, via the Fourier transform, to the multiplication operator $M$ defined on $L^2(\mathbb{R}\to\mathbb{C},\mathrm{d}k)$ by $Mf(k) = -\frac{k^2}{2}f(k)$. The domain of this multiplication operator is the image of $H^2$ under the Fourier transform; furthermore, the operator is closed, self-adjoint, and its spectrum is given by $(-\infty,0]$ (see \cite[Proposition I.4.10]{Engel1999-tt}). Therefore, ${\mcal D}(T^0_0)=H^2\underset{\text{dense}}{\subset} L^2$; the operator $T^0_0$ is closed, self-adjoint, and $\sigma(T^0_0)=(-\infty,0]$. Using Lemma \ref{specdist}, one checks directly that $T^0_0$ is a sectorial operator (of angle $\frac\pi2$) in the sense of \cite[Definition II.4.1]{Engel1999-tt}, and so it generates a bounded analytic semigroup by \cite[Theorem II.4.6]{Engel1999-tt}. The perturbation $T^0_z-T^0_0=z V$ is readily checked to be $T^0_0$-bounded with $T^0_0$-constant $0$ (see \cite[Definition III.2.1]{Engel1999-tt}). Hence, by \cite[Theorem III.2.10]{Engel1999-tt}, $T^0_z$ generates an analytic semigroup. Moreover, ${\mcal D}(T^0_z)=H^2$ and the closedness of $T^0_z$ are inherited from the same properties of $T^0_0$. For $z\in\mathbb{R}$, $zV$ is a real-valued multiplication operator, thus self-adjoint, yielding the self-adjointness of $T^0_z$.
	
\medskip
	
	\noindent\textit{Step 2.} We remark in passing that the graph norm on $\mcal{D}(T^0_z)$ is equivalent to the $H^2$-norm. Indeed, using integration by parts and the Cauchy-Schwarz inequality, we find that $\norm{u}_{H^2}\le C(\norm{u}_{L^2} + \norm{\partial_{xx} u}_{L^2})$ for some $C<\infty$. Moreover, it is easy to see that $\norm{\partial_{xx}u}_{L^2}\leq C(\norm{T^0_{z} u}_{L^2} + \norm{u}_{L^2})$ for some $C<\infty$. Combining these yields $\norm{u}_{H^2}\leq C \norm{u}_{\mcal{D}(T^0_{z})}$ for some $C<\infty$. The reverse inequality is obvious, so we have shown the equivalence.
	
\medskip	
	
	Next, we argue Theorem \ref{DriftGeneration}(b). Recall that by the Sobolev embedding theorem (see, e.g., \cite[Proposition IV.1.3, Corollary IV.1.4]{Taylor2023-zs}), there is a bounded linear map $J:\, H^2 \rightarrow C^1_b$. We claim that $G(t,x) := J(e^{t T^0_z}f)(x)$ has the desired properties. Indeed, since $T^0_z$ generates an analytic semigroup, the latter is immediately infinitely norm differentiable (see \cite[Definition II.4.13]{Engel1999-tt} and the discussion afterwards), and so $e^{t T^0_z}L^2 \subseteq \mcal{D}((T_z^0)^n)$ for all $t>0$ and $n\in\mathbb{N}$. Consequently, $t\mapsto e^{t T^0_z}f$ is continuously differentiable on $(0,\infty)$ with respect to the graph norm on $\mcal{D}(T^0_z)$, and $\partial_t e^{t T^0_z}f = T^0_z e^{t T^0_z}f$. Since the graph norm on $\mcal{D}(T^0_z)$ is equivalent to the $H^2$-norm, and due to the definition of $J$, we have that $G \in C^1((0,\infty);C^1_b)$. Hence, $G$, $\partial_tG$, and $\partial_x G$ exist as jointly continuous functions. Furthermore, for all $t>0$, the second distributional derivative in $x$ of the function $G(t,\cdot)$ is given by 
	\begin{equation}
	\frac12 \partial_{xx} G = \partial_t G - zV G,
	\end{equation}
	which can be chosen as a continuous function on each of $(0,\infty)\times(-\infty,-1)$, $(0,\infty)\times[-1,1]$, and $(0,\infty)\times(1,\infty)$. Thus, $G$ is piecewise $C^{1,2}$.
	
\medskip
	
	\noindent\textit{Step 3.} We turn to the proof of Theorem \ref{DriftGeneration}(c). Fix a $z \in \mathbb{C}$, and define the associated Feynman-Kac semigroup $(S_t)_{t\ge0}$ on $L^2$:
	\begin{equation}\label{FKsemigroup}
		(S_tf)(x) := \mathbb{E}^0_x\bigg[f(X_t)\,
		\exp\bigg(z \int_0^t V(X_s)\,\mathrm{d}s\bigg)\bigg].
	\end{equation}
	Note the inequalities $|S_t f| \leq e^{t|z|}e^{t T^0_0}|f|$ and $\big\|S_t f - e^{t T^0_0} f\big\|_{L^2} \le (e^{t|z|}-1)\norm{f}_{L^2}$, which imply that $(S_t)_{t\ge0}$ is a strongly continuous semigroup on $L^2$. Now, using the Fundamental Theorem of Calculus identity $e^{z\int_0^t V(X_s)\,\mathrm{d}s}=1+z\int_0^t V(X_s)\, e^{z\int_0^s V(X_r)\,\mathrm{d}r}\,\mathrm{d}s$, we find that
	\begin{equation}
	S_t f = e^{t T^0_0} f + z\int_0^t S_s(Ve^{(t-s)T^0_0}f)\,\mathrm{d}s.
	\end{equation}
	Let $\mcal{G}$ be the generator of $(S_t)_{t \geq 0}$. It is then straightforward to see that $\mcal{D}(T^0_0)\subseteq \mcal{D}(\mcal{G})$ and that $\mcal{G}=T^0_z$ on $\mcal{D}(T^0_0)$. Thus, we aim to show that $\mcal{D}(T^0_0)=\mcal{D}(\mcal{G})$.
	
\medskip
	
	As a preliminary step, we argue that $\xi \in \rho(T^0_z)\cap\rho(\mcal{G})$ for all $\xi \in \mathbb{R}$ large enough. Indeed, the contraction property of $(e^{tT^0_0})_{t\ge0}$ implies that the operator norm $\|S_t\|$ of $S_t$ is bounded above by $e^{t|z|}$ for each $t\ge0$. Hence, by \cite[Theorem II.1.10]{Engel1999-tt}, it follows that $\xi\in\rho(\mcal{G})$ for $\xi > |z|$. Also for $\xi >|z|$, we have $\xi \in \rho(T^0_0)$ and, by Lemma \ref{specdist}, $\norm{zV(\xi-T^0_0)^{-1}}\le\frac{|z|}{\xi} <1$. Thus, $\xi-T^0_z = (I-zV(\xi-T^0_0)^{-1})(\xi-T^0_0)$ is invertible, which means that $\xi \in \rho(T^0_z)$.
	
\medskip
	
	We can now show that $\mcal{D}(T^0_0)=\mcal{D}(\mcal{G})$. To this end, we let $f \in \mcal{D}(\mcal{G})$ and choose $\xi \in \rho(T^0_z)\cap\rho(\mcal{G})\neq \varnothing$.  Set $g = (\xi-\mcal{G})f$ and define $h\in\mcal{D}(T^0_0)$ by $(\xi-T^0_z)h = g$. Since $\mcal{D}(T^0_0)\subseteq \mcal{D}(\mcal{G})$ and $\mcal{G}=T^0_z$ on $\mcal{D}(T^0_0)$, we have that $h\in\mcal{D}(\mcal{G})$ and $(\xi-\mcal{G})h = (\xi-T^0_z)h = g = (\xi-\mcal{G})f$. Moreover, $\xi \in \rho(\mcal{G})$, so that $f=h \in \mcal{D}(T^0_0)$, which proves that $\mcal{D}(T^0_0) = \mcal{D}(\mcal{G})$. Thus, 
	$T^0_z = \mcal{G}$, and the semigroup $(e^{tT^0_z})_{t\ge0}$ must be given by \eqref{FKsemigroup}.
	
	Strictly speaking, this yields \eqref{FK repr} for almost every $x\in\mathbb{R}$. To obtain \eqref{FK repr} for all $x\in\mathbb{R}$, let $G_\varepsilon$ denote the right-hand side of \eqref{FKsemigroup} with the lower bound of the integral replaced by $\varepsilon$. From the Markov property of $X$, we see that $G_\varepsilon(t,\cdot)$ is a continuous function for any fixed $t\ge\varepsilon$. It is also readily checked that $|S_tf(x)-G_\varepsilon(t,x)| \leq e^{|z|t}(e^{|z|\varepsilon}-1)\mathbb{E}^0_x[|f(X_t)|]$. By Young's convolution inequality,  $\sup_{x\in\mathbb R} \mathbb{E}^0_x[|f(X_t)|]\leq C\norm{f}_{L^2}$ for some $C<\infty$. Therefore, $x\mapsto S_tf(x)$ is a uniform limit of continuous functions, and so must be a continuous function. Since the equality $S_tf(x)=e^{tT^0_z}f(x)$ holds for almost every $x\in\mathbb{R}$, and both sides are continuous in $x$, it must hold for all $x\in\mathbb{R}$.
	
\medskip
		
	\noindent\textit{Step 4.} It remains to show the last statement in Proposition \ref{MainGeneration}. To this end, we first note that $\Lambda(\lambda)$ is an eigenvalue of $T^0_\lambda$. Indeed, a direct computation based on equation \eqref{Eigenvector} reveals that $\phi_\lambda$ is continuously differentiable, belongs to $H^2$, and is an eigenfunction of $T^0_\lambda$ with eigenvalue $\Lambda(\lambda)>0$. 
	
\medskip

	Clearly, $s(T^0_\lambda) \ge \Lambda(\lambda)$. Suppose for a contradiction that $s(T^0_\lambda) > \Lambda(\lambda)$. We claim that the operator $T^0_{\lambda}-T^0_0 = \lambda V$ is $T^0_0$-compact. To see this, we observe that any $K\subset\mcal{D}(T^0_0)$ bounded in the graph norm is also bounded in $H^2$. Therefore, the respective restrictions of functions in $K$ to $[-1,1]$ are uniformly bounded in the $H^2([-1,1]\to\mathbb{C},\mathrm{d}x)$-norm, and so by the Rellich-Kondrachov Theorem (see \cite[Theorem 5.7.1]{evans2010partial}) admit a sequence $(u_n)_{n\in\mathbb{N}}$ converging strongly to some $u$ in $L^2([-1,1]\to\mathbb{C},\mathrm{d}x)$. Hence, $(\lambda V u_n)_{n\in\mathbb{N}}$ converges strongly in $L^2$ to $\lambda V u$. Moreover, $\mcal{D}(T^0_0)=H^2\subset L^2= \mcal{D}(\lambda V)$, which renders $\lambda V$ a $T^0_0$-compact operator. Thus, $T^0_\lambda$ has the same essential spectrum as $T^0_0$ by Proposition \ref{WeylTheorem}. 
	
	\medskip
	
	Next, we argue that $s(T^0_{\lambda}) < \lambda$. 
	By Proposition \ref{WeylTheorem}, any $\alpha\in\sigma(T^0_\lambda)$ with $\alpha>\Lambda(\lambda)$ must be an eigenvalue of finite algebraic multiplicity. Let $\phi$ be a real eigenfunction with norm $1$ corresponding to such an $\alpha$. Then,
	\begin{equation} \label{eq: less than lambda}
	\alpha = \langle T^0_\lambda \phi, \phi \rangle_{L^2} = - \frac{1}{2}\int_{-\infty}^\infty (\partial_x \phi)^2 \,\mathrm{d}x + \lambda \int_{-1}^1 \phi^2\,\mathrm{d}x \leq \lambda.
	\end{equation}
	The latter inequality must be strict, since otherwise we could infer $\partial_x\phi=0$ almost everywhere on $\mathbb{R}$, $\phi=0$ on $\mathbb{R}\setminus[-1,1]$, and $\phi=0$ on $\mathbb{R}$ due to $\phi \in H^1$. Consequently, $s(T^0_\lambda) < \infty$. Since $\sigma(T^0_\lambda)$ is closed, we have that $s(T^0_{\lambda}) \in \sigma(T^0_\lambda)$, but is not a member of the essential spectrum, and so $s(T^0_\lambda)$ must be an eigenvalue of finite algebraic multiplicity. Therefore, $s(T^0_{\lambda}) < \lambda$ by the above argument. 
	
	\medskip
	
	Furthermore, for each $t>0$, the operator $e^{t T^0_\lambda}$ is positivity-improving (see, e.g., \cite[definition on p.~201]{ReedSimonIV}) by the Feynman-Kac formula \eqref{FK repr}. Thus, in view of \cite[Theorem XIII.44]{ReedSimonIV}, $s(T^0_\lambda)$ is a geometrically simple eigenvalue of $T^0_\lambda$ corresponding to a strictly positive eigenfunction $\phi\in\mathcal{D}(T^0_\lambda)=H^2\hookrightarrow C_b^1$. On $(1,\infty)$, the eigenfunction $\phi$ satisfies $\frac{1}{2}\partial_{xx}\phi = s(T^0_\lambda)\phi\in H^2((1,\infty))$. Hence, $\phi \in H^4((1,\infty))$, and upon iteration, $\phi\in C^\infty((1,\infty))$. Similarly, $\phi\in C^\infty((-1,1))$ and $\phi\in C^\infty((-\infty,-1))$. Moreover, $\phi$ must be even. Indeed, $x\mapsto\phi(-x)$ is a positive eigenfunction with the same norm, and so by simplicity equal to $\phi$. Combining these properties, we find that $\phi$ must be a (positive) constant multiple of 
	\[
	\begin{cases}
		\widetilde{M}e^{-\sqrt{2s(T^0_\lambda)}\,|x|} \hspace{36.5pt}\text{ if } |x|\geq 1, \\
		\cos\big(\sqrt{2(\lambda - s(T^0_{\lambda}))}\, x\big) \text{ if } |x|\leq 1
	\end{cases}
	\]
	for some $\widetilde{M}\in(0,\infty)$. Since $\phi$ and $\partial_x\phi$ are continuous at $1$, we infer that $s(T^0_\lambda)$ solves \eqref{trans1}. This is a contradiction with $s(T_\lambda^0)>\Lambda(\lambda)$ since $\Lambda(\lambda)$ was the only solution to \eqref{trans1} for which the argument of the tangent function therein lies in $\big(0,\frac{\pi}{2}\big)$.
	
\medskip
	
	Lastly, we turn to the algebraic simplicity of $\Lambda(\lambda)$, i.e., the one-dimensionality of the range of the Riesz projection $P^0_\lambda$ (see equation \eqref{RieszDef}), which is the algebraic eigenspace of $\Lambda(\lambda)$. To this end, let $D := (T^0_\lambda - \Lambda(\lambda))P^0_\lambda$ be the associated eigennilpotent, which is bounded due to \cite[proof of Theorem III.6.17, display (6.24)]{Kato1995-kx}. Since $T^0_\lambda$ is a self-adjoint operator, so is $P^0_\lambda$.  Moreover, $T^0_\lambda$ commutes with $P^0_\lambda$ on ${\mathcal D}(T^0_\lambda)$. This implies the self-adjointness of $D$. As a self-adjoint nilpotent operator, $D$ must be identically zero. That is, the range of $P^0_\lambda$ is contained in the kernel of $T^0_\lambda - \Lambda(\lambda)$, and the algebraic simplicity of $\Lambda(\lambda)$ follows from its geometric simplicity. 
\end{proof}

We now complete the proof of Theorem \ref{DriftSpectrum} in the $\mu=0$ case. 

\begin{prop} \label{MainPerturbation}
Let $\mu=0$ and $\lambda>0$. Then, there exists a neighborhood $U\subset\mathbb{C}$ of $\lambda$ such that the statements (a)--(e) in Theorem  \ref{DriftSpectrum} hold. 
\end{prop}

Due to its length, we break the proof of Proposition \ref{MainPerturbation} into a series of lemmas. 

\begin{lemma} \label{topeig}
Let $\mu=0$ and $\lambda>0$. Then, there is a neighborhood $U\subset\mathbb{C}$ of $\lambda$ such that: 
	\begin{enumerate}[(a)]
		\item $\{T^0_z\}_{z\in U}$ is a self-adjoint holomorphic family of type (A) (see Definition \ref{def: type A}).
		\item There exists a ball $D_1\subset\mathbb{C}$ around $\Lambda(\lambda)$ bounded away from $0$ for which $\sigma(T^0_z)\cap \overline{D_1}$ is a singleton for all $z\in U$. We write $\{\Lambda(z)\}=\sigma(T^0_z)\cap \overline{D_1}$.
		\item $\mathrm{Re}(\Lambda(z))=s(T^0_z)$, $\Lambda(z)$ is an algebraically simple eigenvalue of $T^0_z$, and the function $z\mapsto\Lambda(z)$ is holomorphic on $U$.
		\item For all $f \in L^2$, one can choose the Riesz projection $P^0_zf(x)$ to be continuous in $x\in\mathbb{R}$ and holomorphic in $z\in U$. Moreover, $P^0_zf(x)\rightarrow P^0_\lambda f(x) = \langle f, \phi_\lambda \rangle_{L^2}\,\phi_\lambda(x)$ uniformly in $x$ as $z \rightarrow \lambda$.
	\end{enumerate}
\end{lemma}

\begin{proof}
	It follows directly from Definition \ref{def: type A} that $\{T^0_z\}_{z \in \mathbb{C}}$ form a self-adjoint holomorphic family of type (A). By Proposition \ref{MainGeneration}, $\Lambda(\lambda)$ belongs to the discrete spectrum of $T^0_\lambda$. Therefore, we can take an open ball $D_1$ of a radius $\delta_1\in(0,1)$ around $\Lambda(\lambda)$, bounded away from $0$ and such that the rest of $\sigma(T^0_\lambda)$ lies in $\mathbb{C}\setminus\overline{D_1}$. By \cite[Theorem VII.1.7]{Kato1995-kx} (see \cite[Theorem III.6.17]{Kato1995-kx} for the relevant definitions), for $z$ in a neighborhood $U\subset\mathbb{C}$ of $\lambda$, the spectrum of $T^0_z$ is also separated by $\partial D_1$ into two disjoint parts $\sigma'(z)$ and $\sigma''(z)$, with $\sigma'(z)\subset D_1$. Moreover, the projection associated with $\sigma'(z)$ is given by
	\begin{equation}
	P^0_z = \frac{1}{2\pi i}\int_{\partial D_1 }R(\zeta;z)\,\mathrm{d}\zeta:\, L^2 \rightarrow H^2,
	\end{equation}
	where $R(\zeta;z) := (\zeta - T^0_z)^{-1}$ is the resolvent (with the opposite sign convention from \cite{Kato1995-kx}). Thanks to Lemma \ref{CompactResolvent}, upon shrinking $U$ if necessary, each element of $\{P^0_z\}_{z\in U}$ is a bounded linear map from $L^2$ to $H^2$, and therefore from $L^2$ to $C^1_b$. The proof of Lemma \ref{CompactResolvent} implies further that $\{P^0_z\}_{z\in U}$ is bounded-holomorphic from $L^2$ to $C_b^1$. In particular, for each $f\in L^2$, we can choose $P^0_zf(x)$ to be continuous in $x\in\mathbb{R}$ and holomorphic in $z\in U$.    
	
	\medskip
	
	Furthermore, upon shrinking $U$ again if necessary, we find a holomorphic family $\{O_z\}_{z \in U}$ of bounded invertible linear maps such that $P^0_z = O_z P^0_\lambda O_z^{-1}$ (see the construction on \cite[p.~100]{Kato1995-kx} for the finite-dimensional case, which extends to the Banach space setting as noted in \cite[proof of Theorem VII.1.7]{Kato1995-kx}). Therefore, for all $z\in U$, the range of $P^0_z$ is one-dimensional. By \cite[Problem III.6.18]{Kato1995-kx}, the range of $P^0_z$ is spanned by a single eigenfunction whose eigenvalue $\Lambda(z)$ has algebraic multiplicity $1$. Moreover, \cite[Theorem III.6.17]{Kato1995-kx} shows that the spectrum of $T^0_z$ restricted to the range of $P^0_z$ is contained in $D_1$, so that $\Lambda(z)\in D_1$. If there were another element of $\sigma(T^0_z)\cap \overline{D_1}$, it would have to belong to the spectrum of $T^0_z$ restricted to the range of $I-P^0_z$, contradicting the separation of $\sigma'(z)$ and $\sigma''(z)$ by $\partial D_1$.
	
	\medskip
	
	Next, we turn to the proof of $\mathrm{Re}(\Lambda(z))=s(T^0_z)$. If necessary, we shrink $U$ further, to ensure that for all $z\in U$: $|\mathrm{Im}(z)|\leq1$, $0<\mathrm{Re}(z) \leq \lambda+1$, and the compact 
	\begin{equation}
	K := \big\{\alpha \in \mathbb{C}:\,\mathrm{Re}(\alpha) \in [\Lambda(\lambda)-\delta_1,\lambda+1],\,|\mathrm{Im}(\alpha)|\leq 1 \big\}
	\setminus D_1\subset\rho(T^0_z) 
	\end{equation}
	(recall from \eqref{eq: less than lambda} that $\Lambda(\lambda)<\lambda$, note that $K\subset\rho(T^0_\lambda)$, and apply Lemma \ref{CompactResolvent}). Suppose for a contradiction that there were an $\alpha \in \sigma(T^0_z)$ with $\mathrm{Re}(\alpha) > \mathrm{Re}(\Lambda(z))$. As in Step 4 of the proof of Proposition \ref{MainGeneration}, we infer that $\alpha$ is an eigenvalue of finite algebraic multiplicity. Since $\sigma(T^0_z)\cap \overline{D_1}=\{\Lambda(z)\}$, we have $\alpha\notin\overline{D_1}$. 
	Let $\phi$ be an eigenfunction with norm $1$ corresponding to $\alpha$. Then,
	\begin{equation}
	\alpha  = -\frac{1}{2}\int_{-\infty}^\infty |\phi '(x)|^2 \,\mathrm{d}x 
	+ z \int_{-1}^1 |\phi(x)|^2\,\mathrm{d}x,
	\end{equation}
	which yields $\mathrm{Re}(\alpha)\le\mathrm{Re}(z)\le\lambda+1$ and $|\mathrm{Im}(\alpha)| \le |\mathrm{Im}(z)|\le 1$. This means that $\alpha \in K$, contradicting $K\subset\rho(T^0_z)$. We note in passing that this argument implies the stronger statement $\mathrm{Re}(\alpha)<\Lambda(\lambda)-\delta_1$ for all $\alpha\in\sigma(T^0_z)\setminus\{\Lambda(z)\}$.
	
	\medskip
	
	Finally, as observed on \cite[p.~370]{Kato1995-kx}, the family $\{O_z^{-1}T^0_z P^0_z O_z\}_{z\in U}$, acting on the range of $P^0_\lambda$, is bounded-holomorphic. On the other hand, using $O_zP^0_\lambda=P^0_zO_z$, $(P^0_z)^2=P^0_z$, the fact that $T^0_z$ acts as multiplication by $\Lambda(z)$ on the range of $P^0_z$, and $O_z^{-1}P^0_zO_z=P^0_\lambda$, we find that $O_z^{-1}T^0_z P^0_z O_z$ acts as multiplication by $\Lambda(z)$ on the range of $P^0_\lambda$. Consequently, $z\mapsto\Lambda(z)$ is holomorphic on $U$. Lastly, since $\{P^0_z\}_{z\in U}$ is bounded-holomorphic from $L^2$ to $C^1_b$, we have $P^0_zf\to P_\lambda^0 f$ uniformly as $z\to\lambda$. Moreover, $P^0_\lambda$ is the orthogonal projection onto the span of $\phi_\lambda$, and therefore $P^0_\lambda f(x)= \langle f, \phi_\lambda \rangle_{L^2}\,\phi_\lambda(x)$.
\end{proof}

The next lemma is an important ingredient for the proof of Theorem \ref{DriftSpectrum}, statement (d).

\begin{lemma} \label{lowerorderresolvent}
	Let $\mu=0$ and $\lambda>0$. Define the following:
	\begin{itemize}
	\item $\delta_1\in(0,1)$ and the operator $P^0_z$ as in the first paragraph of the proof of Lemma \ref{topeig}; 
	\item $b=\Lambda(\lambda)-\delta_1$ and $A^0_z=T^0_z-b$; 
	\item the (closed) range $M''(z)$ of $Q^0_z:=I-P^0_z$; 
	\item the operator $\widetilde{T}^0_z$ as the restriction of $T^0_z$ to 
	$\mcal{D}(\widetilde{T}^0_z) := \mcal{D}(T^0_z)\cap M''(z)$ and $\widetilde{A}^0_z := \widetilde{T}^0_z - b$. 
	\end{itemize}
	Then, for any $\delta \in \big(0,\frac\pi2\big)$, there exist a constant $M<\infty$ and a neighborhood $U\subset\mathbb{C}$ of $\lambda$ such that for all $z\in U$, the sector
	\begin{equation} \label{sectordef}
		\Sigma_{\delta} := \Big\{\xi \in \mathbb{C}:\, |\arg(\xi)|<\frac\pi2 + \delta \Big\}\setminus\{0\}
	\end{equation}
	is contained in $\rho(\widetilde{A}^0_z)$, and 
	\begin{equation} \label{eq: res est}
	\big\|(\xi - \widetilde{A}^0_z)^{-1}\big\| \leq  \frac{M}{|\xi|},\quad (\xi,z)\in\Sigma_\delta\times U.
	\end{equation}
\end{lemma}

\begin{proof}
	Note that the conclusion of the lemma can be interpreted as a uniform version of the sectorial condition in \cite[Definition II.4.1]{Engel1999-tt}. We start the proof with the neighborhood $U\subset\mathbb{C}$ from Lemma \ref{topeig}, which we will shrink as necessary. By \cite[Theorem III.6.17]{Kato1995-kx}, $L^2 = M'(z) \oplus M''(z)$, where $M'(z)$ is the range of $P^0_z$, and the spectrum of $\widetilde{T}_z^0$ lies outside the closed ball $\overline{D_1}$ of radius $\delta_1$ around $\Lambda(\lambda)$.   
	
	\medskip
	
	We establish \eqref{eq: res est} by considering $\xi$'s in three distinct subregions of $\Sigma_\delta$. First, we consider a ball around $\Lambda(\lambda)-b=\delta_1$. From $\widetilde{A}^0_z = \widetilde{T}^0_z - b$, we have that
	$\sigma(\widetilde{A}^0_z)\cap (\overline{D_1}-b) = \varnothing$. Since the resolvent set $\rho(\widetilde{A}^0_z)$ is open, there exists a neighborhood of $\overline{D_1}-b$ contained in it. By \cite[Theorem III.6.7]{Kato1995-kx}, $(\xi - \widetilde{A}^0_z)^{-1}$ is holomorphic in $\xi$ on this neighborhood of $\overline{D_1}-b$. With $R_z''(\xi) := (\xi - \widetilde{A}^0_z)^{-1}(I-P^0_z)$ and an open ball $D_2$ of radius $\delta_2\in(0,\delta_1)$ around $\Lambda(\lambda)-b$, the Cauchy integral formula implies that
	\begin{equation}
	R_z''(\zeta)=\frac{1}{2 \pi i}\int_{\partial D_1-b}\frac{R_z''(\xi)}{\xi - \zeta}\,\mathrm{d}\xi,\quad \zeta\in D_2.
	\end{equation}
	It follows that
	\begin{equation} \label{Cauchy int bound}
	\norm{R_z''(\zeta)}\leq \frac{\delta_1}{\delta_1 - \delta_2}\sup_{\xi \in \partial D_1-b}\norm{R_z''(\xi)},\quad \zeta\in D_2.
	\end{equation}
	Moreover, $\{P^0_z\}_{z\in U}$ is bounded-holomorphic, and thus the same is true for $\{I-P^0_z\}_{z\in U}$. In addition, by Lemma \ref{CompactResolvent}, $\sup_{z \in U} \sup_{\xi \in \partial D_1-b} \norm{(\xi-A^0_z)^{-1}} < \infty$ upon shrinking $U$ if necessary. Combining these with $(\xi-A^0_z)^{-1}(I-P^0_z) 
	= (\xi - \widetilde{A}^0_z)^{-1}(I-P^0_z)$ and \eqref{Cauchy int bound}, we deduce $\sup_{z \in U} \sup_{\mu \in D_2} \|R''_z(\mu)\|<\infty$.\comment{Decomposing $(\xi-\widetilde{A}^0_z)^{-1}=(\xi-\widetilde{A}^0_z)^{-1}(I-P^0_z)+(\xi-\widetilde{A}^0_z)^{-1}P^0_z$ and recalling that $(\xi-\widetilde{A}^0_z)^{-1}$ acts as multiplication by $(\xi+b)^{-1}$ on the range of $P^0_z$, we arrive at $\sup_{z \in U} \sup_{\xi \in D_2} \|(\xi-\widetilde{A}^0_z)^{-1}\|<\infty$.} Since $Q_z^0 = I - P_z^0$ acts as the identity on $M''(z)$, the resolvent $(\xi - \widetilde{A}_z^0)^{-1}$ is the restriction of $R''_z(\xi)$ to $M''(z)$. Thus,
    \[
    \sup_{z \in U}\sup_{\xi \in D_2}\big|\hspace{-1pt}\big|(\xi - \widetilde{A}_z^0)^{-1}\big|\hspace{-1pt}\big|_{\mcal{L}(M''(z))}\leq \sup_{z \in U}\sup_{\xi \in D_2}\big|\hspace{-1pt}\big|{R''_z(\xi)}\big|\hspace{-1pt}\big|_{\mcal{L}(L^2)}<\infty.
    \]
	
	\medskip
	
	Next, we turn to a subregion of the form $\Sigma_{\delta,R} := \big\{\xi \in \mathbb{C}:\, |\arg(\xi)|<\frac\pi 2 + \delta,\, |\xi| \geq R \big\}$. With $C := \sup_{z\in U} |z|$ and  
	$R := \max \big\{2b, 4C,\frac{2C}{\cos\delta}\big\}$, we claim that $\Sigma_{\delta,R}\subset \rho(A^0_z)$, $z \in U$. Indeed, if $\mathrm{Re}(b + \xi)\geq 0$ and $|\xi| \geq 2b$, then
	\begin{equation} \label{outsidecompact}
		\big\|(b+\xi - T^0_0)^{-1}\big\|= \frac{1}{|\xi+b|}\leq \frac{2}{|\xi|}.
	\end{equation}
	If $\mathrm{Re}(b + \xi) <0$, then $\xi \in \Sigma_{\delta,R}$ and $\mathrm{Re}(\xi)<0$ imply $|\mathrm{Im}(\xi)| \ge |\xi|\cos\delta$, and so
	\begin{equation}
	\norm{(b+\xi - T^0_0)^{-1}} = \frac{1}{|\text{Im}(b+\xi)|}\leq \frac{1}{|\xi|\cos\delta}.
	\end{equation}
	These imply that for all $z \in U$ and $\xi \in \Sigma_{\delta,R}$, $\norm{zV(b+\xi - T^0_0)^{-1}}\le\frac{1}{2}$, and therefore $I-zV(b+\xi - T^0_0)^{-1}$ has an inverse with norm at most $2$. Due to the factorization $\xi - A^0_z = 
	(I - zV(b+\xi-T^0_0)^{-1}) (b+\xi-T^0_0)$, it holds $\xi \in \rho(A^0_z)$, $(\xi,z) \in \Sigma_{\delta,R}\times U$, and 
	\begin{equation}
	\norm{(\xi - A^0_z)^{-1}} \le \frac{2\max(2,(\cos\delta)^{-1})}{|\xi|}.
	\end{equation}
	This also serves as an upper bound for the norm $\big\|(\xi - \widetilde{A}^0_z)^{-1}\big\|$.
	
	\medskip
	
	Finally, consider the (compact) region $K = (\overline{\Sigma_{\delta}} \cap \{ |\xi| \leq R\}) \setminus D_2$. By Lemma \ref{CompactResolvent}, upon shrinking $U$ again if necessary, we obtain a uniform bound on $\norm{(\xi - {A}^0_z)^{-1}}$ over $z \in U$ and $\xi \in K$, which also holds for
 	$\big\|(\xi - \widetilde{A}^0_z)^{-1}\big\|$.
\end{proof}

The following two lemmas complete the proof of Proposition \ref{MainPerturbation}.

\begin{lemma} \label{PointwiseLowerOrder}
	Let $\mu=0$ and $\lambda>0$. Then, there is a neighborhood $U\subset\mathbb{C}$ of $\lambda$ such that for all compact $K\subset \mathbb{R}$, there exist $\varepsilon>0$ and $C=C(K)<\infty$ for which
	\begin{eqnarray}
			&& \big\|e^{tT^0_z}(I-P^0_z)\big\| 
			\le Ce^{t(\Lambda(\lambda)-\varepsilon)},\quad (t,z)\in[0,\infty)\times U, \label{eq: near 0} \\
			&& \big|e^{t T^0_z}(I-P^0_z)f(x)\big| 
			\le C\norm{f}_{L^2} e^{t(\Lambda(\lambda)-\varepsilon)},\quad (f,x,t,z)\in L^2(\mathbb{R})\times K\times[1,\infty)\times U.
	\end{eqnarray}
\end{lemma}

\begin{proof}
	By \cite[Proposition II.4.4]{Engel1999-tt} and our Lemma \ref{lowerorderresolvent}, the semigroup generated by $\widetilde{A}^0_z$ reads
	\begin{equation} \label{DunfordIntegral}
		e^{t \widetilde{A}^0_z} 
		= \frac{1}{2 \pi i}\int_{\gamma} e^{\xi t}(\xi - \widetilde{A}^0_z)^{-1}\,\mathrm{d}\xi,
	\end{equation}
	where we choose $\gamma$ as in \cite[Definition II.4.2]{Engel1999-tt}. Now, in the proof of the uniform convergence of the integral following that definition, the norm bounds of the  integral depend only on the angle $\delta$ and the constant $M_\delta$. Consequently, there exists a constant $C<\infty$ for which $\big\|e^{t \widetilde{A}^0_z}\big\|\leq C$, $(t,z)\in[0,\infty)\times U$. Since $\widetilde{A}^0_z \equiv A^0_z$ on the range of $I-P^0_z$, we have $e^{t\widetilde{A}^0_z}(I-P^0_z) = e^{tA^0_z}(I-P^0_z)$ on $L^2(\mathbb{R})$. Moreover, $e^{tA^0_z}=e^{-bt} \, e^{t T^0_z}$, which results in
	\begin{equation}
	\big\|e^{tT^0_z}(I-P^0_z)\big\| 
	\le Ce^{bt}\sup_{z \in U}\norm{I-P^0_z} 
	= C\sup_{z \in U}\norm{I-P^0_z} \, e^{t(\Lambda(\lambda)-\delta_1)}.
	\end{equation}
	
	\smallskip
	
	Next, we upgrade this to a pointwise bound. (Recall from Proposition \ref{MainGeneration} that we may choose $e^{tT^0_z}(I-P^0_z)f$ as a continuous function of $x$, so that its pointwise evaluation makes sense.) For any $s\in(0,1]$, $g \in L^2(\mathbb{R})$, and $x \in K$, we can bound $|e^{sT^0_z}g(x)|$ by applying the Feynman-Kac formula and the Cauchy-Schwarz inequality as follows:
	\begin{equation}
	\bigg|\mathbb{E}^0_x\bigg[\exp\bigg(z \int_0^s V(X_r)\,\mathrm{d}r\bigg)\, g(X_s) \bigg] \bigg| 
	\le e^{\text{Re}(z)}\int_{\mathbb{R}} p_s(y-x) \, |g(y)|\, \mathrm{d}y 
	\le C' e^{\text{Re}(z)}s^{-\frac14}\norm{g}_{L^2},
	\end{equation}
	where $p_s$ is the Gaussian density with mean $0$ and variance $s$.~Consequently, for any $(f,x,t,z)\in L^2(\mathbb{R})\times K\times[1,\infty)\times U$, 
	\begin{equation}
	\begin{split}
	\big|e^{t T^0_z}(I-P^0_z)f(x)\big|  =\big|e^{\frac{1}{2}T^0_z}e^{(t-\frac{1}{2})T^0_z}(I-P^0_z)f(x)\big| 
	&\le C\,\big\|e^{(t-\frac{1}{2})T^0_z}(I-P^0_z)f\big\|_{L^2} \\
	&\le C\,\|f\|_{L^2}\, e^{t(\Lambda(\lambda)-\delta_1)},
	\end{split}
	\end{equation}
	finishing the proof of the lemma. 
\end{proof}

\begin{lemma} \label{MedFrequencyBounds}
	Let $\mu=0$ and $\lambda>0$. Then, for any $0<m<M< \infty$, there exists an $\varepsilon>0$ such that
	\begin{equation} \label{medfreqoperatorbound}
		\sup_{y\in\mathbb{R}:\, m \leq |y| \leq M} s(T^0_{\lambda+iy}) \le \Lambda(\lambda)-2\varepsilon.
	\end{equation}
	In addition, for any compact $K\subset\mathbb{R}$, there exists a $C=C(K,m,M)<\infty$ such that 
	\begin{eqnarray}
		&& \sup_{y\in\mathbb{R}:\, m \leq |y| \leq M} \big\|e^{t T^0_{\lambda+iy}}\big\|
		\le C e^{t(\Lambda(\lambda)-\varepsilon)}, \quad t\ge0, \label{medfreq1} \\
		&& \sup_{y\in\mathbb{R}:\, m \leq |y| \leq M} \big|e^{t T^0_{\lambda+iy}} f(x)\big| 
		\le C \,\norm{f}_{L^2}\, e^{t(\Lambda(\lambda)-\varepsilon)}, \quad (f,x,t)\in L^2(\mathbb{R}) \times K\times [1,\infty). \label{medfreq2}
	\end{eqnarray}
\end{lemma}

\begin{proof}
	We first argue that $s(T^0_{\lambda+iy})<\Lambda(\lambda)$ for all $y\in\mathbb{R}\setminus\{0\}$. Suppose for contradiction that $s(T^0_{\lambda+iy}) \ge \Lambda(\lambda)$. We claim that there is an eigenvalue $\alpha$ of $T^0_{\lambda+iy}$ with $\mathrm{Re}(\alpha)\ge\Lambda(\lambda)$. If $s(T^0_{\lambda+iy}) > \Lambda(\lambda)$, this follows directly from Proposition \ref{WeylTheorem}. If $s(T^0_{\lambda+iy})=\Lambda(\lambda)$, then there must be a sequence of points $\alpha_n \in \sigma(T^0_{\lambda+iy})$ with $\mathrm{Re}(\alpha_n)\uparrow \Lambda(\lambda)$. By Proposition~\ref{WeylTheorem}, the $\alpha_n$ are eigenvalues of $T^0_{\lambda+iy}$ for all large enough $n$. Let $(\phi_n)_{n\in\mathbb{N}}$ be corresponding eigenfunctions of norm $1$.  We have that
	\begin{eqnarray}
	&& \quad\mathrm{Re}(\alpha_n)  = \mathrm{Re}\,\langle T^0_{\lambda+iy} \phi_n, \phi_n \rangle =  -\frac{1}{2}\int_{-\infty}^\infty |\phi_n'(x)|^2 \,\mathrm{d}x 
	+ \lambda \int_{-1}^1 |\phi_n(x)|^2\,\mathrm{d}x 
	= \langle T^0_\lambda \phi_n, \phi_n \rangle, \\
	&& \quad \mathrm{Im}(\alpha_n) 
	= \mathrm{Im}\,\langle T^0_{\lambda+iy} \phi_n, \phi_n \rangle 
	=  y \int_{-1}^1 |\phi_n(x)|^2\,\mathrm{d}x.
	\end{eqnarray}
	In particular, $(\alpha_n)_{n\in\mathbb{N}}$ is bounded and admits a convergent subsequence. Since $\sigma(T^0_{\lambda+iy})$ is closed, its limit $\alpha$ belongs to $\sigma(T^0_{\lambda+iy})$. In view of Proposition \ref{WeylTheorem}, $\alpha$ is an eigenvalue of $T^0_{\lambda+iy}$ satisfying $\mathrm{Re}(\alpha)=\Lambda(\lambda)$. We write $\phi$ for a corresponding eigenfunction of norm $1$. 
	
	\medskip
	
	The spectral theorem in the form of \cite[Theorem I.4.9]{Engel1999-tt} renders $T^0_{\lambda}$ unitarily equivalent, via some $G_\lambda$, to multiplication $M_{q_\lambda}$ by a real-valued function $q_\lambda$ on some $L^2(\Omega,\mu)$-space. In addition, $\sigma(T^0_\lambda)$ is given by the essential range of $q_\lambda$, and so the essential range of $q_\lambda$ is upper bounded by $\Lambda(\lambda)$. With $\psi := G_\lambda \phi$, we find that 
	\begin{equation}
	\text{Re}(\alpha) = \int_\Omega q_\lambda\,|\psi|^2\,\mathrm{d}\mu 
	\leq \Lambda(\lambda),
	\end{equation}
	which forces $\mathrm{Re}(\alpha)=\Lambda(\lambda)$. But equality in the latter inequality means that $\psi$ is supported on $\{q_\lambda=\Lambda(\lambda)\}$. Hence, $M_{q_\lambda}\psi = \Lambda(\lambda)\psi$. Consequently, $\phi$ is an eigenfunction of $T^0_\lambda$ with eigenvalue $\Lambda(\lambda)$. Simplicity of $\Lambda(\lambda)$ implies that $\phi = c \phi_\lambda$ for some $c\in \mathbb{C}$ with $|c|=1$, where $\phi_\lambda$ is the positive eigenfunction of $T^0_\lambda$ of norm $1$ corresponding to $\Lambda(\lambda)$. Now, combining the strict positivity of $\phi_\lambda$ with the identities $T^0_\lambda \phi_\lambda = \Lambda(\lambda)\phi_\lambda$ and $T^0_{\lambda+iy}\phi_\lambda = \alpha \phi_\lambda$, we arrive at $iyV = \alpha-\Lambda(\lambda)$ almost everywhere, which is the desired contradiction. 
	
	\medskip
	
	Next, we claim that for each $y_0\in[-M,-m]\cup[m,M]$ there are a neighborhood $U\subset\mathbb{R}$ of $y_0$ and an $\varepsilon_0>0$ such that $s(T^0_{\lambda+iy})\le \Lambda(\lambda)-\varepsilon_0$, $y \in U$. Indeed, choose $\varepsilon_0>0$ to ensure $s(T^0_{\lambda+iy_0}) \le \Lambda(\lambda)-2\varepsilon_0$. Then, the compact $K := \{a +ib:\,\Lambda(\lambda)-\varepsilon_0\leq a \leq \lambda,\,|b| \leq 2M\}$ is contained in $\rho(T^0_{\lambda+iy_0})$. Lemma \ref{CompactResolvent} then yields a neighborhood $U\subset\mathbb{R}$ of $y_0$ satisfying $K\subset \rho(T^0_{\lambda + iy})$, $y\in U$. On the other hand, any $\alpha \in \sigma(T^0_{\lambda+iy})$ with $\mathrm{Re}(\alpha)\ge\Lambda(\lambda)-\varepsilon_0$ must lie in $K$ by the  argument in the third paragraph of the proof of Lemma \ref{topeig}. Therefore, all elements of $\sigma(T^0_{\lambda+iy})$, for $y \in U$, have real parts less than $\Lambda(\lambda)-\varepsilon_0$. Using the compactness of $\{y\in\mathbb{R}:\, m \leq |y| \leq M\}$, we obtain the inequality \eqref{medfreqoperatorbound}.
	
	\medskip
	
	To deduce \eqref{medfreq1} and \eqref{medfreq2}, we proceed as in the proof of Lemma \ref{PointwiseLowerOrder}. First, note that if $0<\delta<\arctan(\frac{\varepsilon}{2M})$, then the sector $\Sigma_{\delta}$ (see \eqref{sectordef}) is contained in $\rho(T^0_{\lambda+iy}-(\Lambda(\lambda)-\varepsilon))$ for all $y\in[-M,-m]\cup[m,M]$.
	Following the third paragraph in the proof of Lemma \ref{lowerorderresolvent}, only with $b:=\Lambda(\lambda)-\varepsilon$ and $C:=\sup_{y\in\mathbb{R}:\,m\le|y|\le M}  |\lambda + iy|$, we obtain 
	\begin{equation}
	\norm{(\xi - (T^0_{\lambda+iy}-b))^{-1}} \le \frac{2\max(2,(\cos\delta)^{-1})}{|\xi|},\quad \xi\in\Sigma_{\delta,R},\quad y\in[-M,-m]\cup[m,M]. 
	\end{equation}
	On the other hand, for the compact region $K := \overline{\Sigma_\delta} \cap \{|\xi|\leq R\}$, Lemma \ref{CompactResolvent} shows that each $y_0\in[-M,-m]\cup[m,M]$ admits a neighborhood $U_{y_0}\subset\mathbb{R}$ and a constant $C_{y_0}<\infty$ for which 
	\begin{equation}
	\sup_{y \in U_{y_0}} \, \sup_{\xi \in K}\, 
	\norm{(\xi - (T^0_{\lambda+iy}-b))^{-1}}\leq C_{y_0}.
	\end{equation}
	The compactness of $\{y\in\mathbb{R}:\, m \leq |y| \leq M\}$ then yields a uniform resolvent bound. Using the resolvent bounds as in the first paragraph of the proof of Lemma \ref{PointwiseLowerOrder}, only with $T^0_{\lambda+iy}-b$ in place of $\widetilde{A}^0_z$, we find a constant $C<\infty$ such that 
	\begin{equation}
	e^{-t(\Lambda(\lambda)-\varepsilon)}\big\|e^{t T^0_{\lambda+iy}}\big\|
	= \big\|e^{t (T^0_{\lambda+iy}-(\Lambda(\lambda)-\varepsilon))}\big\| \le C,
	\quad t\ge0,\quad y\in[-M,-m]\cup[m,M]. 
	\end{equation}
	This can be upgraded to the pointwise bound \eqref{medfreq2} by following the last paragraph in the proof of Lemma \ref{PointwiseLowerOrder}.
\end{proof}

We will deduce Theorems \ref{DriftGeneration} and \ref{DriftSpectrum} for $\mu>0$ from the $\mu=0$ case. For this purpose, we define the operator $R_\mu:\, L^2\rightarrow L^2_\mu$, $f\mapsto e^{-\mu x}f$. It is readily checked that $R_\mu$ is unitary and that $R_\mu H^2=H^2_\mu$.

\medskip

\noindent\textit{Proof of Theorems \ref{DriftGeneration} and \ref{DriftSpectrum} for $\mu>0$.}
	The key observation is that the operator $T_z^\mu$ is unitarily equivalent, via $R_\mu$, to the operator $T_z^0-\frac{\mu^2}{2}$. We prove the statements in Theorems \ref{DriftGeneration} and \ref{DriftSpectrum} sequentially. For Theorem \ref{DriftGeneration}(a), note that if a closed operator $A$ with a dense domain $\mcal{D}(A)$ generates an analytic semigroup on a Hilbert space $H$, and $R:\,H\to H'$ is a unitary operator onto a Hilbert space $H'$, then the operator $A':=RAR^{-1}$ with domain $\mcal{D}(A'):= R\mcal{D}(A)$ is closed and densely defined, and generates the analytic semigroup $Re^{\zeta A}R^{-1}$ on $H'$. In our case, $A=T^0_z-\frac{\mu^2}{2}$ and $R = R_\mu$.
	
	For Theorem \ref{DriftGeneration}(b), fix $f\in L^2_\mu$ and take a function $F$ of the described regularity such that $e^{t T^0_z} R_{\mu}^{-1}f(x) = F(t,x)$ for almost every $(t,x)\in(0,\infty)\times\mathbb{R}$. Then the function $G(t,x):=e^{-\mu x-\frac{\mu^2}{2}t} F(t,x)$ has the desired properties. 
	
	For Theorem \ref{DriftGeneration}(c), we infer from the conjugation formula for $e^{tT_z^\mu}$ that 
	\begin{equation}
	e^{t T_z^{\mu}}f(x) = \mathbb{E}^0_x\bigg[e^{\mu B_t - \frac{\mu^2}{2}t} \, f(X_t)\,
	\exp\bigg(z \int_0^t V(X_s)\,\mathrm{d}s\bigg)\bigg],\quad f \in L^2_{\mu}.
	\end{equation}
	Now, note that $e^{\mu B_t - \frac{\mu^2}{2}t}$ is the Girsanov density for the change of measure from $\mathbb{P}_x^0$ to $\mathbb{P}_x^\mu$.
	
	\medskip
	
	We move on to the proof of Theorem \ref{DriftSpectrum}. The self-adjointness of $T^\mu_z$ for $z\in\mathbb{R}$ follows immediately from that of the unitary equivalent $T^0_z-\frac{\mu^2}{2}$. Moreover, the natural conjugation relation between the resolvents of $T^\mu_\lambda$ and $T^0_\lambda-\frac{\mu^2}{2}$ yields \eqref{eq:spectral bound} and the subsequent statement in Theorem \ref{DriftSpectrum}.
	
	For Theorem \ref{DriftSpectrum}(a), it suffices to use the conjugation relation between the resolvents of $T^\mu_z$ and $T^0_z-\frac{\mu^2}{2}$.
	
	For Theorem \ref{DriftSpectrum}(b), we start with the left-hand side of \eqref{eq:decomp} and sequentially apply $e^{t T_z^{\mu}}=R_\mu e^{t(T^0_z-\mu^2/2)} R_\mu^{-1}$, the identity \eqref{eq:decomp} for $\mu=0$, and twice $P^0_z=R_\mu^{-1}P^\mu_zR_\mu$. 
	
	For Theorem \ref{DriftSpectrum}(c), (d), and (e), we note $P^\mu_zf(x)=e^{-\mu x}P_z^0R_\mu^{-1}f(x)$, $e^{tT^\mu_z}(I-P^\mu_z)f(x)=R_\mu e^{-t\mu^2/2}e^{tT^0_z} (I-P^0_z)R_\mu^{-1}f(x)$, and $e^{tT^\mu_{\lambda+iy}}f(x)=R_\mu e^{-t\mu^2/2}e^{tT^0_{\lambda+iy}}R_\mu^{-1}f(x)$, and appeal to Theorem \ref{DriftSpectrum}(c), (d), and (e) for $\mu=0$, respectively. \qed

\subsection{The SDE \eqref{newsde} and the transcendental equation \eqref{trans1}} \label{subsec: SDE trans}

Before analyzing the SDE \eqref{newsde}, we turn our attention to the transcendental equation \eqref{trans1}.  

\begin{lemma}\label{translemma}
Letting $\Lambda(0):=0$ and restricting the argument of $\tan$ in \eqref{trans1} to $\big(0,\frac{\pi}{2}\big)$ for $\lambda>0$, the equation \eqref{trans1} defines a bijective continuous map $\Lambda:\,[0,\infty)\to[0,\infty)$ satisfying $\Lambda\in C^2((0,\infty))$, $\Lambda'((0,\infty))=(0,1)$, and $\Lambda''>0$ on $(0,\infty)$.
\end{lemma}

\begin{proof}
	It is convenient to introduce $\gamma(\lambda):=\sqrt{2(\lambda-\Lambda(\lambda))}$. Then the equation \eqref{trans1} can be recast as $\lambda=\frac{\gamma^2}{2}(\sec\gamma)^2=:f(\gamma)$, and $\Lambda(\lambda)=\lambda-\frac{\gamma^2}{2}$. By a direct computation, we find that $f'(\gamma) = \gamma (1+ \gamma\tan\gamma)(\sec\gamma)^2$, which is positive on $\big(0,\frac\pi 2\big)$. In addition, $\lim_{\gamma\downarrow0} f(\gamma)=0$ and $\lim_{\gamma\uparrow\frac\pi 2} f(\gamma)=\infty$. Thus, for $\lambda>0$, the corresponding $\Lambda(\lambda)$ is obtained by taking the unique solution of $f(\gamma)=\lambda$ in $\big(0,\frac\pi 2\big)$ and setting $\Lambda(\lambda)=\lambda-\frac{\gamma^2}{2}=\frac{\gamma^2}{2}(\tan\gamma)^2=:F(\gamma)$. The resulting map $\Lambda:\,[0,\infty)\to[0,\infty)$ is bijective and continuous, as it is a composition of the strictly increasing continuous maps $F$ and $f^{-1}$.  
	
	\medskip	
	
	Next, we calculate:
	\begin{equation}
	\Lambda'(\lambda) = \frac{\mathrm{d}}{\mathrm{d}\lambda}F(f^{-1}(\lambda)) 
	= \frac{F'(f^{-1}(\lambda))}{f'(f^{-1}(\lambda))} 
	= \frac{\gamma \tan\gamma\,(\tan\gamma  + \gamma (\sec\gamma)^2)}
	{\gamma (1+ \gamma\tan\gamma)(\sec\gamma)^2}
	= \frac{(\sin\gamma)^2 + \gamma \tan\gamma}{1+\gamma \tan \gamma}.
	\end{equation}
	Since $\lim_{\lambda\downarrow0}\gamma(\lambda)=0$, we have that $\lim_{\lambda\downarrow0}\Lambda'(\lambda)=0$. On the other hand, $\lim_{\lambda\uparrow\infty} \gamma(\lambda)=\frac\pi2$, and so $\lim_{\lambda\uparrow\infty} \Lambda'(\lambda)=1$.
	Using the identity $( \frac{g+h}{1+h})' = \frac{g'(1+h)+(1-g)h'}{(1+h)^2}$, we find further that
	\begin{equation}
	\frac{\mathrm{d}}{\mathrm{d}\gamma} \frac{(\sin\gamma)^2 + \gamma \tan\gamma}{1+\gamma \tan \gamma} 
	= \frac{2\sin\gamma\cos\gamma\,(1+\gamma \tan\gamma) + (\cos\gamma)^2\,
		(\tan\gamma + \gamma (\sec\gamma)^2)}{(1+\gamma\tan\gamma)^2},
	\end{equation}
	which means that $\Lambda''(\lambda)$ is this expression divided by $\gamma (1+ \gamma\tan\gamma)(\sec\gamma)^2$. Each term is readily seen to be positive for $\gamma \in \big(0,\frac\pi 2\big)$, so $\Lambda''>0$ on $(0,\infty)$.
\end{proof}

Next, we study the SDE \eqref{newsde}.

\begin{prop} \label{SDEfacts}
	Let $\lambda>0$. Then: 
	\begin{enumerate}[(a)]
	\item For all initial conditions $x\in\mathbb{R}$, the SDE \eqref{newsde} admits a unique strong solution, whose law we denote by $\mathbb{Q}^x_\lambda$. 
	\item For any $x\in\mathbb{R}$ and $t\ge0$, the measure $\mathbb{Q}^x_\lambda$ is absolutely continuous with respect to~$\mathbb{P}_x^\mu$ on $[0,t]$, and the associated density is given by
	\begin{equation} \label{GirsanovDensity}
		M_t := \frac{h_\lambda(X_t)}{h_\lambda(x)}\, e^{-t\Lambda_{\mu}(\lambda)}\,
		\exp\bigg(\lambda \int_0^t V(X_s)\,\mathrm{d}s \bigg),
	\end{equation}
	where $h_\lambda(x):=e^{-\mu x}\phi_\lambda(x)$.
	\item Define $G_t = \frac{1}{\sqrt{t}}\int_0^t V(X_s)-\Lambda'(\lambda)\,\mathrm{d}s$. Then, for all bounded measurable $g:\,\mathbb{R}\to\mathbb{R}$ and all bounded continuous $F:\,\mathbb{R}\to\mathbb{C}$, it holds
	\begin{equation}
	\lim_{t\to\infty} \,\mathbb{E}^{\mathbb{Q}_\lambda^x}[g(X_t)\, F(G_{t})] 
	= \bigg(\int_{-\infty}^\infty g(y)\,\phi_\lambda(y)^2\,\mathrm{d}y\bigg) \,
	\bigg(\int_{-\infty}^\infty F(y)\,\frac{e^{-\frac{y^2}{2\Lambda''(\lambda)}}}{\sqrt{2\pi \Lambda''(\lambda)}}\,\mathrm{d}y\bigg).
	\end{equation}
	\item With $\rho_T(\cdot):=\frac{1}{T}\int_0^T \mathbf{1}_{\cdot}(X_t)\,\mathrm{d}t$, $\rho_T$ converges weakly to $\pi_\lambda$ as $T\to\infty$, $\mathbb{Q}_\lambda^x$-almost surely.
	\end{enumerate}
\end{prop}

\begin{proof}
	\textit{(a).} The explicit formula \eqref{Eigenvector} for $\phi_\lambda$ reveals that $(\log \phi_\lambda)'$ is a bounded continuous function on $\mathbb{R}$ that is $C^1_b$ on each of $(-\infty,-1)$, $(-1,1)$, and $(1,\infty)$. Thus, $(\log \phi_\lambda)'$ is bounded and globally Lipschitz, guaranteeing strong existence and uniqueness of solutions to the SDE \eqref{newsde} by standard theory (see, e.g., \cite[Chapter 5, Theorems 2.5 and 2.9]{Karatzas2004-ky}). 
	
	\medskip
	
	\noindent\textit{(b).} We now identify the densities of $\mathbb{Q}_\lambda^x$ with respect to $\mathbb{P}_x^{\mu}$ on the time intervals $[0,t]$, $t\ge0$. Starting under $\mathbb{P}_x^\mu$, so that $X_t=x+B_t+\mu t$, $t\ge0$ with a standard Brownian motion~$B$, we define the exponential martingale
	\begin{equation} \label{M def}
	M_t = \exp\bigg(\int_0^t \big((\log \phi_\lambda)'(X_s)-\mu\big)\,\mathrm{d}B_s 
	- \frac{1}{2} \int_0^t \big((\log \phi_\lambda)'(X_s) -\mu\big)^2\,\mathrm{d}s\bigg),\quad t\ge0.
	\end{equation}
	Clearly, $M$ satisfies Novikov's condition (see, e.g., \cite[Chapter 3, Corollary 5.13]{Karatzas2004-ky}), and Girsanov's Theorem (see, e.g., \cite[Chapter 3, Theorem 5.1]{Karatzas2004-ky}) shows that on each time interval $[0,t]$, the process $X$ solves \eqref{newsde} weakly after the change of measure induced by~$M_t$. Moreover, by It\^o's formula in the form of \cite[Theorem 2.1, Remark 2.3]{Peskir2005-xp}, we find that 
	\begin{equation} \label{PeskirIto}
	\log \phi_\lambda(X_t) - \log \phi_\lambda(x) = \int_0^t (\log \phi_\lambda)' (X_s)(\mathrm{d}B_s+\mu\,\mathrm{d}s) 
	+ \frac{1}{2}\int_0^t (\log \phi_\lambda)'' (X_s)\,\mathbf{1}_{\{X_s\notin\{-1,1\}\}}\,\mathrm{d}s.
	\end{equation}
	Furthermore, the relation $\frac{1}{2}\phi_\lambda'' + \lambda V \phi_\lambda = \Lambda(\lambda)\phi_\lambda$ on $\mathbb{R}\setminus\{-1,1\}$ reveals that
	\begin{equation}
	\frac{1}{2}(\log \phi_\lambda)'' = -\lambda V+ \Lambda(\lambda)-\frac{1}{2}( (\log \phi_\lambda)')^2 \quad\text{on}\;\;\mathbb{R}\setminus\{-1,1\}. 
	\end{equation}
	Inserting this into \eqref{PeskirIto}, noting that $\mathbf{1}_{\{X_s\notin\{-1,1\}\}}=1$ for almost every $s\in[0,t]$, $\mathbb{P}_x^\mu$-almost surely, and plugging the resulting identity into \eqref{M def}, we get the right-hand side of \eqref{GirsanovDensity}. 
	
	\medskip
	
	\noindent\textit{(c).} 
	The limits in parts (c) and (d) do not involve $\mu$, so we set $\mu=0$ for the rest of the proof. By considering the positive and negative parts of $g$, it suffices to take $g \ge 0$, and we may assume without loss of generality that $g$ is not identically equal to $0$. We first let $F(y)=e^{i\xi y}$ for some $\xi\in\mathbb{R}$. Then, from \eqref{GirsanovDensity} and Theorem \ref{DriftGeneration}(c), we read off that 
	\begin{equation}
		\begin{split}
			&\;\mathbb{E}^{\mathbb{Q}_\lambda^x}\big[g(X_t)\,F(G_t)\big] \\
			&= \phi_\lambda(x)^{-1}\, e^{-t\Lambda(\lambda)-\sqrt{t}i\xi \Lambda'(\lambda)}\, \mathbb{E}^0_x\bigg[g(X_t)\,\phi_\lambda(X_t)\,
			\exp\bigg(\Big(\lambda+i\frac{\xi}{\sqrt{t}}\Big) 
			\int_0^t V(X_s)\,\mathrm{d}s \bigg) \bigg] \\
			&= \phi_\lambda(x)^{-1}\, e^{-t\Lambda(\lambda)-\sqrt{t}i\xi\Lambda'(\lambda)}\,
			e^{t T^0_{\lambda+i\frac{\xi}{\sqrt{t}}}}(g\phi_\lambda)(x).
		\end{split}
	\end{equation}
	
	\smallskip
	
	By Theorem \ref{DriftSpectrum} (b) and (d), we have that $e^{tT^0_{\lambda+i\frac{\xi}{\sqrt{t}}}} (g\phi_\lambda)(x) 
	= e^{t \Lambda(\lambda + i\frac{\xi}{\sqrt{t}})}\, P^0_{\lambda+i\frac{\xi}{\sqrt{t}}}(g\phi_\lambda)(x)$ up to an error which is negligible 
	when multiplied by $e^{-t\Lambda(\lambda)}$. Now, by Theorem \ref{DriftSpectrum} (c) and the holomorphicity of $\Lambda$ at $\lambda$, as $t\to\infty$,
	\begin{equation*}
	\begin{split}
	\phi_\lambda(x)^{-1}\, e^{t (\Lambda(\lambda+i\frac{\xi}{\sqrt{t}}) - \Lambda(\lambda) - i \frac{\xi}{\sqrt{t}}\Lambda'(\lambda))} \, P^0_{\lambda+i\frac{\xi}{\sqrt{t}}}(g\phi_\lambda)(x)
	\longrightarrow & \; e^{-\frac{\xi^2}{2}\Lambda''(\lambda)}\langle g \phi_\lambda,\phi_\lambda \rangle_{L^2} \\ 
	& =\langle g \phi_\lambda,\phi_\lambda \rangle_{L^2} \!\int_{-\infty}^\infty \!\! e^{i\xi y}\, \frac{e^{-\frac{y^2}{2\Lambda''(\lambda)}}}{\sqrt{2\pi \Lambda''(\lambda)}}\,\mathrm{d}y.
	\end{split}
	\end{equation*}
	Taking $\xi=0$, we see that $m_t:=\mathbb{E}^{\mathbb{Q}_\lambda^x}[g(X_t)]\to \langle g\phi_\lambda, \phi_\lambda \rangle_{L^2}$. Finally, we introduce the probability measures
	\begin{equation}
	\nu_t^g(\cdot) := \frac{1}{m_t}\,\mathbb{E}^{\mathbb{Q}_\lambda^x}[g(X_t)\,\mathbf{1}_{\{G_t \in \,\cdot\,\}}], \qquad \nu_\infty^g(\mathrm{d}y) := \frac{e^{-\frac{y^2}{2\Lambda''(\lambda)}}}{\sqrt{2\pi \Lambda''(\lambda)}} \,\mathrm{d}y.
	\end{equation}
	The above calculations show that the Fourier transforms of $\nu_t^g$ converge pointwise
	to the Fourier transform of $\nu_\infty^g$. By the Lévy Continuity Theorem (see, e.g., \cite[Theorem 6.3]{Kallenberg2021-fh}),  $\nu_t^g$ converges weakly to $\nu_\infty^g$. Recalling that $m_t \rightarrow \langle g\phi_\lambda,\phi_\lambda\rangle_{L^2}$ completes the proof of part~(c). 
	
	\medskip
	
	\noindent\textit{(d).} Under $\mathbb{Q}_\lambda^x$, the process $X$ is a Brownian motion with constant drift towards $0$ outside of $(-1,1)$, and so \cite[Section 4.4, Assumption (B)]{Khasminskii2011-yq} with $U:=(-1,1)$ holds (use, e.g., the explicit hitting time density of Brownian motion with drift in \cite[Chapter 3, display (5.12)]{Karatzas2004-ky}). Then, \cite[Chapter 4, Theorems 4.1 and 4.2]{Khasminskii2011-yq} yield a stationary distribution $\eta$ such that for all $x\in\mathbb{R}$ and all $f \in L^1(\eta)$, $\int_{\mathbb{R}} f \,\mathrm{d}\rho_T \to \int_{\mathbb{R}} f \,\mathrm{d}\eta$, $\mathbb{Q}_\lambda^x$-almost surely. For bounded~$f$, the Dominated Convergence Theorem implies $\mathbb{E}^{\mathbb{Q}_\lambda^x}[\int _{\mathbb{R}}  f\, \mathrm{d}\rho_T]\to\int _{\mathbb{R}} f \,\mathrm{d}\eta$. By part (c), we also have $\mathbb{E}^{\mathbb{Q}_\lambda^x}[\int _{\mathbb{R}}  f\, \mathrm{d}\rho_T]\to\int _{\mathbb{R}} f \,\mathrm{d}\pi_\lambda$, and we conclude that $\eta=\pi_\lambda$. It remains to note that $\int_{\mathbb{R}} \mathbf{1}_{(-\infty,q]}\,\mathrm{d}\rho_T\to\int_{\mathbb{R}} \mathbf{1}_{(-\infty,q]}\,\mathrm{d}\pi_\lambda$ for all rational $q$, $\mathbb{Q}_\lambda^x$-almost surely.
\end{proof}

\section{Proof of Theorem \ref{SharpAsymptotics}(a)} \label{sec: super}

We fix a $\mu \ge 0$ throughout this section. The following lemma is the key to our proof of Theorem \ref{SharpAsymptotics}(a).

\begin{lemma} \label{supersmoothfourier}
	Let $\lambda>\Lambda^{-1}\big(\frac{\mu^2}{2}\big)=:\lambda_c$. For $t\ge0$ and $x,\xi\in\mathbb{R}$, define
	\begin{equation} 
	A_t^x(\xi) = \frac{1}{e^{-\mu x}\phi_\lambda(x)}\, e^{-t\Lambda_\mu(\lambda) - ti \xi \Lambda'(\lambda)}\,\mathbb{E}^{\mu}_x\bigg[\exp\bigg((\lambda+i\xi)\int_0^t V(X_s)\,\mathrm{d}s\bigg)\bigg].
	\end{equation}
	In addition, for any $\widehat{\psi} \in C_c(\mathbb{R})$ and $b\in\mathbb{R}$, define $\widehat{\psi}_b(\xi)=e^{-i\xi b}\,\widehat{\psi}(\xi)$. Then for any compact $K\subset \mathbb{R}$, as $t \rightarrow \infty$,
	\begin{equation} \label{A asymp}
	\sup_{x \in K}\,\sup_{b \in \mathbb{R}}\,\bigg|\sqrt{\frac{t}{2\pi}}\int_{-\infty}^\infty \widehat{\psi}_b(\xi)\, A_t^x(\xi)\,\mathrm{d}\xi - \gamma_\lambda\,\widehat{\psi}(0)\, e^{-\frac{b^2}{2t\Lambda''(\lambda)}}\bigg| \longrightarrow 0,
	\end{equation}
	where $\gamma_\lambda := \frac{\lambda}{\Lambda_{\mu}(\lambda)\sqrt{\Lambda''(\lambda)}}\,\langle R_{\mu}^{-1}V,\phi_{\lambda}\rangle_{L^2}$.
	Moreover,
	\begin{equation}
	\sup_{t \geq 1}\,\sup_{x \in K}\,\sup_{b \in \mathbb{R}}\,\bigg|\sqrt{\frac{t}{2\pi}} \int_{-\infty}^{\infty} \widehat{\psi}_b(\xi)\, A_t^x(\xi)\,\mathrm{d}\xi\bigg|<\infty.
	\end{equation}
\end{lemma}

\begin{proof}
	We will sequentially break off negligible pieces of the integral in \eqref{A asymp} until we arrive at its leading-order term, checking the uniform boundedness of the pieces for $t\ge1$ along the way. Let $z:= \lambda+i\xi$ and $h_{\lambda}(x):=e^{-\mu x}\phi_{\lambda}(x)$. First, we note that the Fundamental Theorem of Calculus, Fubini's Theorem, and Theorem \ref{DriftGeneration}(c) yield
	\begin{equation} \label{Duhamel}
		\mathbb{E}^\mu_x\bigg[\exp\bigg(z\int_{0}^t V(X_s)\,\mathrm{d}s\bigg)\bigg] 
		= 1 + z\int_0^t 
		\mathbb{E}^\mu_x\Big[V(X_s)\, e^{z\int_0^s V(X_r)\,\mathrm{d}r }\Big]\,\mathrm{d}s 
		= 1 + z\int_0^t e^{s T_z^{\mu}}V(x)\,\mathrm{d}s.
	\end{equation}
	The contribution of the term $1$ to the first expression in \eqref{A asymp} can be bounded as follows: 
	\begin{equation}\label{1bound}
	\sup_{x \in K}\,\sup_{b \in \mathbb{R}}\,\bigg| \sqrt{\frac{t}{2\pi}}\, h_\lambda(x)^{-1} \int_{-\infty}^\infty 
	e^{-t\Lambda_\mu(\lambda) - ti\xi\Lambda'(\lambda)}
	\,\widehat{\psi}_b(\xi)\,\mathrm{d}\xi \bigg| 
	\leq C(K)\big\|\widehat{\psi}\big\|_{L^1}\sqrt{t}e^{-t\Lambda_{\mu}(\lambda)}\to 0.
	\end{equation}
	Moreover, the Feynman-Kac formula implies that $\sup_{0\le s\le 1,\,\xi\in\mathbb{R},\,x\in K} |e^{sT_z^\mu}V(x)|\le e^\lambda$, so
	\begin{equation}
		\begin{split}
	&\,\sup_{x \in K}\sup_{b \in \mathbb{R}} \bigg|\sqrt{\frac{t}{2\pi}}\, h_\lambda(x)^{-1} \int_{-\infty}^\infty e^{-t\Lambda_{\mu}(\lambda)- ti \xi \Lambda'(\lambda)}\, z
	\int_0^1 e^{sT_z^\mu}V(x)\,\mathrm{d}s\,\widehat{\psi}_b(\xi)\,\mathrm{d}\xi \bigg| \\
	&\le C(K)\big\|z\widehat{\psi}\big\|_{L^1}\sqrt{t}e^{-t\Lambda_{\mu}(\lambda)}\to 0.
		\end{split}
	\end{equation}
	
	\smallskip
	
	Next, we recall $\Lambda''(\lambda)>0$ from Lemma \ref{translemma}. Hence, we can choose a $\delta>0$ small enough to ensure that on the open ball of radius $\delta$ around $\lambda$, the results of Theorem \ref{DriftSpectrum}(a)--(e) hold, and $\mathrm{Re}(\Lambda_{\mu}(z))>0$ along with $\mathrm{Re}(\Lambda_{\mu}(\lambda + i \xi)-\Lambda_{\mu}(\lambda) - i \xi\Lambda'(\lambda))\le -c \xi^2$ for some $c>0$.  
	Let $[-M,M]\supset\mathrm{supp}(\widehat{\psi})$. By Theorem \ref{DriftSpectrum}(e), there exist $\varepsilon>0$ and $C<\infty$, depending on $K$, $\delta$, and $M$, such that for all $\xi\in[-M,-\delta]\cup[\delta,M]$, $s \geq 1$, and $x \in K$,
	\begin{equation}
	\big| z e^{sT_z^{\mu}}V(x)\big| \leq C e^{s(\Lambda_{\mu}(\lambda)-\varepsilon)}.
	\end{equation}
	Without loss of generality, we may assume that $\varepsilon\le\Lambda_{\mu}(\lambda)/2$. It follows that
	\begin{equation}
	\bigg| z\int_1^t e^{sT_{z}^{\mu}}V(x)\,\mathrm{d}s \bigg| 
	\le C' e^{t(\Lambda_{\mu}(\lambda) - \varepsilon)}.
	\end{equation}
	We conclude that 
	\begin{equation}
	\sup_{x \in K}\sup_{b \in \mathbb{R}}\bigg |\sqrt{\frac{t}{2\pi}}\, h_\lambda(x)^{-1} \int_{\delta \leq |\xi| \leq M} \! e^{-t\Lambda_{\mu}(\lambda)-ti \xi \Lambda'(\lambda)} z \int_1^t \! e^{s T_z^\mu}V(x)\,\mathrm{d}s \,\widehat{\psi}_b(\xi)\,\mathrm{d}\xi \bigg| 
	\le C\big\|\widehat{\psi}\big\|_{L^1}\sqrt{t}e^{-t\varepsilon}
	\end{equation}
	which vanishes as $t\to\infty$. An analogous argument, relying on Theorem \ref{DriftSpectrum}(d), yields
	\begin{equation}
	\sup_{x \in K}\sup_{b \in \mathbb{R}} \bigg|\sqrt{\frac{t}{2\pi}}\, h_\lambda(x)^{-1} \int_{-\delta}^\delta e^{-t\Lambda_{\mu}(\lambda)-ti \xi \Lambda'(\lambda)} z 
	\int_1^t e^{s T_z^{\mu}}Q_z^{\mu}V(x)\,\mathrm{d}s \, \widehat{\psi}_b(\xi)\,\mathrm{d}\xi \bigg| 
	\le C\big\|\widehat{\psi}\big\|_{L^1}\sqrt{t}e^{-t\varepsilon}.
	\end{equation}
	
	\smallskip
	
	It remains to analyze the contribution of the term
	\begin{equation}
	z \int_1^t e^{s \Lambda_{\mu}(z)}P_z^{\mu}V(x)\,\mathrm{d}s = \frac{z}{\Lambda_{\mu}(z)}\big(e^{t\Lambda_{\mu}(z)}-e^{\Lambda_\mu(z)}\big)P_z^{\mu}V(x)
	\end{equation}
	to the $\xi$-integral over $(-\delta,\delta)$. The quantity $-\frac{z}{\Lambda_{\mu}(z)}e^{\Lambda_{\mu}(z)}P_z^{\mu}V(x)$ is bounded in absolute value by a constant thanks to Theorem \ref{DriftSpectrum}(a) and (c), so its contribution can be estimated similarly to \eqref{1bound}. Thus, we are left to study
	\begin{equation}
	\sqrt{\frac{t}{2\pi}}\, h_{\lambda}(x)^{-1}\!\!
	\int_{-\delta}^\delta \! \frac{z}{\Lambda_{\mu}(z)}P^{\mu}_zV(x)\, e^{-i\xi b} \widehat{\psi}(\xi) \, e^{t\eta(\xi)}\,\mathrm{d}\xi
	= \frac{h_{\lambda}(x)^{-1}}{\sqrt{2\pi}}
	\!\! \int_{-\sqrt{t}\delta}^{\sqrt{t}\delta} \! a_x\Big(\frac{u}{\sqrt{t}}\Big)\, 
	e^{-i\frac{u}{\sqrt{t}} b}e^{t\eta(\frac{u}{\sqrt{t}})}\mathrm{d}u,
	\end{equation}
	where $\eta(\xi) := \Lambda_{\mu}(\lambda + i \xi) - \Lambda_{\mu}(\lambda)-i\xi \Lambda'(\lambda)$ and $a_x(\xi) := \frac{\lambda + i \xi}{\Lambda_{\mu}(\lambda+ i \xi)}P^{\mu}_{\lambda + i\xi}V(x)\,\widehat{\psi}(\xi)$.
	
	\medskip
	
	We split the domain of integration in the last integral into $R^1_t:=[-t^{1/8}\delta,t^{1/8}\delta]$ and $R^2_t:=(-\sqrt{t}\delta,\sqrt{t}\delta)\setminus R^1_t$. On $R^2_t$, the quantity $a_x(u/\sqrt{t})$ is bounded in absolute value by a constant thanks to Theorem \ref{DriftSpectrum}(a) and (c), and $\mathrm{Re}(\eta(u/\sqrt{t}))\le -cu^2/t$. Thus,
	\begin{equation}
	\sup_{x \in K}\,\sup_{b \in \mathbb{R}}\, 
	\bigg|\frac{h_{\lambda}(x)^{-1}}{\sqrt{2\pi}} 
	\int_{R^2_t} a_x\Big(\frac{u}{\sqrt{t}}\Big)\, e^{-i\frac{u}{\sqrt{t}} b}e^{t\eta(\frac{u}{\sqrt{t}})}\,\mathrm{d}u \bigg|
	\le C\int_{t^{1/8}\delta}^\infty e^{-cu^2}\,\mathrm{d}u\to 0.
	\end{equation}
	Moreover, by Theorem \ref{DriftSpectrum}(a) and (c), $a_x(\xi) \to a_x(0)$
	uniformly over $x\in K$. Therefore,
	\begin{equation}
	\sup_{x\in K}\sup_{b\in\mathbb{R}}\bigg|\frac{h_{\lambda}(x)^{-1}}{\sqrt{2\pi}} 
	\!\!\int_{R^1_t}\! \Big(a_x\Big(\frac{u}{\sqrt{t}}\Big)-a_x(0)\!\Big)\, e^{-i\frac{u}{\sqrt{t}} b}e^{t\eta(\frac{u}{\sqrt{t}})}\mathrm{d}u \bigg|
	\le C\sup_{x \in K} \sup_{|\xi|\le t^{-\frac{3}{8}}\delta} |a_x(\xi)-a_x(0)| \rightarrow 0,
	\end{equation}
	where we reused the inequality $\mathrm{Re}(\eta(u/\sqrt{t}))\le -cu^2/t$. Furthermore, 
	the holomorphicity of $\Lambda_\mu$ (see Theorem \ref{DriftSpectrum}(a)) yields 
	a constant $C<\infty$ such that $|\eta(\xi) + \frac{1}{2}\xi^2\Lambda''(\lambda)|\leq C|\xi|^3$, $\xi\in(-\delta,\delta)$. Thus,  $|e^{t\eta(u/\sqrt{t})}-e^{-\frac{1}{2}u^2\Lambda''(\lambda)}|\le e^{-\frac{1}{2}u^2\Lambda''(\lambda)}|e^{Ct^{-1/8}\delta^3}-1|$, $u \in R^1_t$, and so
	\begin{equation}
	\sup_{x \in K}\,\sup_{b \in \mathbb{R}}\, \bigg|\frac{h_{\lambda}(x)^{-1}}{\sqrt{2\pi}}\, a_x(0)
	\int_{R^1_t} e^{-i\frac{u}{\sqrt{t}} b}\big(e^{t\eta(\frac{u}{\sqrt{t}})}-e^{-\frac{1}{2}u^2\Lambda''(\lambda)}\big)\,\mathrm{d}u \bigg|
	\le C\big|e^{Ct^{-\frac{1}{8}}}-1\big|\to 0.
	\end{equation}
	
	\smallskip
	
	Finally, we may add in the integral over $\mathbb{R}\setminus R^1_t$ thanks to the bound 
	\begin{equation}
	\sup_{x \in K}\,\sup_{b \in \mathbb{R}}\, \bigg|\frac{h_{\lambda}(x)^{-1}}{\sqrt{2\pi}}\, a_x(0)
	\int_{\mathbb{R}\setminus R^1_t} e^{-i\frac{u}{\sqrt{t}} b}\, e^{-\frac{1}{2}u^2\Lambda''(\lambda)}\,\mathrm{d}u \bigg|
	\le C\int_{t^{1/8}\delta}^\infty e^{-\frac{1}{2}u^2\Lambda''(\lambda)}\,\mathrm{d}u\to0.
	\end{equation}
	Using the definition of $a_x(0)$ and the identity $P_\lambda^\mu V(x) = \langle R_{\mu}^{-1}V,\phi_{\lambda}\rangle_{L^2}\,h_{\lambda}(x)$, we conclude that $\sqrt{\frac{t}{2\pi}}\int_{-\infty}^\infty \widehat{\psi}_b(\xi)\, A_t^x(\xi)\,\mathrm{d}\xi$ equals 
	\begin{equation}
	\frac{h_{\lambda}(x)^{-1}}{\sqrt{2\pi}}\, a_x(0) \int_{\mathbb{R}}e^{-i\frac{u}{\sqrt{t}}b}e^{-\frac{1}{2}u^2\Lambda''(\lambda)}\,\mathrm{d}u 
	= \frac{\lambda}{\Lambda_{\mu}(\lambda)\sqrt{\Lambda''(\lambda)}}
	\,\langle R_\mu^{-1}V, \phi_\lambda \rangle_{L^2} \,\widehat{\psi}(0)\, e^{-\frac{b^2}{2t\Lambda''(\lambda)}},
	\end{equation}
	up to an $o(1)$ error uniform over $(x,b) \in K\times\mathbb{R}$.
\end{proof}

The next two lemmas provide the bridge between Lemma \ref{supersmoothfourier} and Theorem \ref{SharpAsymptotics}(a).

\begin{lemma} \label{uniformbound}
	For $t\ge0$ and $x\in\mathbb{R}$, and with $Z_t := \int_0^t V(X_s)-\Lambda'(\lambda)\,\mathrm{d}s$, define the measures
	\begin{equation} \label{nu def}
	\nu_t^x(\cdot) 
	= \sqrt{t}\,\mathbb{E}^{\mathbb{Q}^x_\lambda}\Big[
	\big(e^{-\mu X_t}\phi_{\lambda}(X_t)\big)^{-1}\,\mathbf{1}_{\{Z_t\in\,\cdot\,\}}\Big].
	\end{equation}
	Then, for any compact $K\subset\mathbb{R}$,
	\begin{equation}
	\sup_{t \geq 1}\,\sup_{x \in K}\,\sup_{b \in \mathbb{R}}\,\nu_t^x([b,b+1])<\infty.
	\end{equation}
\end{lemma}

\begin{proof}
	Preparing for the eventual use of Lemma \ref{supersmoothfourier}, we apply Proposition \ref{SDEfacts}(b) to obtain
	\begin{equation} \label{Lemma 3.2 prep}
	\int_{-\infty}^\infty e^{i\xi y}\,\nu_t^x(\mathrm{d}y) 
	= \sqrt{t}\,\mathbb{E}^{\mathbb{Q}^x_\lambda}
	\Big[\big(e^{-\mu X_t}\phi_{\lambda}(X_t)\big)^{-1}\, e^{i\xi Z_t}\Big] = \sqrt{t}A_t^x(\xi).
	\end{equation}
	Now, we take a standard (even) mollifier $\widehat{p}$ in the sense of \cite[first definition in Section C.5]{evans2010partial}. Under the Fourier transform convention $\widehat{\psi}(\xi) = (2\pi)^{-1/2}\int_{-\infty}^\infty e^{-i\xi y}\,\psi(y)\,\mathrm{d}y$, the inverse Fourier transform $p$ of $\widehat{p}$ is continuous and $p(0) = (2\pi)^{-1/2} \int_{-\infty}^\infty \widehat{p}(\xi)\,\mathrm{d}\xi > 0$. Moreover, there exist $a,R\in(0,\infty)$ such that 
	$\psi(y) := ap(y/R)^2 \ge \mathbf{1}_{[0,1]}(y)$, $y\in\mathbb{R}$. In addition, from $\widehat{p}\in C_c^\infty(\mathbb{R})$ it follows that $\widehat{\psi}(\xi) = \frac{aR}{\sqrt{2\pi}} (\widehat{p}*\widehat{p})(R\xi)$ is $C_c^\infty$, where $*$ denotes convolution. Thus, $\psi_b(y) := \psi(y-b) \ge \mathbf{1}_{[b,b+1]}(y)$, $y\in\mathbb{R}$, with Fourier transform $\widehat{\psi}_b(\xi) = e^{-i\xi b}\,\widehat{\psi}(\xi)$. By the identity \eqref{Lemma 3.2 prep}, we have that
	\begin{equation} \label{Fourier eq}
	\nu_t^x([b,b+1]) \le \int_{-\infty}^\infty \psi_b(y)\,\nu_t^x(\mathrm{d}y) 
	= \sqrt{\frac{t}{2\pi}} \int_{-\infty}^\infty \widehat{\psi}_b(\xi)\, A_t^x(\xi)
	\,\mathrm{d}\xi.
	\end{equation}
	The result now follows from Lemma \ref{supersmoothfourier}.
\end{proof}

\begin{lemma}\label{fourierupgrade}  
	Fix a compact $K\subset\mathbb{R}$ and $M>0$. Then, for any bounded and compactly supported $\psi:\,\mathbb{R}\to\mathbb{R}$ that is continuous almost everywhere, it holds
	\begin{equation}\label{fourierupgradeeq}
	\lim_{t\to\infty}\,\sup_{x \in K}\,\sup_{|b|\leq M\sqrt{t}}\,\bigg|\int_{-\infty}^\infty \psi(y-b)\,\nu_t^x(\mathrm{d}y)- \frac{\gamma_\lambda}{\sqrt{2\pi}}\, \int_{-\infty}^\infty \psi(y)\,\mathrm{d}y\,e^{-\frac{b^2}{2t\Lambda''(\lambda)}}\bigg|=0,
	\end{equation}
	with $\gamma_\lambda$ of Lemma \ref{supersmoothfourier}. The limit in \eqref{fourierupgradeeq} also applies to $\psi(y) = e^{-cy}\,\mathbf{1}_{[0,\infty)}(y)$ for $c>0$.
\end{lemma}

\begin{proof}
	If $\widehat{\psi}\in C_c(\mathbb{R})$, then the limit in \eqref{fourierupgradeeq} follows from the equality in \eqref{Fourier eq} (recall that $\psi_b(y)=\psi(y-b)$) and \eqref{A asymp}. Next, we turn to $\psi \in C_c(\mathbb{R})$. It suffices to show that for every $\psi \in C_c(\mathbb{R})$, there are functions $\psi^\varepsilon$ with $\widehat{\psi^\varepsilon} \in C_c(\mathbb{R})$  such that as $\varepsilon \downarrow 0$, $\norm{\psi^\varepsilon-\psi}_{L^1}\rightarrow 0$ and
	\begin{equation}
	\sup_{x\in K,\, |b|\leq M\sqrt{t},\, t\geq 1} \bigg|\int_{-\infty}^\infty \psi(y-b)-\psi^\varepsilon(y-b)\,\nu_t^x(\mathrm{d}y) \bigg| \longrightarrow 0.
	\end{equation}
	Given a $\psi\in C_c(\mathbb{R})$, note that $\widehat{\psi}\in C^\infty(\mathbb{R})$ (see, e.g., \cite[Theorem 8.22]{Folland1999-li}). Take a $\widehat{p} \in C_c^\infty(\mathbb{R})$ with $\widehat{p}(0)=(2\pi)^{-1/2}$. Then,  the inverse Fourier transform $p$ of $\widehat{p}$ is a Schwartz function (see, e.g., \cite[Corollary 8.23]{Folland1999-li}) which integrates to $1$. With $p^\varepsilon(y) := \varepsilon^{-1} p(\varepsilon^{-1}y)$, the functions $\psi^\varepsilon := \psi * p^\varepsilon$ converge uniformly to $\psi$ (see, e.g., \cite[Theorem 8.14]{Folland1999-li}). In addition, $\widehat{\psi^\varepsilon}(\xi) = \sqrt{2\pi}\,\widehat{\psi}(\xi)\,\widehat{p}(\varepsilon\xi)$ is $C_c^{\infty}$. Let $[-N,N]\supset\mathrm{supp}(\psi)$ with $N \in \mathbb{N}$. Then, for $|y|\ge 2N$ and $\varepsilon\in(0,1)$, we find using the Schwartz property of $p$ that 
	\begin{equation}\label{psibound}
		|\psi^\varepsilon(y)| 
		\le C \int_{-N}^N \big|\varepsilon^{-1}p(\varepsilon^{-1}(y-y'))\big|\,\mathrm{d}y' 
		\le C'N\frac{\varepsilon^{-1}}{(1+\varepsilon^{-1}(|y|-N))^2}
		\le C'N\frac{1}{(|y|-N)^2}.
	\end{equation}
	
	Writing $\nu_t^x(b+\mathrm{d}y)$ for the integration with respect to the measure that assigns mass $\nu_t^x(b+A)$ to $A\subset\mathbb{R}$, Lemma \ref{uniformbound} yields, uniformly over $x \in K$, $|b|\le M\sqrt{t}$, and $t \geq 1$,
	\begin{equation}
		\begin{aligned}[t]
			\bigg|\int_{-\infty}^\infty \psi(y-b)-\psi^\varepsilon(y-b)\,\nu_t^x(\mathrm{d}y) \bigg| 
			&= \bigg|\int_{-\infty}^\infty \psi(y)-\psi^\varepsilon(y)\,\nu_t^x(b+\mathrm{d}y) \bigg| \\ 
			&\le 2NC \norm{\psi - \psi^\varepsilon}_{L^{\infty}([-2N,2N])} 
			+ C \sum_{|k|\ge 2N} \sup_{[k,k+1]} |\psi^\varepsilon|.
		\end{aligned}
	\end{equation}
	The last expression goes to $0$ thanks to the third and fourth bounds in \eqref{psibound} and the Dominated Convergence Theorem. Finally, $\norm{\psi^\varepsilon-\psi}_{L^1}\rightarrow 0$ by \cite[Theorem 8.14]{Folland1999-li}.
	
	\medskip
	
	Now, let $\psi$ be bounded, compactly supported, and continuous almost everywhere. At this point, it is enough to demonstrate that for all $\varepsilon>0$, there exist $f,g\in C_c(\mathbb{R})$ with $f\le\psi\le g$ and $\int_{-\infty}^\infty g(y)\,\mathrm{d}y-\int_{-\infty}^\infty f(y)\,\mathrm{d}y<\varepsilon$. To this end, we observe that $\psi$ is Riemann integrable by \cite[Theorem 2.28]{Folland1999-li}, and so there exist compactly supported step functions $f_1$, $g_1$ with $f_1 \le \psi\le g_1$ and $\int_{-\infty}^\infty g_1(y)\,\mathrm{d}y-\int_{-\infty}^\infty f_1(y)\,\mathrm{d}y < \frac\varepsilon 2$. It is then straightforward to construct the desired, e.g., piecewise linear, $f,g\in C_c(\mathbb{R})$. 
	
	\medskip

	Lastly, consider $\psi(y) = e^{-cy}\,\mathbf{1}_{[0,\infty)}(y)$ for $c>0$. Define $\psi_R(y) = e^{-cy}\,\mathbf{1}_{[0,R]}(y)$, for which the limit in \eqref{fourierupgradeeq} holds. We see by Lemma \ref{uniformbound} that, uniformly over $(t,x,b)\in [1,\infty)\times K\times\mathbb{R}$, $\int_{-\infty}^\infty |\psi(y-b) - \psi_R(y-b)|\,\nu_t^x(\mathrm{d}y) \le C\frac{e^{-c\lfloor R\rfloor }}{1-e^{-c}}$. We conclude via $\norm{\psi - \psi_R}_{L^1}= c^{-1}e^{-cR}$. 
\end{proof}

We are now ready for the proof of Theorem \ref{SharpAsymptotics}(a). 

\medskip

\noindent\textit{Proof of Theorem \ref{SharpAsymptotics}(a).}
	Changing the measure to $\mathbb{Q}_{\lambda}^x$ by means of Proposition \ref{SDEfacts}(b), and recalling \eqref{nu def} and $\Lambda'(\lambda)=\beta$, we find that 
	\begin{equation}
		\begin{aligned}[t]
			\mathbb{P}^{\mu}_x\bigg( \int_0^T V(X_t)\,\mathrm{d}t \geq \beta T + y\bigg) &= \frac{e^{-\mu x}\phi_{\lambda}(x)}{\sqrt{T}}\, 
			e^{-T(\lambda \beta - \Lambda_{\mu}(\lambda))}\, 
			\int_y^\infty e^{-\lambda y'}\,\nu_T^x(\mathrm{d}y') \\
			&= \frac{e^{-\mu x}\phi_{\lambda}(x)}{\sqrt{T}}\,
			e^{-T(\lambda \beta - \Lambda_{\mu}(\lambda))}\, e^{-\lambda y}\, \int_{-\infty}^\infty \psi(y'-y)\,\nu_T^x(\mathrm{d}y'),
		\end{aligned}
	\end{equation}
	where $\psi(y'):=e^{-\lambda y'}\,\mathbf{1}_{[0,\infty)}(y')$. It remains to apply Lemma \ref{fourierupgrade}. \qed 

\section{Proof of Theorem \ref{SharpAsymptotics}(b), (c)} \label{sec: sub}

Our proof of Theorem \ref{SharpAsymptotics}(b), (c) is structured in the same way as that of part (a), though the asymptotic analysis in the following analogue of Lemma \ref{supersmoothfourier} is more involved. 

\begin{lemma} \label{CriticalAbounds}
	Let $\mu\ge0$ and $\beta\in(0,\beta_c]$. For $t\ge0$ and $x,\xi\in\mathbb{R}$, define
	\begin{equation}
	A_t^x(\xi) = \frac{1}{e^{-\mu x}\phi_{\lambda_c}(x)}\, e^{ - ti\xi\beta}\, \mathbb{E}^\mu_x\bigg[\exp\bigg((\lambda_c+i\xi)\int_0^t V(X_s)\,\mathrm{d}s\bigg)\bigg].
	\end{equation}
	As before, for any $\widehat{\psi}\in C^\infty_c(\mathbb{R})$ and $b\in\mathbb{R}$, let $\widehat{\psi}_b(\xi):=e^{-i\xi b}\widehat{\psi}(\xi)$. If $\beta = \beta_c$, then for any compact $K\subset \mathbb{R}$ and $M\in(0,\infty)$, as $t\to\infty$,
	\begin{equation}\label{eq: crit asympt}
	\sup_{x \in K}\,\sup_{|b|\leq M\sqrt{t}}\,\bigg|\int_{-\infty}^\infty \widehat{\psi}_b(\xi)\, A_t^x(\xi)\,\mathrm{d}\xi 
	- \frac{2\pi\lambda_c}{\beta_c}\, \langle R_{\mu}^{-1}V,\phi_{\lambda_c}\rangle_{L^2} \,\widehat{\psi}(0)\,\Phi\Big(\!-\frac{b}{\sigma_c \sqrt{t}}\Big) \bigg| \longrightarrow 0,
	\end{equation}
	where $\Phi$ is the standard Gaussian cumulative distribution function and $\sigma_c^2:= \Lambda''(\lambda_c)>0$. If $\beta \in (0,\beta_c)$, then for any compact $K\subset \mathbb{R}$ and $M\in(0,\infty)$, as $t\to\infty$,
	\begin{equation}\label{eq: subcrit asympt}
	\sup_{x \in K}\,\sup_{|b|\leq M\sqrt{t}}\,\bigg|\int_{-\infty}^\infty \widehat{\psi}_b(\xi)\, A_t^x(\xi)\,\mathrm{d}\xi - \frac{2\pi\lambda_c}{\beta_c}\,\langle R_{\mu}^{-1}V,\phi_{\lambda_c}\rangle_{L^2}\,\widehat{\psi}(0)\bigg| \longrightarrow 0.
	\end{equation}
	Moreover, for any $\beta \in (0,\beta_c]$ and compact $K\subset\mathbb{R}$,
	\begin{equation} \label{sec4 unif}
	\sup_{t \geq 1}\,\sup_{x\in K}\,\sup_{b \in \mathbb{R}}\,\bigg|\int_{-\infty}^\infty  \widehat{\psi}_b(\xi)\, A_t^x(\xi)\,\mathrm{d}\xi\bigg|<\infty.
	\end{equation}
\end{lemma}

\begin{proof}
	Let $(-R,R)\supset\mathrm{supp}(\widehat{\psi})$, take $\delta>0$ such that Theorem \ref{DriftSpectrum}(a)--(e) hold on the open ball of radius $\delta$ around $\lambda_c$, and pick a function $g \in C_c^\infty(\mathbb{R})$ with $g\equiv1$ on $\big(\!-\frac{\delta}{2},\frac{\delta}{2}\big)$ and $g\equiv0$ on $\mathbb{R}\setminus(-\delta,\delta)$. We write $z = \lambda_c + i \xi$ and $h_{\lambda_c}(x) = e^{-\mu x}\phi_{\lambda_c}(x)$. By equation \eqref{Duhamel}, with $\lambda=\lambda_c$, we obtain the decomposition 
	\begin{equation} \label{eq: decomp}
		\begin{split}
			h_{\lambda_c}(x)\, A_t^x(\xi) & = e^{-ti\xi\beta}\bigg(1 + z\int_0^t e^{sT^{\mu}_{z}}V(x)\,\mathrm{d}s \bigg) \\
			& = e^{-ti\xi \beta} 
			+e^{-ti\xi \beta}(1-g(\xi))z\int_0^t e^{sT^{\mu}_z}V(x)\,\mathrm{d}s \\
			&\quad + e^{-ti\xi \beta}g(\xi)z\int_0^t e^{sT^{\mu}_z}Q^{\mu}_zV(x)\,\mathrm{d}s 
			+e^{-ti\xi \beta}g(\xi)z\int_0^t e^{s\Lambda_{\mu}(z)}\, P^{\mu}_zV(x)\,\mathrm{d}s.
		\end{split}
	\end{equation}
	We may study $h_{\lambda_c}(x)\, A_t^x(\xi)$ instead of $A_t^x(\xi)$, and so we analyze the terms coming from these four summands in a series of steps.
	
	\medskip
	
	\noindent\textit{Step 1.} For the term coming from the first summand, we integrate by parts and find that
	\begin{equation}
	\sup_{|b| \leq M\sqrt{t}}\,\bigg| \int_{-R}^R e^{-i\xi b}\widehat{\psi}(\xi)\, e^{-ti\xi \beta}\,\mathrm{d}\xi \bigg| \le \sup_{|b| \leq M\sqrt{t}}\,\frac{1}{|\beta t + b|}\,\|\widehat{\psi}'\|_{L^1}\rightarrow 0 \quad\text{as}\quad t \rightarrow \infty.
	\end{equation}
	Moreover, the contribution of this term to the left-hand side of \eqref{sec4 unif} is bounded by $\|\widehat{\psi}\|_{L^1}$.
	
	\medskip
	
	\noindent\textit{Step 2.} Consider the term arising from the second summand on the right-hand side of \eqref{eq: decomp}. Applying Theorem \ref{DriftSpectrum}(e), with $m:=\frac{\delta}{2}$ and $M:=R$ therein, recalling that $\Lambda_\mu(\lambda_c)=0$, and using the bound $|e^{sT_z^{\mu}}V(x)|\le e^{\lambda_c}$ for $0 \leq s \leq 1$, we find that 
	\begin{equation}
	G_t^x(\xi) := \int_0^t e^{sT_z^{\mu}}V(x)\,\mathrm{d}s
	\end{equation}
	is well-defined for all $t\in[0,\infty]$ when $\frac\delta 2\le |\xi|\le R$. We also have  $|G_{\infty}^x(\xi) - G_t^x(\xi)| \le Ce^{-t\varepsilon}$ and $|G_{\infty}^x(\xi)|<\infty$, uniformly over $\xi\in\big[-R,-\frac\delta 2\big]\cup\big[\frac\delta 2,R\big]$, $x\in K$, and $t\in[1,\infty)$. Hence, we can estimate the contribution to the left-hand side of \eqref{sec4 unif} as follows:
	\begin{equation} 
		\begin{split}
	&\, \bigg|\int_{-\infty}^\infty e^{-i\xi b}\widehat{\psi}(\xi)\, e^{-ti\xi\beta}(1-g(\xi))zG_t^x(\xi)\,
	\mathrm{d}\xi\bigg| \\
	&\le C\|\widehat{\psi}\|_{L^1}\, e^{-t\varepsilon}+\bigg|\int_{-\infty}^\infty e^{-i\xi b}\widehat{\psi}(\xi)\,e^{-ti\xi \beta}(1-g(\xi))zG^x_\infty(\xi)\,\mathrm{d}\xi \bigg|,
		\end{split}
	\end{equation}
	which is bounded uniformly over $b \in \mathbb{R}$, $x\in K$, and $t\in[1,\infty)$.
	
	\medskip
	
	We now show that the final integral is negligible in the regime of \eqref{eq: crit asympt}, \eqref{eq: subcrit asympt}. To this end, we control $\frac{\mathrm{d}G_\infty^x}{\mathrm{d}\xi}$ uniformly over $x\in K$. 
	As preparation, we verify that when $\frac{\delta}{2}\leq |\xi| \leq R$, 
	\begin{equation}\label{eq: res iden}
	\int_0^\infty e^{sT_z^{\mu}}f\,\mathrm{d}s = (0-T_z^{\mu})^{-1}f,\quad f\in L^2_\mu. 
	\end{equation}
	Indeed, note that the integral on the left-hand side is well-defined by \eqref{medfreq1} (and its $\mu>0$ analogue -- recall the proof of Theorem \ref{DriftSpectrum}(e)). To see the identity \eqref{eq: res iden}, we start with $f \in \mcal{D}(T_z^\mu)$. Then, $s \mapsto e^{sT_z^\mu}f$ is strongly differentiable with derivative $T_z^{\mu}e^{sT_z^{\mu}}f$, and so
	\begin{equation}
	T_z^\mu \int_0^t e^{sT_z^\mu}f\,\mathrm{d}s = e^{tT_z^\mu}f - f,\quad t\in[0,\infty),\quad f\in \mcal{D}(T_z^\mu).
	\end{equation}
	Due to  \eqref{medfreq1} (and its $\mu>0$ analogue), the norms of $T_z^\mu \int_t^\infty e^{sT_z^\mu}f\,\mathrm{d}s$ and of $e^{tT_z^\mu}f$ decay exponentially as $t\to\infty$, and we see the identity \eqref{eq: res iden} upon applying the closedness of $T_z^\mu$.
	
	\medskip
	
	As further preparation, we claim the identity
	\begin{equation} \label{G infty der}
	\frac{\mathrm{d}}{\mathrm{d}\xi}(0-T_z^{\mu})^{-1} = i(0-T_z^{\mu})^{-1}V(0-T_z^{\mu})^{-1}.
	\end{equation}
	Indeed, using $T_{(\eta)}^{-1}-T_{(\xi)}^{-1}=T_{(\eta)}^{-1}(T_{(\xi)}-T_{(\eta)})T_{(\xi)}^{-1}$, $\eta,\xi\in\mathbb{R}$ where $T_{(\eta)}:=-T^\mu_{\lambda_c+i\eta}$ and $T_{(\xi)}:=-T^\mu_{\lambda_c+i\xi}$, we find that
	\begin{equation}
	(-T^\mu_{\lambda_c+i\eta})^{-1} - (-T^\mu_{\lambda_c+i\xi})^{-1}
	=i(\eta-\xi)(-T^\mu_{\lambda_c+i\eta})^{-1}V (-T^\mu_{\lambda_c+i\xi})^{-1},\quad \eta,\xi\in\mathbb{R}.
	\end{equation}
	Dividing both sides by $\eta-\xi$, taking $\eta\to\xi$, and recalling the continuity of $\eta\mapsto (-T^\mu_{\lambda_c+i\eta})^{-1}$ (see the proof of Lemma \ref{CompactResolvent}), we arrive at the identity \eqref{G infty der}. 
	
	\medskip
	
	We are now ready to bound $\frac{\mathrm{d}G_\infty^x}{\mathrm{d}\xi}$ for $\frac\delta 2\le |\xi|\le R$ and $x\in K$. By Lemma \ref{CompactResolvent} and a compactness argument, the operators $(0-T_z^{\mu})^{-1}:\, L^2_\mu\to H^2_\mu$ are bounded uniformly over $\xi\in\big[-R,-\frac\delta 2\big]\cup\big[\frac\delta 2,R\big]$. Combining this with \eqref{G infty der}, we deduce that $R_\mu^{-1} \frac{\mathrm{d}}{\mathrm{d}\xi}(0-T_z^\mu)^{-1}$, compositions of the form $L^2_\mu\to L^2_\mu\to L^2_\mu\to H^2_\mu\to H^2\hookrightarrow C^1_b$, are also bounded uniformly over $\xi\in\big[-R,-\frac\delta 2\big]\cup\big[\frac\delta 2,R\big]$. Therefore, 
	\begin{equation}\label{Ginftybounds}
		\begin{split}
			\sup_{x\in K}\,\sup_{\frac\delta2\leq |\xi|\leq R}\, \Big|\frac{\mathrm{d}G_\infty^x}{\mathrm{d}\xi}(\xi)\Big| 
			&= \sup_{x\in K}\,\sup_{\frac\delta2\leq |\xi|\leq R}
			\Big|R_\mu R_\mu^{-1}\frac{\mathrm{d}}{\mathrm{d}\xi} (0-T_z^\mu)^{-1}V(x)\Big| \\
			& \le C\sup_{x \in K}\, e^{-\mu x}\,
			\sup_{\frac\delta 2 \le |\xi| \le R}\,
			\Big\|R_\mu^{-1}\frac{\mathrm{d}}{\mathrm{d}\xi}(0-T_z^\mu)^{-1}\Big\|_{\mcal{L}(L^2_{\mu},C^1_b)}.
		\end{split}
	\end{equation}
	Together with the uniform bound on $G_\infty^x(\xi)$ proven above, this shows that the $\xi$-derivative of $\widehat{\psi}(\xi)(1-g(\xi))zG_\infty^x(\xi)$ is bounded uniformly for $\frac\delta 2\le |\xi|\le R$ and $x\in K$. Thus, upon integrating by parts, we find that
	\begin{equation}
	\sup_{|b| \leq M\sqrt{t}}\,\sup_{x \in K}\,\bigg|\int_{-\infty}^\infty e^{-i\xi b}e^{-ti\xi \beta}\,
	\widehat{\psi}(\xi)(1-g(\xi))zG_\infty^x(\xi)\,\mathrm{d}\xi \bigg| \le 
	\sup_{|b| \leq M\sqrt{t}}\,\frac{C}{|\beta t + b|}\rightarrow 0 \quad\text{as}\quad t \rightarrow \infty.
	\end{equation}
	
	\smallskip
	
	\noindent\textit{Step 3.} Step 3 is similar to Step 2. Let $\xi\in[-\delta,\delta]$, and
	\begin{equation}
	F_t^x(\xi) := \int_0^t e^{sT_z^\mu} Q_z^\mu V(x)\,\mathrm{d}s,
	\end{equation}
	which is well-defined for all $t \in [0,\infty]$. Indeed, for $s \in (0,1)$, we note that in $e^{sT_z^\mu}Q_z^\mu V(x) = e^{sT_z^\mu}V(x) - e^{s\Lambda_\mu(z)}P_z^{\mu}V(x)$, the first term can be bounded by $e^{\lambda_c}$ thanks to Theorem \ref{DriftGeneration}(c), and the second term is uniformly bounded (also over $(\xi,x)\in[-\delta,\delta]\times K$ by Theorem \ref{DriftSpectrum}(a), (c)). On the other hand, for $s\in[1,\infty)$, we apply Theorem \ref{DriftSpectrum}(d) to infer that $F_t^x(\xi)$ is well-defined for all $t \in [0,\infty]$, 
	\begin{equation}
	\sup_{x \in K,\,|\xi|\le\delta} |F_\infty^x(\xi) - F_t^x(\xi)| \leq Ce^{-\varepsilon t},\quad t\ge1, 
	\end{equation}
	and $F_\infty^x$ is uniformly bounded for $|\xi|\le\delta$ and $x \in K$. In particular, with respect to \eqref{sec4 unif}, 
	\begin{equation}
	\sup_{t \geq 1}\,\sup_{x\in K}\,\sup_{b \in \mathbb{R}}\,\bigg|\int_{-\infty}^\infty e^{-i\xi b}\widehat{\psi}(\xi)\, e^{-ti\xi\beta}g(\xi)zF_t^x(\xi)\,\mathrm{d}\xi \bigg| < \infty.
	\end{equation}
	
	\smallskip
	
	Next, recall the restricted operator $\widetilde{T}^0_z$ from Lemma \ref{lowerorderresolvent}, and define the analogous $\widetilde{T}_z^\mu=R_\mu(\widetilde{T}^0_z - \frac{\mu^2}{2})R_\mu^{-1}$. Following the proof of \eqref{eq: res iden} and using \eqref{eq: near 0} (and its $\mu>0$ analogue), we show that 
	\begin{equation}\label{F repr}
	\int_0^\infty e^{s T_z^\mu}Q_z^\mu f\,\mathrm{d}s
	=(0-\widetilde{T}_z^{\mu})^{-1}Q_z^{\mu}f,\quad f\in L^2_\mu. 
	\end{equation}
	As in the proof of Lemma \ref{lowerorderresolvent}, we choose a contour $\Gamma$ around $0$ such that
	\begin{equation}\label{Smallcontour}
		(0-\widetilde{T}_z^\mu)^{-1}Q_z^\mu = \frac{1}{2\pi i}\int_{\Gamma} 
		\frac{(\zeta - T_z^\mu)^{-1}Q_z^\mu}{\zeta}\,\mathrm{d}\zeta.
	\end{equation}
	Similarly to Step 2, the norms of $(\zeta - T_z^\mu)^{-1}$ and of its $\xi$-derivative, as operators from the range of $Q_z^\mu$ to $H^2_\mu$, are bounded uniformly over $\zeta\in\Gamma$ and $\xi\in[-\delta,\delta]$. In addition, we recall that $Q_z^\mu=I-P_z^\mu:\,L^2_\mu\to L^2_\mu$ is bounded-holomorphic, which means that 
	both $Q_z^\mu$ and $\frac{\mathrm{d}}{\mathrm{d}\xi} Q_z^\mu$ are uniformly bounded for $|\xi|\le \delta$. Overall, the norm of $\frac{\mathrm{d}}{\mathrm{d}\xi}((0-\widetilde{T}_z^\mu)^{-1}Q_z^\mu):\,L^2_\mu\to H^2_\mu$ is uniformly bounded for $|\xi|\le \delta$. Thus, we obtain an analogue of \eqref{Ginftybounds}, where we replace $\frac\delta 2\le|\xi|\le R$ by $|\xi|\le\delta$, and $G^x_\infty$ by $F^x_\infty$. Integration by parts again yields 
	\begin{equation}
	\sup_{|b| \leq M\sqrt{t}}\,\sup_{x \in K}\,\bigg|\int_{-\infty}^\infty e^{-i\xi b}e^{-ti\xi\beta}\, \widehat{\psi}(\xi)g(\xi)zF_\infty^x(\xi)\,\mathrm{d}\xi \bigg| 
	\le \sup_{|b| \leq M\sqrt{t}}\,\frac{C}{|\beta t + b|}\rightarrow 0 \quad\text{as}\quad t \rightarrow \infty.
	\end{equation}
	
	\smallskip
	
	\noindent\textit{Step 4.} Let us recast the term arising from the fourth summand on the right-hand side of~\eqref{eq: decomp}. Define
	\begin{equation}
	\xi\mapsto a_x(\xi) = \widehat{\psi}(\xi)g(\xi)zP_z^\mu V(x) \in C_c^\infty(\mathbb{R};
	\mathbb{C}).
	\end{equation}
	Also, let $\eta(\xi):=\Lambda_{\mu}(\lambda_c+i\xi) - i \xi \beta_c$, which by Theorem \ref{DriftSpectrum}(a) has the properties $\eta(\xi) = -\frac{1}{2}\sigma_c^2\xi^2 + O(|\xi|^3)$ and $\mathrm{Re}(\eta(\xi)) \le - c\xi^2$ for $|\xi| \le \delta$ (upon decreasing $\delta>0$ if necessary). We seek to understand 
	\begin{equation} \label{relevantintegral}
		\begin{split}
		\int_{-\delta}^\delta \int_0^t e^{-i\xi b}e^{-ti\xi\beta}e^{s\Lambda_{\mu}(z)}\, a_x(\xi) \,\mathrm{d}s \,\mathrm{d}\xi 
		& = \int_{-\delta}^\delta \int_0^t e^{-i\xi b}e^{-(\alpha t -  s)i\xi\beta_c}e^{s\eta(\xi)}\, a_x(\xi) \,\mathrm{d}s \,\mathrm{d}\xi \\
		& = \int_0^t \kappa^x_s((\alpha t - s)\beta_c + b)\,\mathrm{d}s,
		\end{split}
	\end{equation}
	where we have defined $\alpha = \beta/\beta_c$ and
	\begin{equation}\label{Kernel}
		\kappa^x_s(y)= \int_{-\delta}^\delta e^{-i\xi y}e^{s\eta(\xi)}\, a_x(\xi)\,\mathrm{d}\xi = \frac{1}{\sqrt{s}} \int_{-\delta \sqrt{s}}^{\delta \sqrt{s}} e^{-iu y/\sqrt{s}}\, e^{s\eta(u/\sqrt{s})}\, a_x\Big(\frac{u}{\sqrt{s}} \Big)\,\mathrm{d}u,\quad s>0.
	\end{equation}
	
	\smallskip
	
	We first consider
	\begin{equation}\label{01integral}
		\int_{-\delta}^\delta \int_0^1 e^{-i\xi b}e^{-(\alpha t - s)i\xi\beta_c}e^{s\eta(\xi)}\, a_x(\xi) \,\mathrm{d}s \,\mathrm{d}\xi.
	\end{equation}
	The uniform boundedness of the contribution to \eqref{sec4 unif} follows directly from the uniform boundedness of the integrand. Since $\frac{\mathrm{d}}{\mathrm{d}\xi}(e^{si\xi\beta_c}e^{s\eta(\xi)} a_x(\xi))$ is uniformly bounded over $(s,\xi,x)$ in $[0,1]\times[-\delta,\delta]\times K$, integration by parts in $\xi$ again shows that the term in \eqref{01integral} is negligible in the regime of \eqref{eq: crit asympt}, \eqref{eq: subcrit asympt}.

	\medskip
	
	\noindent\textit{Step 5.} Before proceeding further, we establish the estimate
	\begin{equation}\label{KernelBound}
		|\kappa^x_s(y)|\le \frac{C}{\sqrt{s}} \Big(1+\frac{|y|}{\sqrt{s}}\Big)^{-2},\quad (x,s,y)\in K\times[1,\infty)\times\mathbb{R}. 
	\end{equation}
	Let $u\mapsto h^x_s(u) := a_x(u/\sqrt{s})\, e^{s\eta(u/\sqrt{s})}\in C_c^\infty(\mathbb{R})$. Since $a_x$ is bounded and $\mathrm{Re}(\eta(\xi))\leq -c\xi^2$ on the support of $a_x$, we have that $|h^x_s(u)|\le Ce^{-cu^2}$, and so $|\kappa^x_s(y)|\le C/\sqrt{s}$. Moreover,
	\begin{equation}
		\begin{aligned}[t]
			(h^x_s)''(u) =  \frac{1}{s}\, a_x''\Big(\frac{u}{\sqrt{s}}\Big)\, e^{s\eta(u/\sqrt{s})} 
			+2a_x'\Big(\frac{u}{\sqrt{s}}\Big)\,\eta'\Big(\frac{u}{\sqrt{s}}\Big)\,
			e^{s\eta(u/\sqrt{s})} 
		 + a_x\Big(\frac{u}{\sqrt{s}}\Big)\,\eta''\Big(\frac{u}{\sqrt{s}}\Big)\, e^{s\eta(u/\sqrt{s})} \\
		 +\, a_x\Big(\frac{u}{\sqrt{s}}\Big) s  \eta'\Big(\frac{u}{\sqrt{s}}\Big)^2 \, e^{s\eta(u/\sqrt{s})}.
		\end{aligned}
	\end{equation}
	Thanks to the smoothness of $a_x$ and $\eta$, and $\mathrm{supp}(a_x)\subset[-\delta,\delta]$, the first three summands are bounded in absolute value by $Ce^{-cu^2}$ for $s\geq 1$. For the fourth summand, we recall that $\eta$ is twice differentiable with $\eta'(0)=0$, which implies the existence of a constant $C<\infty$ such that $|\eta'(\xi)|\leq C|\xi|$ for $\xi\in\mathrm{supp}(a_x)$. Therefore, the fourth summand is bounded in absolute value by $Cu^2e^{-cu^2}$. With $\theta := y/\sqrt{s}$, we find via integration by parts that 
	\begin{equation}
	|\sqrt{s}\kappa^x_s(y)| = \bigg|\int_{-\delta\sqrt{s}}^{\delta\sqrt{s}} e^{-i\theta u}\, h^x_s(u)\,\mathrm{d}u\bigg| = \frac{1}{\theta^2}\bigg|\int_{-\delta \sqrt{s}}^{\delta\sqrt{s}} e^{-i\theta u}\, (h^x_s)''(u)\,\mathrm{d}u\bigg| \leq \frac{C}{\theta^2},\quad \theta\neq0.
	\end{equation}
	Since we have already shown $|\kappa^x_s(y)| \leq C/\sqrt{s}$, the estimate \eqref{KernelBound} follows from the inequality $\min(1,\theta^{-2})\leq 4(1+|\theta|)^{-2}$.
	
	\medskip
	
	\noindent\textit{Step 6.} We can now complete the proof of \eqref{sec4 unif}. By \eqref{KernelBound}, it holds
	\begin{equation}
		\begin{split}
			&\,\bigg|\int_{-\delta}^\delta \int_1^t e^{-i\xi b}e^{-(\alpha t -  s)i\xi\beta_c}e^{s\eta(\xi)}\, a_x(\xi)\,\mathrm{d}s \,\mathrm{d}\xi\bigg|  = \bigg|\int_1^t \kappa^x_s((\alpha t - s)\beta_c + b)\,\mathrm{d}s\bigg| \\
			&\le C\int_1^t \frac{1}{\sqrt{s}}\bigg(1 + \frac{|(\alpha t-s)\beta_c + b|}{\sqrt{s}} \bigg)^{-2}\,\mathrm{d}s
			\le C'\sup_{r \in \mathbb{R}}\,\int_1^\infty \frac{1}{\sqrt{s}}
			\bigg(1 + \frac{|s - r|}{\sqrt{s}} \bigg)^{-2}\,\mathrm{d}s.
		\end{split}
	\end{equation}
	If $r\le 0$, then $|s-r|\ge s$, and so we can use the estimate
	\begin{equation}
	\int_1^\infty \frac{1}{\sqrt{s}}(1+\sqrt{s})^{-2}\,\mathrm{d}s  \le \int_1^\infty s^{-3/2}\,\mathrm{d}s < \infty.
	\end{equation}
	If $r > 0$, we make the substitution $v=\sqrt{s}$ to rewrite the integral as 
	\begin{equation}
	2\int_1^\infty \bigg(1 + \frac{|v^2-r|}{v}\bigg)^{-2}\,\mathrm{d}v.
	\end{equation}
	Upon ignoring the harmless factor of $2$, we make the substitution $q = j_r(v) := v - r/v$. Note that $j'_r\ge 1$, which leads to
	\begin{equation}
	\int_1^\infty (1+|j_r(v)|)^{-2}\,\mathrm{d}v = \int_{j_r(1)}^\infty (1+|q|)^{-2}\, j'_r(j_r^{-1}(q))^{-1}\,\mathrm{d}q \le \int_{-\infty}^\infty (1+|q|)^{-2}\,\mathrm{d}q < \infty.
	\end{equation}
	
	\smallskip
	
	\noindent\textit{Step 7.} Let $D>M/\beta_c$, and define $r_t = \alpha t + b/\beta_c$ so that
	$(r_t-s)\beta_c=(\alpha t - s)\beta_c+b$. We split the difference between the right-hand side of \eqref{relevantintegral} and the negligible term in \eqref{01integral} according to the value of $D$, and will take $t\to\infty$ followed by $D\to\infty$. First, consider
	\begin{equation}
	\int_1^{r_t - D\sqrt{t}} \kappa^x_s((r_t - s)\beta_c)\,\mathrm{d}s,\quad\text{and if $\alpha<1$, also}\;\int_{r_t + D\sqrt{t}}^t \kappa^x_s((r_t - s)\beta_c)\,\mathrm{d}s.
	\end{equation}
	Using the estimate \eqref{KernelBound}, we find that 
	\begin{equation}
	\bigg| \int_1^{r_t - D\sqrt{t}} \kappa^x_s((r_t-s)\beta_c)\,\mathrm{d}s\bigg| 
	\le C\int_1^{r_t-D\sqrt{t}} \frac{1}{\sqrt{s}}\bigg(1+\frac{r_t-s}{\sqrt{s}}\bigg)^{-2}\,\mathrm{d}s.
	\end{equation}
	For $s\in[1,r_t/2]$, it holds $1 + (r_t-s)/\sqrt{s}\ge r_t/(2\sqrt{s})$, and so for $t\ge0$ large enough,
	\begin{equation}
	\int_1^{\frac{r_t}{2}} \frac{1}{\sqrt{s}}\left(1+\frac{r_t-s}{\sqrt{s}}\right)^{-2}\,\mathrm{d}s
	\le \frac{C}{t^2}\int_1^{\frac{r_t}{2}} \sqrt{s}\,\mathrm{d}s \le \frac{C}{\sqrt{t}}.
	\end{equation}
	For $s\in[r_t/2,r_t - D\sqrt{t}]$, we have $1/\sqrt{s}\in[1/\sqrt{t},C/\sqrt{t}]$ for $t\ge0$ large enough. Combining this with the substitution $q:=(r_t-s)/\sqrt{t}$, we find for $t\ge0$ large enough,
	\begin{equation} \label{C/(1+D) bound}
		\begin{aligned}[t]
	\int_{\frac{r_t}{2}}^{r_t - D\sqrt{t}} \frac{1}{\sqrt{s}}\bigg(1+\frac{r_t-s}{\sqrt{s}}\bigg)^{-2}\, \mathrm{d}s 
	\le \frac{C}{\sqrt{t}}\int_{\frac{r_t}{2}}^{r_t - D\sqrt{t}} 
	\bigg(1 + \frac{r_t-s}{\sqrt{t}}\bigg)^{-2}\,\mathrm{d}s 
	&= C \int_D^{\frac{r_t }{2\sqrt{t}}} (1+q)^{-2}\,\mathrm{d}q \\
	&\le \frac{C}{1+D}.
		\end{aligned}
	\end{equation}
	If $\alpha<1$, then we proceed as in \eqref{C/(1+D) bound} to obtain, for $t\ge0$ large enough,
	\begin{equation}
	\bigg| \int_{r_t + D\sqrt{t}}^t \kappa^x_s((r_t - s)\beta_c)\,\mathrm{d}s\bigg|\leq \frac{C}{1+D}.
	\end{equation}
	
	\smallskip
	
	\noindent\textit{Step 8.} If $\alpha<1$, we are left to study the term
	\begin{equation}
	\int_{-\delta}^\delta \int_{r_t-D\sqrt{t} }^{r_t+D\sqrt{t}} e^{-(r_t-s)i\xi\beta_c} e^{s\eta(\xi)}\, a_x(\xi)\,\mathrm{d}s \,\mathrm{d}\xi =\int_{-\sqrt{t}\delta}^{\sqrt{t}\delta} \int_{-D}^D e^{-i\beta_c uv} e^{(r_t-v\sqrt{t})\eta(\frac{u}{\sqrt{t}})} \, a_x\Big(\frac{u}{\sqrt{t}}\Big)\,\mathrm{d}v \,\mathrm{d}u,
	\end{equation}
	and if $\alpha=1$, the term
	\begin{equation}
	\int_{-\delta}^\delta \int_{r_t-D\sqrt{t}}^t e^{-(r_t - s)i\xi\beta_c}e^{s\eta(\xi)}\,a_x(\xi) \,\mathrm{d}s \,\mathrm{d}\xi
	=\int_{-\sqrt{t}\delta}^{\sqrt{t}\delta} \int_{\frac{b}{\beta_c \sqrt{t}}}^D e^{-i\beta_c uv} e^{(r_t-v\sqrt{t})\eta(\frac{u}{\sqrt{t}})}\, a_x\Big(\frac{u}{\sqrt{t}}\Big)\,\mathrm{d}v \,\mathrm{d}u, 
	\end{equation}
	where we have made the substitution $(u,v):= (\sqrt{t}\xi,(r_t-s)/\sqrt{t})$. Split the domain of integration in $u$ into $R^1_t:=[-t^{1/8}\delta,t^{1/8}\delta]$ and $R^2_t:=(-\sqrt{t}\delta,\sqrt{t}\delta)\setminus R^1_t$. We claim that the contribution from $R_t^2$ is negligible, and we detail only the $\alpha<1$ case since the $\alpha=1$ case is similar. In view of Theorem \ref{DriftSpectrum}(c) and $|b|\le M\sqrt{t}$, we have the bounds $|a_x(u/\sqrt{t})| \le C$, $r_t - v\sqrt{t}\geq c t$, and $\mathrm{Re}(\eta(u/\sqrt{t})) \le -ct^{-3/4}$ for $t\ge0$ large enough. For such $t$, it holds 
	\begin{equation}
	\sup_{x\in K}\,\sup_{|b|\le M\sqrt{t}}\, \bigg|\int_{R^2_t} \int_{-D}^D e^{-i\beta_c uv} e^{(r_t-v\sqrt{t})\eta(\frac{u}{\sqrt{t}})} \, a_x\Big(\frac{u}{\sqrt{t}}\Big)\, \mathrm{d}v \,\mathrm{d}u\bigg| 
	\le C\sqrt{t}De^{-ct^{1/4}}.
	\end{equation}
	
	\smallskip
	
	Turning to the integral over $R^1_t$, we again detail only the $\alpha<1$ case since the $\alpha=1$ case is similar. In that integral, we replace $a_x(u/\sqrt{t})$ with $a_x(0)$ by using
	\begin{equation}
		\begin{split}
	&\,\sup_{x\in K}\,\sup_{|b|\le M\sqrt{t}}\, \bigg|\int_{R^1_t} \int_{-D}^D e^{-i\beta_c uv} e^{(r_t-v\sqrt{t})\eta(\frac{u}{\sqrt{t}})}\,\bigg(a_x\Big(\frac{u}{\sqrt{t}}\Big)-a_x(0)\bigg)\,
	\mathrm{d}v \,\mathrm{d}u\bigg| \\ 
	& \le CD\sup_{x \in K}\,\sup_{|\xi| \le t^{-3/8}\delta}\,|a_x(\xi)-a_x(0)|,
		\end{split}
	\end{equation}
	where we relied on $(r_t-v\sqrt{t})\,\mathrm{Re}(\eta(u/\sqrt{t}))\le -cu^2$ for $t\ge0$ large enough. Next, we seek to replace $(r_t-v\sqrt{t})\eta(u/\sqrt{t})$ with $-\alpha\sigma_c^2 u^2/2$. Recall from Step 4 that $\eta(\xi)+\frac{1}{2}\sigma_c^2\xi^2=O(|\xi|^3)$, and so $|e^{(r_t-v\sqrt{t})\eta(u/\sqrt{t})}-e^{-\frac{\alpha}{2}\sigma_c^2u^2}|\le e^{-\frac{\alpha}{2}\sigma_c^2u^2}|e^{Ct^{-1/8}+Ct^{-1/4}}-1|$ uniformly over $|b|\leq M\sqrt{t}$. Therefore, it holds
	\begin{equation}
		\begin{split}
	&\, \sup_{x \in K}\,\sup_{|b|\leq M\sqrt{t}}
	\,\bigg|a_x(0)\int_{R^1_t} \int_{-D}^D e^{-i\beta_c uv} \big(e^{(r_t-v\sqrt{t})\eta(\frac{u}{\sqrt{t}})} - e^{-\frac{\alpha}{2}\sigma_c^2 u^2}\big)\, \mathrm{d}v \,\mathrm{d}u\bigg| \\ 
	& \leq CD\,\big|e^{Ct^{-1/8}+Ct^{-1/4}}-1\big|.
		\end{split}
	\end{equation}
	Finally, we replace $R^1_t$ with $\mathbb{R}$ by exploiting the integrability of $u\mapsto e^{-\frac{\alpha}{2}\sigma_c^2 u^2}$. 
	
	\medskip
	
	At this point, it remains to analyze 
	\begin{equation}
	a_x(0)\int_{-\infty}^\infty \int_{-D}^D e^{-i\beta_c u v}\, e^{-\frac{\alpha}{2}\sigma_c^2 u^2} \,\mathrm{d}v \,\mathrm{d}u
	=a_x(0)\,\frac{\sqrt{2\pi}}{\sqrt{\alpha}\sigma_c} \int_{-D}^D 
	e^{-\frac{\beta_c^2 v^2}{2\alpha\sigma_c^2}}\,\mathrm{d}v,
	\end{equation} 
	if $\alpha<1$, and
	\begin{equation}
	a_x(0)\int_{-\infty}^\infty \int_{\frac{b}{\beta_c \sqrt{t}}}^D 
	e^{-i\beta_c u v}\, e^{-\frac{1}{2}\sigma_c^2 u^2} \,\mathrm{d}v \,\mathrm{d}u =a_x(0)\,\frac{\sqrt{2\pi}}{\sigma_c} \int_{\frac{b}{\beta_c \sqrt{t}}}^D 
	e^{-\frac{\beta_c^2 v^2}{2\sigma_c^2}}\,\mathrm{d}v,
	\end{equation}
	if $\alpha=1$. These terms admit the respective approximations 
	\begin{equation}
	a_x(0)\,\frac{\sqrt{2\pi}}{\sqrt{\alpha}\sigma_c} \int_{-\infty}^\infty 
	e^{-\frac{\beta_c^2 v^2}{2\alpha\sigma_c^2}}\,\mathrm{d}v 
	= a_x(0)\,\frac{2\pi}{\beta_c},
	\end{equation}
	if $\alpha < 1$, and
	\begin{equation}
	a_x(0)\,\frac{\sqrt{2\pi}}{\sigma_c} \int_{\frac{b}{\beta_c \sqrt{t}}}^\infty 
	e^{-\frac{\beta_c^2 v^2}{2\sigma_c^2}}\,\mathrm{d}v 
	= a_x(0)\,\frac{2\pi}{\beta_c}\,\Phi\Big(\! -\frac{b}{\sigma_c\sqrt{t}}\Big),
	\end{equation}
	if $\alpha = 1$, up to an error $\iota(D)\to0$ as $D\to\infty$, uniformly over $x\in K$ and $|b|\le M\sqrt{t}$.

	\medskip

	In summary, we have shown that for any $D>M/\beta_c$,
	\begin{equation}
	\limsup_{t\to\infty}\,\sup_{x \in K}\,\sup_{|b|\le M\sqrt{t}}\,
	\bigg|\int_{-\infty}^\infty \widehat{\psi}_b(\xi)\, A_t^x(\xi)\,\mathrm{d}\xi - h_{\lambda_c}(x)^{-1}\, a_x(0)\,\frac{2\pi}{\beta_c}\bigg|
	\le \frac{C}{1+D}+\iota(D),
	\end{equation}
	if $\alpha<1$, and 
	\begin{equation}
	\limsup_{t\to\infty}\,\sup_{x \in K}\,\sup_{|b|\le M\sqrt{t}}\,
	\bigg|\int_{-\infty}^\infty \widehat{\psi}_b(\xi)\, A_t^x(\xi)\,\mathrm{d}\xi - h_{\lambda_c}(x)^{-1}\, a_x(0)\,\frac{2\pi}{\beta_c}\,\Phi\Big(\!-\frac{b}{\sigma_c\sqrt{t}}\Big)\bigg|
	\le \frac{C}{1+D}+\iota(D),
	\end{equation}
	if $\alpha = 1$. We conclude by taking $D\to\infty$ and using $a_x(0) = \widehat{\psi}(0)\lambda_c\langle R_{\mu}^{-1}V, \phi_{\lambda_c}\rangle_{L^2} h_{\lambda_c}(x)$.
\end{proof}

\begin{lemma}\label{subcritvague}
	Let $\mu\ge0$ and $\beta\in(0,\beta_c]$. Fix a compact $K\subset\mathbb{R}$ and $M>0$, and with $Z_t:=\int_0^t V(X_s)-\beta\,\mathrm{d}s$, define the measures
	\begin{equation}
	\nu_t^x(\cdot) = \mathbb{E}^{\mathbb{Q}^x_{\lambda_c}}
	\Big[\big(e^{-\mu X_t}\,\phi_{\lambda_c}(X_t)\big)^{-1}\,\mathbf{1}_{\{Z_t\in\,\cdot\,\}}\Big].
	\end{equation}
	Then, 
	\begin{equation}
	\sup_{t \geq 1}\,\sup_{x \in K}\,\sup_{b \in \mathbb{R}}\,\nu_t^x([b,b+1])<\infty.
	\end{equation}
	Moreover, for any bounded and compactly supported $\psi:\,\mathbb{R}\to\mathbb{R}$ that is continuous almost everywhere, it holds
	\begin{equation} \label{eq: crit asymp}
	\lim_{t\to\infty}\,\sup_{x\in K}\,\sup_{|b|\le M\sqrt{t}}\,
	\bigg|\int_{-\infty}^\infty \psi(y-b)\,\nu_t^x(\mathrm{d}y) - \frac{\lambda_c}{\beta_c}\,\Phi\Big(\!-\frac{b}{\sigma_c\sqrt{t}}\Big)\,
	\langle R_{\mu}^{-1}V,\phi_{\lambda_c}\rangle_{L^2} 
	\int_{-\infty}^\infty \psi(y)\,\mathrm{d}y\bigg| = 0,
	\end{equation}
	if $\beta = \beta_c$, and
	\begin{equation} \label{eq: subcrit asymp}
	\lim_{t\to\infty}\,\sup_{x\in K}\,\sup_{|b|\le M\sqrt{t}}\, 
	\bigg|\int_{-\infty}^\infty \psi(y-b)\,\nu_t^x(\mathrm{d}y) - 
	\frac{\lambda_c}{\beta_c}\,\langle R_{\mu}^{-1}V,\phi_{\lambda_c}\rangle_{L^2}
	\int_{-\infty}^\infty \psi(y)\,\mathrm{d}y\bigg| = 0,
	\end{equation}
	if $\beta\in(0,\beta_c)$. The limits in \eqref{eq: crit asymp} and \eqref{eq: subcrit asymp} also apply to $\psi(y)=e^{-cy}\,\mathbf{1}_{[0,\infty)}(y)$ for $c>0$.
\end{lemma}

\begin{proof}
	It suffices to follow the proofs of Lemmas \ref{uniformbound} and \ref{fourierupgrade}, and to use Lemma \ref{CriticalAbounds} instead of Lemma \ref{supersmoothfourier}.
\end{proof}

We are now ready to prove Theorem \ref{SharpAsymptotics}(b), (c). 

\medskip

\noindent\textit{Proof of Theorem \ref{SharpAsymptotics}(b), (c).} By a change of measure to $\mathbb{Q}_{\lambda_c}^x$ (see Proposition~\ref{SDEfacts}(b)),
\begin{equation}
	\begin{split}
		\mathbb{P}^{\mu}_x\bigg(\int_0^T V(X_t)\,\mathrm{d}t \ge \beta T + y\bigg) &
		= {e^{-\mu x}\phi_{\lambda_c}(x)}\, e^{-T\lambda_c \beta}\,
		\int_y^\infty e^{-\lambda_c y'}\,\nu_T^x(\mathrm{d}y') \\
		&= e^{-\mu x}\phi_{\lambda_c}(x)\, e^{-T\lambda_c \beta}\, 
		e^{-\lambda_c y}\,\int_{-\infty}^\infty \psi(y'-y)\,\nu_T^x(\mathrm{d}y'),
	\end{split}
\end{equation}
where $\psi(y') := e^{-\lambda_cy'}\,\mathbf{1}_{[0,\infty)}(y')$. It remains to apply Lemma \ref{subcritvague}. \qed

\section{Proofs of Theorems \ref{MainTheorem} and \ref{Macroscopic}} \label{sec: main}

In this section, we deduce the statements in Theorems \ref{MainTheorem}, \ref{Macroscopic} from Theorem \ref{SharpAsymptotics} in the following order: Theorem \ref{MainTheorem}(a), Theorem \ref{Macroscopic}(a), Theorem \ref{MainTheorem}(b), and Theorem \ref{Macroscopic}(b). 

\medskip

\noindent\textit{Proof of Theorem \ref{MainTheorem}(a).} We start with the case $\alpha = 0$. Let $L_t^T := \int_t^T V(X_s)\,\mathrm{d}s$. By Theorem \ref{SharpAsymptotics}, for any $x,x',y\in\mathbb{R}$,
\begin{equation} \label{problimit}
\lim_{T\to\infty} \frac{\mathbb{P}^\mu_{x'}(L_0^{T-S} \ge \beta (T-S) + \beta S -y)} {\mathbb{P}^\mu_x(L_0^T\geq \beta T)} = e^{-\mu(x'-x)}\,\frac{\phi_\lambda(x')}{\phi_\lambda(x)}\,\exp\Big(S\Big(\frac{\mu^2}{2} - \Lambda(\lambda)\Big) + \lambda y\Big).
\end{equation}
Let $E\subset C([0,S])$ be Borel and let $(\mcal{F}_t)_{t \ge 0}$ denote the canonical filtration.  Then, 
\begin{equation}
\mathbb{P}_x^{\mu}(E | L_0^T \ge \beta T) = \frac{\mathbb{E}_x^{\mu}[\mathbf{1}_E \, \mathbf{1}_{\{L_0^T \ge \beta T\}}]}{\mathbb{P}_x^{\mu}\left(L_0^T \geq \beta T\right)} 
= \mathbb{E}_x^\mu\bigg[ \mathbf{1}_E \,\frac{\mathbb{P}_x^\mu(L_S^T \ge  \beta (T-S) + \beta S -  L_0^S | \mcal{F}_S)}{\mathbb{P}_x^{\mu}(L_0^T \geq \beta T)}\bigg].
\end{equation}
Together with Fatou's Lemma, the Markov property of $X$ under $\mathbb{P}^\mu_x$, and \eqref{problimit}, this yields
\begin{equation}
\liminf_{T \to \infty}\,\mathbb{P}_x^\mu(E | L_0^T \ge \beta T) \ge \mathbb{E}_x^{\mu}\bigg[ \mathbf{1}_E \, e^{-\mu(X_S - x)}\, \frac{\phi_\lambda(X_S)}{\phi_\lambda(x)}\,
\exp\Big(S\Big(\frac{\mu^2}{2}-\Lambda(\lambda)\Big) + \lambda L_0^S\Big)\bigg] = \mathbb{Q}_{\lambda}^x(E),
\end{equation}
where the equality follows from \eqref{GirsanovDensity}. We conclude by the Portmanteau Theorem.

\medskip

We move on to the case $\alpha \in (0,1)$. Let $E'\subset C([0,S])$ be Borel and define $E = \{X_{[\alpha T,\alpha T + S]} \in E' \}$. Set $r_T = \alpha T + S$ and $p_x(T,y)= \mathbb{P}^\mu_x(L_0^T \geq \beta T + y)$. Write
$E_N = E \cap \{X_{r_T} \in [-N,N]\}$. Take $g_N:\,\mathbb{R} \to [0,1]$ to be a $C_c^\infty$ function with $\mathrm{supp}(g_N)\subset[-N,N]$ and $g_N\equiv 1$ on $[-N/2,N/2]$. Finally, recall that $Z_t = L_0^t - \beta t$ and $G_t = Z_t/\sqrt{t}$. To analyze $\mathbb{E}_x^\mu[\mathbf{1}_E \,\mathbf{1}_{\{L_0^T \geq \beta T \}}]$, we use the following chain of inequalities based on \eqref{GirsanovDensity}:
\begin{equation}\label{ineqchain}
\begin{aligned}
&\;\mathbb{E}_x^\mu[\mathbf{1}_E \,\mathbf{1}_{\{L_0^T \geq \beta T\}}] = \mathbb{E}_x^\mu[\mathbf{1}_E \, p_{X_{r_T}}(T-r_T,-Z_{r_T})] \\
&= h_\lambda(x)\, e^{r_T\Lambda_\mu(\lambda) - \lambda\beta r_T}\, \mathbb{E}^{\mathbb{Q}_{\lambda}^x}\big[\exp(-\lambda Z_{r_T})\, h_\lambda(X_{r_T})^{-1}\,\mathbf{1}_E \, p_{X_{r_T}}(T-r_T,-Z_{r_T})\big] \\
&\ge  h_\lambda(x)\, e^{r_T\Lambda_\mu(\lambda) - \lambda\beta r_T}\, \mathbb{E}^{\mathbb{Q}_\lambda^x}\big[\mathbf{1}_{E_N}\, g_N(G_{r_T})\,
\exp(-\lambda Z_{r_T})\, h_\lambda(X_{r_T})^{-1}\, p_{X_{r_T}}(T-r_T,-Z_{r_T})\big] \\
&\ge \frac{\sqrt{T}}{\sqrt{T-r_T}}\,a(\kappa N,T-r_T)\, b_x(T)\, \mathbb{E}^{\mathbb{Q}_\lambda^x}\Big[\mathbf{1}_{E_N}\, g_N(G_{r_T})\,
\exp\Big(\!-\frac{G_{r_T}^2}{2 \Lambda''(\lambda)}\,\frac{r_T}{T-r_T}\Big)\Big] \\
&\ge \frac{\sqrt{T}}{\sqrt{T-r_T}}\,
a(\kappa N,T - r_T)\, b_x(T)\, c\Big(\alpha,\frac{S}{T},N\Big)\, \mathbb{E}^{\mathbb{Q}_\lambda^x}[\mathbf{1}_{E_N}\, g_N(G_{r_T})\,
\exp\Big(\!-\frac{G_{r_T}^2}{2 \Lambda''(\lambda)}\,\frac{\alpha}{1-\alpha}\Big)\Big],
\end{aligned}
\end{equation}
where $\kappa\ge1$ ensures $|Z_{r_T}|\le\kappa N\sqrt{T-r_T}$ whenever $g_N(G_{r_T})>0$ and $T\ge0$ is large~enough, and we have defined 
\begin{eqnarray}
&& b_x(T) = h_\lambda(x)\, 
\frac{e^{-T(\lambda \beta-\Lambda_{\mu}(\lambda))}} {\Lambda_{\mu}(\lambda)\sqrt{2\pi T\Lambda''(\lambda)}}\,\int_{-1}^1\phi_{\lambda}(r)e^{\mu r}\,\mathrm{d}r, \\
&& a(N,T) = \inf_{|x|\leq N,\,|y| \leq N\sqrt{T}} \Big(p_x(T,y)\, b_x(T)^{-1}\,
e^{\lambda y + \frac{y^2}{2\Lambda''(\lambda)T}}\Big), \\
&& c(\alpha,\iota,N)=\exp\Big(\!-\frac{N^2}{2\Lambda''(\lambda)}\,
\Big|\frac{\alpha + \iota}{1-\alpha - \iota}-\frac{\alpha}{1-\alpha} \Big|\Big).
\end{eqnarray}

\smallskip

By Theorem \ref{SharpAsymptotics}(a), $a(\kappa N,T-r_T)$ and $b_x(T)/\mathbb{P}_x^\mu(L_0^T \ge \beta T)$ converge to $1$ as $T \to \infty$. Clearly, $c(\alpha,S/T,N)$ also converges to $1$ as $T \to \infty$. Combining these with \eqref{ineqchain}, we obtain
\begin{equation} \label{two liminfs}
\liminf_{T\to\infty}\,\mathbb{P}_x^\mu(E | L_0^T \ge \beta T) \ge 
\frac{1}{\sqrt{1-\alpha}}\,\liminf_{T \to\infty}\, \mathbb{E}^{\mathbb{Q}_\lambda^x}\Big[\mathbf{1}_{E_N}\, g_N(G_{r_T})\,
\exp\Big(\!-\frac{G_{r_T}^2}{2 \Lambda''(\lambda)}\,\frac{\alpha}{1-\alpha}\Big)\Big].
\end{equation}
Since $|G_{r_T}-G_{\alpha T}|\le C/\sqrt{T}$, the latter limit inferior is unchanged if we replace $r_T$ by $\alpha T$. In addition, the Markov property of $X$ under $\mathbb{Q}_\lambda^x$ (cf.~Proposition \ref{SDEfacts}(a)) yields 
\begin{equation}\label{supercritexpectation}
\mathbb{E}^{\mathbb{Q}_\lambda^x}\Big[\mathbf{1}_{E_N} g_N(G_{\alpha T})
\exp\Big(\!-\frac{G_{\alpha T}^2}{2 \Lambda''(\lambda)}\,\frac{\alpha}{1\!-\!\alpha}\Big)\Big] = \mathbb{E}^{\mathbb{Q}_\lambda^x}\Big[e_N(X_{\alpha T})\, g_N(G_{\alpha T})
\exp\Big(\!-\frac{G_{\alpha T}^2}{2 \Lambda''(\lambda)}\frac{\alpha}{1\!-\!\alpha}\Big)\Big],
\end{equation}
where $e_N(x) := \mathbb{Q}_\lambda^x\big(X_{[0,S]}\in E',\, X_S\in[-N,N]\big)$ is a bounded measurable function. Thus, by Proposition \ref{SDEfacts}(c), the right-hand side of \eqref{two liminfs} evaluates to 
\begin{equation}
\frac{1}{\sqrt{1-\alpha}}\,
\mathbb{Q_\lambda^{\pi_\lambda}}\big(X_{[0,S]}\in E',\, X_S\in[-N,N]\big)\,
\int_{-\infty}^\infty g_N(y)\, e^{-\frac{y^2}{2\Lambda''(\lambda)}\,\frac{\alpha}{1-\alpha}}\,
\frac{e^{-\frac{y^2}{2\Lambda''(\lambda)}}}{\sqrt{2\pi\Lambda''(\lambda)}}\,\mathrm{d}y.
\end{equation}
Finally, we take the $N\to\infty$ limit (obtaining a lower bound on the left-hand side of \eqref{two liminfs}). The probability converges to $\mathbb{Q}_{\lambda}^{\pi_\lambda}(X_{[0,S]}\in E')$, and by the dominated convergence theorem,
\begin{equation} \label{easy DCT}
\int_{-\infty}^\infty g_N(y)\, e^{-\frac{y^2}{2\Lambda''(\lambda)}\,\frac{\alpha}{1-\alpha}}\,
\frac{e^{-\frac{y^2}{2\Lambda''(\lambda)}}}{\sqrt{2\pi\Lambda''(\lambda)}}\,\mathrm{d}y \to \int_{-\infty}^\infty e^{-\frac{y^2}{2\Lambda''(\lambda)}\,\frac{\alpha}{1-\alpha}}\,
\frac{e^{-\frac{y^2}{2\Lambda''(\lambda)}}}{\sqrt{2\pi\Lambda''(\lambda)}}\,\mathrm{d}y = \sqrt{1-\alpha}.
\end{equation}
It follows that
\begin{equation}
\liminf_{T \to \infty}\,\mathbb{P}_x^\mu(E | L_0^T \ge \beta T) \ge \mathbb{Q}_\lambda^{\pi_{\lambda}}(X_{[0,S]}\in E'),
\end{equation}
and we conclude again by the Portmanteau Theorem. \qed

\medskip

\noindent\textit{Proof of Theorem \ref{Macroscopic}(a).} We start with the case $\alpha \in (0,1)$. Let $E_T\in\mathcal{F}_{\alpha T}$, $T\ge0$ be any events satisfying $\mathbb{Q}_\lambda^x(E_T)\to 1$ as $T\to\infty$. Arguing as in the proof of \eqref{two liminfs}, we lower bound $\mathbb{P}_x^{\mu}({E_T}|{L_0^T\geq \beta T})$ by 
\begin{equation}\label{empmeasureineqchain}
\begin{aligned}
&\,\frac{\sqrt{T}}{\sqrt{T-\alpha T}}\,(1-o(1))\, \mathbb{E}^{\mathbb{Q}_\lambda^x}\Big[\mathbf{1}_{E_T}\,\mathbf{1}_{\{X_{\alpha T}\in[-N,N]\}}\,
g_N(G_{\alpha T})\,
\exp\Big(\!-\frac{G_{\alpha T}^2}{2\Lambda''(\lambda)}\frac{\alpha}{1-\alpha}\Big)\Big] \\
&\ge\frac{1}{\sqrt{1-\alpha}}\,(1-o(1))\,
\bigg(\mathbb{E}^{\mathbb{Q}_\lambda^x}\Big[\mathbf{1}_{\{X_{\alpha T}\in[-N,N]\}}\,
g_N(G_{\alpha T})\,
\exp\Big(\!-\frac{G_{\alpha T}^2}{2\Lambda''(\lambda)}\frac{\alpha}{1-\alpha}\Big)\Big] 
- \mathbb{Q}^x_\lambda(E_T^c) \!\bigg).
\end{aligned}
\end{equation}
Taking $T\to\infty$ via Proposition \ref{SDEfacts}(c) together with $\mathbb{Q}_\lambda^x(E_T)\to 1$, and subsequently letting $N\to\infty$ as in \eqref{easy DCT}, we arrive at a lower bound of $1$. Hence, to obtain Theorem \ref{Macroscopic}(a), it suffices to define $E_T$ as the event that $\rho_{\alpha T}$ belongs to a neighborhood of $\pi_\lambda$ and to apply Proposition \ref{SDEfacts}(d). For the case $\alpha =1$, fix an $\varepsilon\in(0,1)$ and set $\alpha'=1-\varepsilon/2$. Since the L\'{e}vy-Prokhorov distance between $\rho_{\alpha' T}$ and $\rho_T$ is bounded above by $\varepsilon/2$, it holds 
\begin{equation}
\liminf_{T\to\infty}\,\mathbb{P}_x^\mu(\rho_T \in B_{\varepsilon}(\pi_\lambda)|L^T_0\ge\beta T) \ge \liminf_{T\to\infty}\,\mathbb{P}_x^\mu\big(\rho_{\alpha 'T} \in B_{\frac \varepsilon 2}(\pi_\lambda)\big|L^T_0\ge\beta T\big)\ge1, 
\end{equation}
where $B_{\varepsilon}(\pi_\lambda)$ is the L\'{e}vy-Prokhorov open $\varepsilon$-ball around $\pi_\lambda$. \qed

\medskip

\noindent\textit{Proof of Theorem \ref{MainTheorem}(b), Remark \ref{rem: stronger}(a).}~Take an event of the form $E = \{X_{[\alpha T,\alpha T+S]} \in E'\}$. If $\alpha=0$ or $\beta=\beta_c$, then we can repeat the proof of Theorem \ref{MainTheorem}(a), only relying on Theorem \ref{SharpAsymptotics}(b), (c) instead of Theorem \ref{SharpAsymptotics}(a). Hence, we turn to the case $\beta\in(0,\beta_c)$ and $0<\alpha<\frac{\beta}{\beta_c}$. Let $r_T = \alpha T + S$ and $p_x(T,y,\upsilon)=\mathbb{P}^\mu_x(L_0^T\ge \upsilon T + y)$. As before, write $E_N=E\cap\{X_{r_T}\in[-N,N]\}$, take $g_N:\,\mathbb{R} \to [0,1]$ to be a $C_c^\infty$ function with $\mathrm{supp}(g_N)\subset[-N,N]$ and $g_N\equiv 1$ on $[-N/2,N/2]$, and set $Z_t=L_0^t-\beta_c t$ and $G_t=Z_t/\sqrt{t}$. Define further $\chi_T = \beta_c\,\frac{(\beta/\beta_c)T-r_T}{T-r_T}$, $\chi= \beta_c\,\frac{\beta/\beta_c-\alpha}{1-\alpha}<\beta_c$, $\varepsilon_T = \chi - \chi_T=O(T^{-1})$, and $\Theta_T = -Z_{r_T}-(T-r_T)\varepsilon_T$. 
We employ the following chain of inequalities: 
\begin{equation} \label{subcritexpectation}
\begin{aligned}
&\;\mathbb{E}_x^\mu[\mathbf{1}_E \,\mathbf{1}_{\{L_0^T \ge \beta T \}}] = \mathbb{E}_x^\mu[\mathbf{1}_E \, p_{X_{r_T}}(T-r_T,-Z_{r_T},\chi_T)] \\
&= h_{\lambda_c}(x)\, e^{-\lambda_c \beta_c r_T}\, \mathbb{E}^{\mathbb{Q}_{\lambda_c}^x}\big[\exp(-\lambda_cZ_{r_T})\, h_{\lambda_c}(X_{r_T})^{-1}\,
\mathbf{1}_E \, p_{X_{r_T}}(T-r_T,-Z_{r_T}-(T-r_T)\varepsilon_T,\chi)\big] \\
&\geq h_{\lambda_c}(x)\, e^{-\lambda_c \beta_c r_T}\, \mathbb{E}^{\mathbb{Q}_{\lambda_c}^x}\big[\mathbf{1}_{E_N}\, g_N(G_{r_T})\, \exp(-\lambda_cZ_{r_T})\, h_{\lambda_c}(X_{r_T})^{-1}\, p_{X_{r_T}}(T-r_T,\Theta_T,\chi)\big] \\
& \ge b_x(T,\beta)\, a(\kappa N,T-r_T,\chi)\, \mathbb{E}^{\mathbb{Q}_{\lambda_c}^x}[\mathbf{1}_{E_N}\, g_N(G_{r_T})],
\end{aligned}
\end{equation}
where $\kappa \geq 1$ ensures $\max(|Z_{r_T}+\beta_c S|,|\Theta_T|)\le \kappa N \sqrt{T-r_T}$ whenever $g_N(G_{r_T})>0$ and $T\ge0$ is large enough, and we have defined
\begin{equation}
b_x(T,\upsilon) = h_{\lambda_c}(x)\, 
\frac{e^{-T\lambda_c\upsilon}}{\beta_c}\, 
\int_{-1}^1\phi_{\lambda_c}(r)\, e^{\mu r}\,\mathrm{d}r,\quad
a(N,T,\upsilon) = \inf_{|x|\leq N,\,|y|\leq N\sqrt{T}}\, 
\frac{p_x(T,y,\upsilon)\, e^{\lambda_c y}}{b_x(T,\upsilon)}.
\end{equation}

\smallskip

By Theorem \ref{SharpAsymptotics}(c), $b_x(T,\beta)/{\mathbb P}^\mu_x(L^T_0\ge\beta T)$ and $a(\kappa N,T-r_T,\chi)$ converge to $1$ as $T\to\infty$ (note that $\beta,\chi\in(0,\beta_c)$). Moreover, arguing as after \eqref{two liminfs} in the proof of Theorem \ref{MainTheorem}(a), we find that $\mathbb{E}^{\mathbb{Q}_{\lambda_c}^x}[\mathbf{1}_{E_N}\, g_N(G_{r_T})]$ tends to $\mathbb{Q}_{\lambda_c}^{\pi_{\lambda_c}}(X_{[0,S]}\in E')$ when we take $T\to\infty$ and then $N\to\infty$. Consequently, $\liminf_{T\to\infty}\,\mathbb{P}_x^\mu(E|L_0^T \ge \beta T)\ge \mathbb{Q}_{\lambda_c}^{\pi_{\lambda_c}}(X_{[0,S]}\in E')$, and we conclude again by the Portmanteau Theorem.  

\medskip

We proceed to the proof of Remark \ref{rem: stronger}(a). Set $\tau_T = \inf\{t \ge 0:\, L_0^t \ge \beta T\}$. Then, $\{\tau_T \le \gamma T\} = \{L_0^{\gamma T} \ge \beta T \}$ for any $\gamma>0$. For any $\gamma\in\big(\frac{\beta}{\beta_c},\alpha_2\big)$, it holds therefore $\{\tau_T\le \gamma T\} \subset \{ L_0^T \ge \beta T\}$. By Theorem \ref{SharpAsymptotics}(c), as $T\to\infty$, 
\begin{equation} \label{event approx}
\mathbb{P}_x^{\mu}(\tau_T\le \gamma T | L_0^T \geq \beta T) \rightarrow 1.
\end{equation}
Let $X^{(1)}_T := X_{[\alpha_1T,\alpha_1 T+S_1]}$ and $X^{(2)}_T := (X_{\alpha_2 T + t} - X_{\alpha_2 T})_{0 \leq t \leq S_2}$. Then, for Borel $E^{(1)}$, $E^{(2)}$:
\begin{equation}
\begin{aligned}
&\;\mathbb{P}^\mu_x(X^{(1)}_T \in E^{(1)}, X^{(2)}_T \in E^{(2)} | L_0^T \ge \beta T) \\
&= \mathbb{P}^\mu_x(X^{(1)}_T \in E^{(1)}, X^{(2)}_T \in E^{(2)} | L_0^T \geq \beta T,\tau_T\le \gamma T) + o(1) \\
&= \mathbb{P}_x^\mu(X^{(2)}_T \in E^{(2)} | X^{(1)}_T \in E^{(1)}, \tau_T\le \gamma T)\, \mathbb{P}_x^\mu(X^{(1)}_T \in E^{(1)} | \tau_T\le \gamma T) + o(1).
\end{aligned}
\end{equation}
Since $\{X^{(1)}_T \in E^{(1)}, \tau_T\le \gamma T\}\in\mcal{F}_{\gamma T}$ for $T\ge0$ large enough and $X$ has independent stationary increments, the first probability in the last line equals $\mathbb{P}_x^\mu((X_t-X_0)_{0 \le t \le S_2} \!\in\! E^{(2)})$. For the second probability, we observe that conditioning on $\{\tau_T\le \gamma T\}$ is the same as conditioning on $\big\{L_0^{\gamma T}\ge \frac{\beta}{\gamma}\gamma T\big\}$. Now, $\frac{\beta}{\gamma}\in(\alpha_1\beta_c,\beta_c)$ by construction, and so the previous parts of this proof imply that $\mathbb{P}_x^\mu(X^{(1)}_T \in E^{(1)} | \tau_T\le \gamma T)\to \mathbb{Q}^\pi_{\lambda_c}(X_{[0,S_1]} \in E^{(1)})$, with $\pi=\delta_x$ if $\alpha_1=0$, and $\pi=\pi_{\lambda_c}$ if $\alpha_1>0$. We conclude via the Portmanteau Theorem. \qed

\medskip

\noindent\textit{Proof of Theorem \ref{Macroscopic}(b).} The result for $\alpha\in\big(0,\frac{\beta}{\beta_c}\big]$ can be inferred from the proof of Theorem \ref{MainTheorem}(b) in the same way we deduced Theorem \ref{Macroscopic}(a) from the proof of Theorem \ref{MainTheorem}(a). Hence, we turn to the case $\alpha\in\big(\frac{\beta}{\beta_c},1\big]$. Note that
\begin{equation}
\rho_{\alpha T} = \frac{\beta}{\alpha\beta_c}\rho_{(\beta/\beta_c)T} 
+ \frac{1}{\alpha T}\int_{\frac{\beta}{\beta_c}T}^{\alpha T} \mathbf{1}_{\cdot}(X_t)\,\mathrm{d}t.
\end{equation}
Conditioned on $\{L_0^T \geq \beta T \}$, the first term on the right-hand side converges weakly in probability to $\frac{\beta}{\alpha \beta_c}\pi_{\lambda_c}$, and so it suffices to show that for all $N \in \mathbb{N}$ and $\varepsilon_0\in(0,1)$,
\begin{equation}
\mathbb{P}_x^\mu\bigg(\frac{1}{T}\int_{\frac{\beta}{\beta_c} T}^{\alpha T} \mathbf{1}_{[-N,N]}(X_t)\,\mathrm{d}t > \varepsilon_0 \bigg | L_0^T \ge \beta T\bigg) \to 0.
\end{equation}

\smallskip

Pick $\varepsilon\in(0,\varepsilon_0/2)$ such that $\gamma :=\beta/\beta_c + \varepsilon 
< \alpha$, and recall \eqref{event approx}. On the event $\{\tau_T \le \gamma T\}$, 
\begin{equation}
\frac{1}{T}\int_{\frac{\beta}{\beta_c} T}^{\alpha T} \mathbf{1}_{[-N,N]}(X_t)\,\mathrm{d}t \le \frac{\varepsilon_0}{2} +\frac{1}{T}\int_{\tau_T}^\infty \mathbf{1}_{[-N,N]}(X_t)\,\mathrm{d}t.
\end{equation}
We conclude that
\begin{equation}
\begin{split}
&\;\mathbb{P}_x^\mu\bigg(\frac{1}{T}\int_{\frac{\beta}{\beta_c} T}^{\alpha T} \mathbf{1}_{[-N,N]}(X_t)\,\mathrm{d}t > \varepsilon_0 \bigg | L_0^T \ge \beta T\bigg) \\
&= \mathbb{P}_x^\mu\bigg(\frac{1}{T}\int_{\frac{\beta}{\beta_c} T}^{\alpha T} \mathbf{1}_{[-N,N]}(X_t)\,\mathrm{d}t > \varepsilon_0 \bigg | \tau_T \le \gamma T\bigg) + o(1)\\
&\le  \mathbb{P}_x^\mu\bigg(\frac{1}{T}\int_{\tau_T}^\infty \mathbf{1}_{[-N,N]}(X_t)\,\mathrm{d}t > \frac{\varepsilon_0}{2} \bigg | \tau_T \le \gamma T\bigg) + o(1)\\
&\leq \sup_{y \in [-1,1]}\,\mathbb{P}_y^\mu\bigg(\int_0^\infty \mathbf{1}_{[-N,N]}(X_t)\,\mathrm{d}t > \frac{\varepsilon_0}{2}T \bigg) + o(1),
\end{split}
\end{equation}
where the final inequality relies on the strong Markov property. The convergence of the last expression to $0$ is readily checked via Markov's inequality. \qed

\section{Other corollaries of Theorem \ref{SharpAsymptotics}} \label{corollaries}

We proceed to the proof of the large deviations principle announced in Remark \ref{sharprmk}(b).

\begin{cor}\label{ldp} 
Let $\mu\ge0$. The family $\{A_T\}_{T>0}$ satisfies a large deviations principle with rate $T$ and good rate function given by the Fenchel-Legendre transform of $J(\lambda) := \big(\Lambda(\lambda)-\frac{\mu^2}{2}\big)^+$, where we set $\Lambda(\lambda)=0$ for $\lambda<0$.
\end{cor}

\begin{proof}
Define the (lower semicontinuous) Fenchel-Legendre transform of $J$: 
\begin{equation} \label{J star}
J^*(\beta) = \sup_{\lambda \in \mathbb{R}}\, (\lambda\beta - J(\lambda)) 
= \max\Big(\sup_{\lambda \leq \lambda_c}\lambda \beta,\,\sup_{\lambda\ge\lambda_c} \Big(\lambda \beta - \Lambda(\lambda) +\frac{\mu^2}{2}\Big)\Big).
\end{equation}
Recall Lemma \ref{translemma}. If $\beta\in[\beta_c,1)$, we choose $\lambda\ge\lambda_c$ such that $\Lambda'(\lambda)=\beta$, which makes $\lambda$ the unique critical point of the strictly concave $\lambda \beta -\Lambda(\lambda) + \frac{\mu^2}{2}$ on $[\lambda_c,\infty)$, i.e., its unique global maximum. In particular, 
\begin{equation}
\lambda \beta - \Lambda(\lambda) + \frac{\mu^2}{2} \geq \lambda_c \beta = \sup_{\lambda \leq \lambda_c}\lambda \beta,
\end{equation}
and so $J^*(\beta) = \lambda \beta - \Lambda(\lambda)+\frac{\mu^2}{2}$. If $\beta \in [0,\beta_c)$, then the first-order characterization of convexity yields
\begin{equation}
\Lambda(\lambda) - \frac{\mu^2}{2}\geq \beta_c(\lambda - \lambda_c) \geq \beta (\lambda - \lambda_c),\quad \lambda \geq \lambda_c. 
\end{equation}
Consequently, 
\begin{equation}
\sup_{\lambda \geq \lambda_c} \Big(\lambda \beta - \Lambda(\lambda) + \frac{\mu^2}{2}\Big) \le \lambda_c \beta,
\end{equation}
and $J^*(\beta) = \lambda_c\beta$. If $\beta<0$, we see from the right-hand side of \eqref{J star} that $J^*(\beta)=\infty$. If $\beta>1$, then $\Lambda'(\lambda)<1$ for $\lambda\ge\lambda_c$ implies that $J^*(\beta)=\infty$. Thus,
\begin{equation}
J^*(\beta) = \begin{cases}
    \lambda_c \beta, & \beta \in [0,\beta_c),\\
    \lambda \beta - \Lambda(\lambda)+\frac{\mu^2}{2}, & \beta \in [\beta_c,1),\;\;\,\text{where}\;\Lambda'(\lambda)=\beta, \\
    \infty, & \beta \in (-\infty,0)\cup(1,\infty).
\end{cases}
\end{equation}

\smallskip

For the large deviation upper bound, pick a closed $E\subset\mathbb{R}$. Since $\{A_T\}_{T>0}$ takes values in $[0,1]$, we may assume that $E\subset[0,1]$. If $\inf E \in [0,1)$, Theorem \ref{SharpAsymptotics} yields
\begin{equation}
\limsup_{T \to \infty}\,\frac{1}{T}\log \mathbb{P}(A_T \in E) \le \limsup_{T \to \infty}\frac{1}{T}\log \mathbb{P}(A_T \ge \inf E) = - J^*(\inf E) = -\inf_{\beta \in E} J^*(\beta).
\end{equation}
If $E = \{1\}$, then for any $\varepsilon\in(0,1)$, we infer from Theorem \ref{SharpAsymptotics} that 
\begin{equation*}
\limsup_{T \to\infty}\,\frac{1}{T}\log \mathbb{P}(A_T \in E) \le 
\limsup_{T \to\infty}\,\frac{1}{T}\log \mathbb{P}(A_T \ge 1 - \varepsilon) = -J^*(1-\varepsilon).
\end{equation*}
It remains to send $\varepsilon\to0$ and to use the lower semicontinuity of $J^*$. 

For the large deviation lower bound, let $E\subset\mathbb{R}$ be open. We may assume the existence of $[\underline{e}_n,\overline{e}_n)\subset E\cap[0,1]$ such that $1>\overline{e}_n>\underline{e}_n\to\inf E\cap[0,1]$. By Theorem \ref{SharpAsymptotics},
\begin{equation}
\liminf_{T \to \infty}\,\frac{1}{T}\log \mathbb{P}(A_T \in E) \ge 
\liminf_{T \to \infty}\,\frac{1}{T} \log\big(\mathbb{P}(A_T \ge \underline{e}_n) 
- \mathbb{P}(A_T \ge \overline{e}_n)\big) = -J^*(\underline{e}_n),
\end{equation}
where we have used the fact that $J^*$ is strictly increasing on $[0,1)$. We conclude by sending $n\to\infty$ and by relying on the continuity and monotonicity of $J^*$ on $[0,1]$ (recall \eqref{J star}).  
\end{proof}

We turn to the (conditional) limit theorem for $\{L_0^T-\beta T\}_{T\ge0}$ announced in Remark~\ref{sharprmk}(c). 

\begin{cor} \label{explimit}
Let $\mu\ge0$.
\begin{enumerate}[(a)]
    \item If $\beta \in (\beta_c,1)$, then the conditional law of $L_0^T-\beta T$ given $L_0^T \ge \beta T$ converges to the exponential distribution with parameter $\lambda>0$, where $\Lambda'(\lambda)=\beta$.
    \item If $\beta \in (0,\beta_c]$, then the conditional law of $L_0^T-\beta T$ given $L_0^T \ge \beta T$ converges to the exponential distribution with parameter $\lambda_c$.
\end{enumerate}
\end{cor}

\begin{proof}
    For any fixed $y\ge 0$, we write $\mathbb{P}(L_0^T -\beta T\geq y |L_0^T \geq \beta T)$ as a ratio, and apply Theorem \ref{SharpAsymptotics} to the numerator and the denominator.
\end{proof}

We close with a limit theorem for $\tau_T = \inf \{t \ge 0\!: L_0^t \ge \beta T \}$ announced in Remark~\ref{sharprmk}(d).

\begin{cor}
Let $\mu\ge0$.
\begin{enumerate}[(a)]
    \item If $\beta \in (\beta_c,1)$, then the conditional law of $T-\tau_T$ given $L_0^T \geq \beta T$ converges to the exponential distribution with parameter $\Lambda_\mu(\lambda)>0$, where $\Lambda'(\lambda)=\beta$. 
    \item  If $\beta=\beta_c\in(0,1)$, then the conditional law of $(T-\tau_T)/\sqrt{T}$ given $L_0^T \geq \beta T$ converges to the one-sided normal distribution with mean $0$ and variance $\Lambda''(\lambda_c)/\beta_c^2>0$.
    \item If $\beta \in (0,\beta_c)$, then the conditional law of $(\tau_T-(\beta/\beta_c)T)/\sqrt{T}$ given $L_0^T \geq \beta T$ converges to the normal distribution with mean $0$ and variance $(\beta/\beta_c)\Lambda''(\lambda_c)/\beta_c^2>0$. 
\end{enumerate}
\end{cor}

   \begin{proof}
    	\textit{(a).} Note that for any $y \ge 0$,
    	\begin{equation}
    	\mathbb{P}(T - \tau_T\ge y | L_0^T \ge \beta T) = \frac{\mathbb{P}(L_0^{T-y} \geq \beta T)}{\mathbb{P}(L_0^T \geq \beta T)}.
    	\end{equation}
    	It remains to apply Theorem \ref{SharpAsymptotics}(a).
    	
    	\medskip
    
        \noindent\textit{(b), (c).} Fix $y\in\mathbb{R}$ and let $r_T:=(\beta/\beta_c) T + y\sqrt{T}$. 
        By $ \{\tau_T\le r_T\} = \{L_0^{r_T}\ge \beta_c r_T -\beta_c y \sqrt{T} \}$, Theorem \ref{SharpAsymptotics}(b) implies that 
        \begin{equation}
        \begin{split}
        \mathbb{P}(\tau_T\le (\beta/\beta_c) T + y\sqrt{T}) &\sim C(x)\, 
        e^{-r_T\lambda_c \beta_c}\, e^{\lambda_c \beta_c y \sqrt{T}}
        \,\Phi\bigg(\frac{\beta_c y\sqrt{T}}{\sqrt{r_T\Lambda''(\lambda_c)}} \bigg) \\
        & \sim C(x)\, e^{-T\lambda_c \beta}\,
        \Phi\bigg(\frac{\beta_c y}{\sqrt{(\beta/\beta_c)\Lambda''(\lambda_c)}} \bigg).
        \end{split}
        \end{equation}
        Together with $\mathbb{P}(L_0^T \ge \beta T) \sim C(x)\, e^{-T\lambda_c\beta}$ (Theorem \ref{SharpAsymptotics}(c)), this proves part (c). Lastly, for $\beta=\beta_c\in(0,1)$, 
       	Theorem \ref{SharpAsymptotics}(b) yields $\mathbb{P}(L_0^T \ge \beta T) \sim \frac{1}{2}C(x)\, e^{-T\lambda_c\beta_c}$, which proves part~(b) upon observing that only $y\le 0$ are relevant. 
   \end{proof}

\appendix 

\section{Generalizations}

In this appendix, we show several results necessary to extend our proof of Theorem \ref{MainTheorem} to other Markov processes and additive functionals.~For the Markov process, we let $\mcal{X}$ be the state space, use $\nu$ for a measure on $\mcal{X}$, and write $L^2 = L^2(\mcal{X},\nu)$.~We make the following regularity assumptions throughout. 

\begin{asmp}\label{appass}	
	\begin{enumerate}[(a)]
		\item $\mathcal{X}$ is a locally compact, separable metric space, and $\nu$ is a Radon measure on $\mathcal{X}$ with full support such that $\nu(\mcal{X})=\infty$.
		\item $(\Omega,(\mcal{F}_t)_{t \geq 0},(X_t)_{t \geq 0},\{\mathbb{P}_x\}_{x \in \mcal{X}})$ is a Hunt process on $\mcal{X}$ (with the Borel $\sigma$-algebra) in the sense of \cite[Appendix A.2]{Fukushima1994-eo}. 
		\item The semigroup $P_t f(x) = \mathbb{E}_x[f(X_t)]$ on $L^2$ is well-defined, strongly continuous, self-adjoint, and irreducible in the sense of \cite[Section 1.3, and p.~55]{Fukushima1994-eo}.
		\item For all $(t,x)\in(0,\infty)\times\mcal{X}$, the law of $X_t$ under $\mathbb{P}_x$ is absolutely continuous with respect to $\nu$, and its density $p_t(x,\cdot)$ belongs to $L^2$. Furthermore, 
		\begin{equation} \label{trans dens}
		P_t f(x) :=\int_{\mcal{X}} p_t(x,y)\, f(y)\,\nu(\mathrm{d}y),\quad (t,f)\in(0,\infty)\times L^2
		\end{equation}
		are continuous in $x$, and $x\mapsto\|p_t(x,\cdot)\|_{L^2}$, $t\in(0,\infty)$ are locally bounded.
		\item Let $\mcal{L}$ be the generator associated with the semigroup $(P_t)_{t \geq 0}$ on $L^2$. Then, for every function $f \in L^\infty(\mcal{X})$ with compact support, the multiplication operator $M_f$ on $L^2$ defined by $M_f g = fg$ is compact from $(\mcal{D}(\mcal{L}),\norm{\cdot}_{\mcal{D}(\mcal{L})})$ to $L^2$.
	\end{enumerate}
\end{asmp}
The definition of irreducibility in \cite{Fukushima1994-eo} may appear opaque. However, given the other conditions in Assumption \ref{appass}, it will become clear from the proofs that this condition is equivalent to the fact that $p_t(x,y)>0$ for $\nu\times\nu$-almost every $(x,y)$. 
Under Assumption \ref{appass}, we prove a Feynman-Kac formula.

\begin{prop} \label{appFeynman}
	Let $0\leq V \in L^\infty(\mathcal{X})$ have compact support, and let $\mcal{L}$ be the generator of a process $X$ satisfying Assumption \ref{appass}. Then, for every $z \in \mathbb{C}$, the operator $\mcal{L}_z := \mcal{L}+zV$ with domain $\mcal{D}(\mcal{L}_z) := \mcal{D}(\mcal{L})$ is closed and generates an analytic semigroup. Furthermore, for every $(t,f)\in(0,\infty)\times L^2$, $e^{t \mcal{L}_z}f$ admits a continuous representative. Finally, we have the Feynman-Kac formula
	\begin{equation}\label{appFK}
		e^{t \mcal{L}_z}f(x) = \mathbb{E}_x\bigg[f(X_t)\,
		\exp\bigg(z \int_0^t V(X_s)\,\mathrm{d}s\bigg)\bigg],\quad x\in\mcal{X}.
	\end{equation}
\end{prop}

\begin{proof}
	As the generator of a strongly continuous contraction semigroup, $\mcal{L}$ is densely defined and closed by the Hille-Yosida Theorem (see, e.g., \cite[Theorem II.3.5]{Engel1999-tt}), which implies that $\mcal{L}_z$ is closed. 
	As mentioned in \cite[Lemma 1.3.1]{Fukushima1994-eo}, $\mcal{L}$
	must be a non-positive definite self-adjoint operator. Hence, the spectral theorem in the form of \cite[Theorem I.4.9]{Engel1999-tt} implies that $\sigma(\mcal{L})\subset(-\infty,0]$. One can follow Step 1 in the proof of Proposition \ref{MainGeneration} to conclude that $\mcal{L} $ and $\mcal{L}_z$ generate analytic semigroups. One can also repeat Step 3 in the proof of Proposition \ref{MainGeneration} to obtain the Feynman-Kac formula, but using Assumption \ref{appass}(d) to upgrade from $\nu$-almost every $x$ to every $x\in\mcal{X}$.
\end{proof}

Going forward, we make the following assumption on the potential $V$.

\begin{asmp}\label{potentialassumption}
	The potential $0\!\le\! V \!\in\! L^\infty(\mathcal{X})$ has compact support and is not identically~$0$.
\end{asmp}

\begin{thm}\label{Generalspectrum}
		Suppose that Assumptions \ref{appass} and \ref{potentialassumption} hold.~For $z\in\mathbb{R}$, let $\Lambda(z):=\sup\,\sigma(\mcal{L}_z)$. 
	\begin{enumerate}[(a)]
		\item If the semigroup $(e^{t\mcal{L}})_{t\ge0}$ is recurrent in the sense of \cite[display~(1.6.2)]{Fukushima1994-eo}, then for all $\lambda>0$, $\Lambda(\lambda)>0$ is an isolated, algebraically simple eigenvalue corresponding to a continuous, strictly positive eigenfunction $\phi_\lambda$ of norm $1$. Moreover, $\lim_{\lambda\downarrow 0}\Lambda(\lambda)=\Lambda(0)=0$. In this case, we define $\lambda_0 =\lambda_c = 0$.
		\item If $\sigma(\mcal{L})\subset (-\infty,-\varepsilon]$ for some $\varepsilon>0$, then there exists a $\lambda_0 >0$ with $\Lambda(\lambda_0)=-\varepsilon/2$ such that for all $\lambda\ge\lambda_0$, $\Lambda(\lambda)$ is an isolated, algebraically simple eigenvalue corresponding to a continuous, strictly positive eigenfunction $\phi_{\lambda}$ of norm $1$. Moreover, there is a unique $\lambda_c >\lambda_0$ with $\Lambda(\lambda_c)=0$.
	\end{enumerate}
\end{thm}

\begin{proof}
	Let $\mcal{E}$ denote the Dirichlet form associated with $\mcal{L}$. From the spectral theorem in the form of \cite[Theorem I.4.9]{Engel1999-tt} we infer that for all $\lambda \geq 0$,
	\begin{equation}\label{dirichlet}
	\begin{split}
	\Lambda(\lambda) = \sup_{f \in \mcal{D}(\mcal{L}),\,\norm{f}_{L^2}=1} 
	\langle \mcal{L}_\lambda f,f \rangle_{L^2} 
	&= \sup_{f \in \mcal{D}(\mcal{L}),\,\norm{f}_{L^2}=1} 
	\bigg(\langle\mcal{L} f, f \rangle_{L^2}+\lambda\int_{\mcal{X}} V \, |f|^2\,\mathrm{d}\nu\bigg) \\
	&= \sup_{f \in \mcal{D}(\mcal{E}),\,\norm{f}_{L^2}=1} 
	\bigg(-\mcal{E}(f,f) + \lambda\int_{\mcal{X}} V \,|f|^2\,\mathrm{d}\nu\bigg).
	\end{split}
	\end{equation}
	The last equality follows from the fact that $\mcal{D}(\mcal{L})$ (which is also the range of the resolvent) is dense in $\mcal{D}(\mcal{E})$ when the latter space is equipped with the inner product $\mcal{E}_1(f,g):=\mcal{E}(f,g) + \langle f,g\rangle_{L^2}$ (cf.~\cite[proof of Theorem 1.3.1]{Fukushima1994-eo}). The identity \eqref{dirichlet} shows, in particular, that $\Lambda$ is Lipschitz on $[0,\infty)$. Indeed, for any $f\in\mcal{D}(\mcal{E})$ with $\|f\|_{L^2}=1$ and 
	$\lambda,\lambda'\ge 0$, 
	\begin{equation}
	-\mcal{E}(f,f) + \lambda \int_{\mcal{X}} V \,|f|^2\,\mathrm{d}\nu 
	=- \mcal{E}(f,f) + \lambda'\int_{\mcal{X}} V \,|f|^2\,\mathrm{d}\nu + (\lambda-\lambda')\int_{\mcal{X}} V \,|f|^2\,\mathrm{d}\nu. 
	\end{equation}	
	Taking the supremum over $f$, we find that $\Lambda(\lambda)\le\Lambda(\lambda')+|\lambda-\lambda'|\,\|V\|_\infty$. Reversing the roles of $\lambda$ and $\lambda'$ yields the Lipschitz property of $\Lambda$.
	
	\medskip
	
	Next, consider the case in part (a) and let $\lambda>0$. Since the semigroup is recurrent, \cite[Theorem 1.6.3]{Fukushima1994-eo} implies the existence of $f_n \in \mcal{D}(\mcal{E})$ such that $f_n\to1$ $\nu$-almost everywhere and $\mcal{E}(f_n,f_n)\to 0$. By Fatou's lemma, 
	\begin{equation}
	\liminf_{n \to \infty} \,\bigg(-\mcal{E}(f_n,f_n) 
	+ \lambda \int_{\mcal{X}} V \,|f_n|^2\,\mathrm{d}\nu\bigg) 
	\ge \lambda \int_{\mcal{X}} V \,\mathrm{d}\nu > 0.
	\end{equation}
	Hence, the term inside the limit inferior is positive for some $n$ which, upon dividing it by $\|f_n\|_{L^2}^2$, implies $\Lambda(\lambda)>0$ . Then, by Assumption \ref{appass}(e) and Proposition \ref{WeylTheorem}, $\Lambda(\lambda)$ must be an isolated eigenvalue of finite algebraic multiplicity. Since $\Lambda(\lambda)\le\lambda\|V\|_\infty$ by \eqref{dirichlet}, it holds that $\lim_{\lambda\downarrow0} \Lambda(\lambda)=0$. Moreover, $\Lambda(0)=0$. Indeed, it is clear that $\Lambda(0)\le 0$. To see $\Lambda(0)\ge0$, we use Fatou's lemma: 
	\begin{equation}
	\liminf_{n \to \infty} \norm{f_n}_{L^2}^2 
	\ge \liminf_{n \to \infty} \frac{1}{\norm{V}_\infty}\int_{\mcal{X}} V \,|f_n|^2\,\mathrm{d}\nu 
	\ge \frac{1}{\norm{V}_\infty} \int_{\mcal{X}} V \,\mathrm{d}\nu > 0.
	\end{equation}
	Therefore, $\mcal{E}\big(\frac{f_n}{\norm{f_n}_{L^2}},\frac{f_n}{\norm{f_n}_{L^2}}\big) 
	\to 0$, and so $\Lambda(0)\geq 0$ thanks to \eqref{dirichlet}.
	
	\medskip
	
	Now, consider the case in part (b). Since $\mcal{D}(\mcal{L})$ is dense in $L^2$, there exists an $f \in\mcal{D}(\mcal{L})$ with $\norm{f}_{L^2}=1$ such that $\int_{\mcal{X}} V \, |f|^2\,\mathrm{d}\nu>0$. Together with \eqref{dirichlet}, this shows that $\Lambda(\cdot)$ 
	is non-decreasing on $[0,\infty)$ and that $\Lambda(\lambda)\to\infty$ as $\lambda\to\infty$. Set $\lambda_0 = \inf\{\lambda \ge0:\,\Lambda(\lambda)\ge -\varepsilon/2\}$. Since $\Lambda(0)\le -\varepsilon$, $\Lambda(\cdot)$ is non-decreasing and continuous, and $\lim_{\lambda\to\infty} \Lambda(\lambda)=\infty$, it follows that $\lambda_0\in(0,\infty)$ and $\Lambda(\lambda_0)=-\varepsilon/2$. For $\lambda\geq \lambda_0$, Proposition \ref{WeylTheorem} implies that $\Lambda(\lambda)$ is an isolated eigenvalue of finite algebraic multiplicity. Let $f_\lambda$ be a corresponding eigenfunction of norm $1$.
	Note that $Vf_\lambda$ is not identically $0$ because $\lambda \ge \lambda_0$. We claim that $\Lambda(\cdot)$ is strictly increasing on $[\lambda_0,\infty)$. Indeed, for any $\lambda'>\lambda\ge\lambda_0$,
	\begin{equation}
	\Lambda(\lambda') - \Lambda(\lambda)\ge (\lambda' - \lambda) \int_{\mcal{X}} V \,|f_\lambda|^2\,\mathrm{d}\nu > 0.
	\end{equation}
	This, along with the continuity of $\Lambda(\cdot)$, yields a unique $\lambda_c>\lambda_0$ such that $\Lambda(\lambda_c)=0$.
	
	\medskip
	
	Finally, we argue that in both cases, the eigenfunction is continuous, strictly positive, and unique up to multiplication by constants. To this end, we check first that the semigroup is positivity improving in the following sense: For every $0\le f \in L^2$ which is not identically~$0$, $P_t f$ is positive $\nu$-almost everywhere. We aim to use \cite[Theorem 3]{MR83}. (While the authors never define a holomorphic semigroup, they need that for every $f,g \in L^2$, the map $t \mapsto \langle P_t f,g \rangle_{L^2}$ extends to a holomorphic map on a sector $\{z\in\mathbb{C}:\,|\arg(z)|<\delta\}\setminus\{0\}$. This follows from the fact that $\mcal{L}$ generates an analytic semigroup.) 
	
	\medskip
	
	We start by verifying that for every $0\le f,g \in L^2$ which are not identically $0$, there exists a $t>0$ such that $\langle P_t f, g \rangle_{L^2} > 0$. Suppose for contradiction that there are such $f$, $g$ satisfying $\langle P_t f ,g \rangle_{L^2}=0$ for all $t>0$. For $\alpha>0$, it holds that
	\begin{equation} \label{semi -> res}
	(\alpha-\mcal{L})^{-1} f = \int_0^\infty e^{-\alpha t}\, P_tf \,\mathrm{d}t,
	\end{equation}
	and so $\langle (\alpha-\mcal{L})^{-1} f, g \rangle_{L^2} = 0$. Fix an $\alpha>0$, and let $E:=\{(\alpha-\mcal{L})^{-1} f > 0\}$. Since $f$ is not identically $0$, we have $\nu(E)>0$. We claim that $E$ is invariant in the sense defined above \cite[Lemma 1.6.1]{Fukushima1994-eo}. To see the claim, fix any $h\in L^2$ such that $h=0$ $\nu$-almost everywhere on $E^c$. Define $h_n = h \,\mathbf{1}_{\{(\alpha-\mcal{L})^{-1} f \ge 1/n\}}\,\mathbf{1}_{\{|h|\leq n\}}$, and note that $|h_n|\leq n^2(\alpha-\mcal{L})^{-1} f$ and $h_n\to h$ in $L^2$. Therefore, by using \eqref{semi -> res}, we find that on $E^c$,
	\begin{equation}
	|P_t h_n| \le P_t|h_n| \le n^2 P_t(\alpha-\mcal{L})^{-1}f \le n^2 e^{\alpha t}(\alpha-\mcal{L})^{-1}f  = 0,\quad t>0.
	\end{equation}
	Since $P_t$ is a bounded operator on $L^2$, and so $P_t h_n \to P_t h$ in $L^2$, it follows that $P_t h = 0$ on~$E^c$. Thus, $E$ is invariant. However, $\{g>0\}\subset E^c$, hence $\nu(E^c)>0$. Consequently, $E$ is an invariant set for which $\nu(E)>0$ and $\nu(E^c)>0$, contradicting irreducibility. 
	
	\medskip
	
	In view of \cite[Proposition 1, Theorem 3]{MR83}, for all $0\le f,g \in L^2$ which are not identically~$0$, we have $\langle P_t f, g \rangle_{L^2}>0$ for all $t>0$. Therefore, for all $0\le f\in L^2$ which are not identically~$0$, it holds that $P_t f >0$ $\nu$-almost everywhere for all $t>0$. This implies that for $\lambda>0$, the semigroup $(e^{t\mcal{L}_\lambda})_{t\ge0}$ is also positivity improving, since by the Feynman-Kac formula \eqref{appFK},
	\begin{equation}
	e^{t \mcal{L}_\lambda}f(x)
	=\mathbb{E}_x\big[f(X_t)\, e^{\lambda\int_0^t V(X_s)\,\mathrm{d}s}\big]
	\ge P_tf(x),\quad t>0. 
	\end{equation}
	Consequently, it follows from \cite[Theorem XIII.44]{ReedSimonIV} that $\Lambda(\lambda)$ is a geometrically simple eigenvalue corresponding to a $\nu$-almost everywhere positive eigenfunction $\phi_\lambda$. Moreover, the algebraic simplicity of $\Lambda(\lambda)$ can be deduced by repeating the last paragraph in the proof of Proposition \ref{MainGeneration}. The identity $\phi_\lambda=e^{-\Lambda(\lambda)t}\, e^{t\mcal{L}_\lambda}\phi_\lambda$ and Lemma \ref{appFeynman} show further that $\phi_\lambda$ can be chosen as a continuous function. Lastly, the identity \eqref{trans dens} and $\phi_\lambda>0$ $\nu$-almost everywhere imply that for all $t>0$, the continuous $x\mapsto P_t\phi_\lambda(x)$ is positive for all $x\in\mcal{X}$. Thus, the chain of inequalities
	\begin{equation} \label{FK eigen}
	\phi_\lambda(x)= e^{-\Lambda(\lambda)t}
	\,\mathbb{E}_x\big[\phi_\lambda(X_t)\,e^{\lambda\int_0^tV(X_s)\,\mathrm{d}s}\big] 
	\ge e^{-\Lambda(\lambda)t}\, P_t\phi_\lambda(x) > 0,\quad (t,x)\in(0,\infty)\times\mcal{X}
	\end{equation}
	concludes the proof. 
\end{proof}

Theorem \ref{Generalspectrum} gives the analogues of the statements (a), (b), (d), and~(e) in Proposition~\ref{MainPerturbation} for $\lambda>\lambda_0$. In particular, for $z\in\mathbb{C}$ close to $\lambda$, the operator $\mcal{L}_z$ has an isolated, algebraically simple eigenvalue $\Lambda(z)$, and the map $z \mapsto\Lambda(z)$ is holomorphic. For statement (c), we fix a compact $K\subset\mcal{X}$ and use $\|e^{t\mcal{L}_z}\|_{B(L^2,C(K))} \le e^{t|z|\norm{V}_\infty}\,\sup_{x \in K}\norm{p_t(x,\cdot)}_{L^2}$ (resulting from the Feynman-Kac formula \eqref{appFK}) to deduce that $z\mapsto e^{t\mcal{L}_z}$ maps into $B(L^2,C(K))$, the space of bounded linear operators $L^2\to C(K)$. Moreover, with
\begin{equation}
	(e^{t\mcal{L}_z})'f(x):=\mathbb{E}_x\bigg[f(X_t)\int_0^t V(X_s)\,\mathrm{d}s\,
	\exp\bigg(z \int_0^t V(X_s)\,\mathrm{d}s\bigg)\bigg],\quad (f,x)\in L^2\times\mcal{X},
\end{equation}
the formula \eqref{appFK} and the inequality $\big|\frac{e^{h\ell}-1-h\ell}{h}\big|\le |h|\,\frac{\ell^2}{2}\, e^{|h|\ell}$, $|h|\le 1$, $\ell\ge0$ imply that
\begin{equation}
\sup_{x\in K}\,\bigg|\frac{e^{t\mcal{L}_{z+h}}f(x)-e^{t\mcal{L}_z}f(x)}{h} -(e^{t\mcal{L}_z})'f(x)\bigg|
\le |h|\,\frac{t^2\|V\|_\infty^2}{2}\,e^{t\|V\|_\infty}\,\sup_{x\in K} e^{t\mcal{L}_{\mathrm{Re}(z)}}|f|(x),
\end{equation}
and so $z\mapsto e^{t\mcal{L}_z}$ is a holomorphic $B(L^2,C(K))$-valued function. Let $P^0_z$ denote the Riesz projection associated with $\Lambda(z)$. Then, $P^0_z = e^{-\Lambda(z)t}\, e^{t\mcal{L}_z}P^0_z$ for all $t>0$. Since $P^0_z$ is holomorphic as an operator on $L^2$, $e^{t\mcal{L}_z}$ is holomorphic as an operator from $L^2$ to $C(K)$, and $z\mapsto\Lambda(z)$ is holomorphic, we have that for any $f \in L^2$, $P^0_zf$ can be chosen as a continuous function that converges locally uniformly to $P^0_\lambda f$ as $z \rightarrow \lambda$. Finally, $P^0_\lambda$ is the orthogonal projection onto the span of $\phi_\lambda$, and so $P^0_\lambda f = \langle f, \phi_\lambda \rangle_{L^2}\, \phi_\lambda$. Thus, the analogue of the statement (c) in Proposition \ref{MainPerturbation} holds.

\medskip

We also have analogues of Lemma \ref{translemma} and Proposition \ref{SDEfacts}.

\begin{lemma} \label{generallambdaprop}
	Suppose that Assumptions \ref{appass} and \ref{potentialassumption} hold. Then, $\Lambda\in C^2((\lambda_0,\infty))$ and $\Lambda''>0$. Moreover, in the case of Theorem \ref{Generalspectrum}(a), $\Lambda'$ is a bijection from $(0,\infty)$ to $(0,\|V\|_\infty)$; and in the case of Theorem \ref{Generalspectrum}(b), it is a bijection from $(\lambda_c,\infty)$ to $(\Lambda'(\lambda_c),\|V\|_\infty)$. 
\end{lemma}

\begin{proof}
	To find the derivatives of $\Lambda$, we use the formulas in \cite[Chapter II, display (2.33)]{Kato1995-kx} which remain valid in our case due to \cite[Remark VII.2.10]{Kato1995-kx}. Fix $\lambda>\lambda_0$, let $P^0_\lambda$ denote the Riesz projection onto the eigenspace of $\Lambda(\lambda)$, set $Q^0_\lambda = I-P^0_\lambda$, and write $\widetilde{\mcal{L}}_\lambda$ for the operator $\mcal{L}_\lambda$ restricted to $\mcal{D}(\mcal{L})\cap\mathrm{Range}(Q^0_\lambda)$. The spectrum of $\widetilde{\mcal{L}}_\lambda-\Lambda(\lambda)$ is upper bounded by $-\delta$ for some $\delta>0$. Therefore, $(\Lambda(\lambda)-\widetilde{\mcal{L}}_\lambda)^{-1}$ is a bounded, positive definite, self-adjoint operator. Define the reduced resolvent $S_\lambda=(\Lambda(\lambda)-\widetilde{\mcal{L}}_\lambda)^{-1}Q^0_\lambda = Q^0_\lambda(\Lambda(\lambda)-\widetilde{\mcal{L}}_\lambda)^{-1}Q^0_\lambda$. The formulas in \cite[Chapter II, display (2.33)]{Kato1995-kx} yield
	\begin{equation}
	\Lambda'(\lambda)= \mathrm{trace}(VP^0_\lambda),\quad 
	\Lambda''(\lambda)=  2\,\mathrm{trace}(VS_\lambda VP^0_\lambda).
	\end{equation}
	With $P^0_\lambda$ being the orthogonal projection onto the span of $\phi_\lambda$, these simplify to
	\begin{equation} \label{Lambda der form}
	\Lambda'(\lambda) = \langle V \phi_\lambda, \phi_\lambda \rangle_{L^2}, \quad \Lambda''(\lambda)= 2\langle VS_\lambda V \phi_\lambda,\phi_\lambda\rangle_{L^2} 
	=2\langle (\Lambda(\lambda)-\widetilde{\mcal{L}}_\lambda)^{-1}Q^0_\lambda V\phi_\lambda,Q^0_\lambda V\phi_\lambda \rangle_{L^2},
	\end{equation}
	where
	\begin{equation}
	Q^0_\lambda V \phi_\lambda 
	= V\phi_\lambda - \langle V \phi_\lambda, \phi_\lambda \rangle_{L^2} \phi_\lambda 
	= (V-\Lambda'(\lambda))\,\phi_\lambda.
	\end{equation}
	Since $\phi_\lambda$ is strictly positive, $V$ has compact support, $\nu(\mcal{X})=\infty$, and $V$ is not identically~$0$, we have $Q^0_\lambda V \phi_\lambda\neq0$. Thus, the final expression in \eqref{Lambda der form} implies that $\Lambda''(\lambda)>0$.
	
	\medskip
	
 	Turning to $\Lambda'$, we read off from $\Lambda'(\lambda)=\langle V\phi_\lambda, \phi_\lambda \rangle_{L^2}$ that $\Lambda'(\lambda) \in [0,\|V\|_\infty]$. Moreover, $\Lambda'$ is strictly increasing on $(\lambda_0,\infty)$ since $\Lambda''>0$. To show that $\lim_{\lambda \to \infty} \Lambda'(\lambda)=\|V\|_\infty$, it suffices to verify that $\lim_{\lambda \to \infty} \Lambda(\lambda)/\lambda=\|V\|_\infty$. The identity \eqref{dirichlet} implies directly that $\Lambda(\lambda)\le \lambda\|V\|_\infty$. For a lower bound, let $v\in(0,\|V\|_\infty)$ and $\mcal{X}_v := \{V>\|V\|_\infty-v\}$. Since $V$ has compact support, $\nu(\mcal{X}_v)\in(0,\infty)$. Set $g_v = \mathbf{1}_{\mcal{X}_v}/\sqrt{\nu(\mcal{X}_v)}$, and note that $\langle V g_v, g_v \rangle_{L^2} \ge \|V\|_\infty-v$. Since $\mcal{D}(\mcal{L})$ is dense, there is an $f_v \in \mcal{D}(\mcal{L})$ such that $\norm{f_v}_{L^2}=1$ and $\langle V f_v, f_v \rangle_{L^2} \ge \|V\|_\infty-2v$. Then, $\Lambda(\lambda) \geq \langle \mcal{L}f_v,f_v \rangle_{L^2} + \lambda(\|V\|_\infty-2v)$ by \eqref{dirichlet}. The arbitrariness of $v\in(0,\|V\|_\infty)$ reveals that $\liminf_{\lambda\rightarrow \infty} \Lambda(\lambda)/\lambda\ge\|V\|_\infty$. 
	
	\medskip
	
	In the case of Theorem \ref{Generalspectrum}(b), $\lim_{\lambda\downarrow \lambda_c} \Lambda'(\lambda)=\Lambda'(\lambda_c)$ follows from the holomorphicity of $\Lambda(\cdot)$. It remains to prove that in the case of Theorem \ref{Generalspectrum}(a), $\lim_{\lambda \downarrow 0}\Lambda'(\lambda)=0$. It suffices to establish the existence of a sequence $\lambda_n\downarrow 0$ such that $\Lambda'(\lambda_n)\downarrow 0$. Start with any $\lambda_n \downarrow 0$. It is clear from $\lim_{\lambda\downarrow0} \Lambda(\lambda)=0$ and the boundedness of $V$ that $\mcal{L}\phi_{\lambda_n} = \Lambda(\lambda_n)\phi_{\lambda_n} - \lambda_n V \phi_{\lambda_n}$ converges to $0$ in $L^2$. In addition, since $\|\phi_{\lambda_n}\|_{L^2}=1$, the Banach–Alaoglu theorem reveals that $(\phi_{\lambda_n})_{n\in\mathbb{N}}$ converges weakly along a subsequence to some $\phi\in L^2$. Due to Mazur's lemma and the closedness of $\mcal{L}$, we find that $\phi\in\mcal{D}(\mcal{L})$ and $\mcal{L}\phi=0$. We claim that $\phi=0$. Indeed, $P_t \phi = \phi$ for any $t>0$, and so 
	\begin{equation} \label{const eigen}
	0 = 2\,(\norm{\phi}_{L^2}^2 - \langle P_t \phi,\phi \rangle_{L^2}) = \int_{\mcal{X}} \int_{\mcal{X}} |\phi(x)-\phi(y)|^2\, p_t(x,y)\, \nu(\mathrm{d}x)\,\nu(\mathrm{d}y).
	\end{equation}
	On the other hand, $P_{\frac t 2}(p_{\frac t 2}(x,\cdot))(y) = p_t(x,y)$ for all $x\in\mcal{X}$ and $\nu$-almost every $y\in\mcal{X}$ by the Chapman-Kolmogorov equation. Since $P_{\frac t 2}$ is positivity improving, $p_t(x,y)>0$ for all $x\in\mcal{X}$ and $\nu$-almost every $y\in\mcal{X}$. Together with \eqref{const eigen}, this yields the constancy of $\phi$. However, due to $\phi\in L^2$ and $\nu(\mcal{X})=\infty$, the only possibility is $\phi=0$. Finally, since $\|\phi_{\lambda_n}\|_{L^2}=1$ and $\|\mcal{L}\phi_{\lambda_n}\|_{L^2}\to 0$, the sequence $(\phi_{\lambda_n})_{n\in\mathbb{N}}$ is bounded in the graph norm. By Assumption~\ref{appass}(e), after passing to a further subsequence, $V\phi_{\lambda_n}$ converges (strongly) in $L^2$. Together with the weak convergence $\phi_{\lambda_n}\to0$, we find that $\Lambda'(\lambda_n)=\langle V\phi_{\lambda_n}, \phi_{\lambda_n} \rangle_{L^2} \to 0$.
\end{proof}

\begin{prop} 
	Suppose that Assumptions \ref{appass} and \ref{potentialassumption} hold. Let $\lambda>\lambda_0$ and set $\pi_\lambda(\mathrm{d}x) = \phi_\lambda(x)^2\, \nu(\mathrm{d}x)$. Then:
	\begin{enumerate}[(a)]
	 \item The semigroup $(Q^\lambda_t)_{t \ge 0}$ defined by $Q_t^\lambda f = \phi_\lambda^{-1}e^{-\Lambda(\lambda)t}\, e^{t\mcal{L}_\lambda}(\phi_\lambda f)$ is strongly continuous on $L^2(\pi_\lambda)$. For any $(t,f) \in (0,\infty)\times L^2(\pi_\lambda)$, one can choose $Q_t^\lambda f$ to be continuous.
	\item For any $x\in\mcal{X}$ and $t\ge0$, the probability measure that has the density
	\begin{equation}
		M_t := \frac{\phi_\lambda(X_t)}{\phi_\lambda(x)}\, e^{-t\Lambda(\lambda)}\, \,\exp\bigg(\lambda\int_0^t V(X_s)\,\mathrm{d}s\bigg)
	\end{equation}
	with respect to the law of $X_{[0,t]}$ defines a Markov process $Y^\lambda$ realizing the semigroup~$Q^\lambda$.
	\item Define $G_t = \frac{1}{\sqrt{t}}\,
	\int_0^t V(Y^\lambda_s)-\Lambda'(\lambda)\,\mathrm{d}s$. Then, for all bounded measurable $g:\,\mcal{X}\to\mathbb{R}$, and all bounded continuous $F:\,\mathbb{R}\to\mathbb{C}$, we have
	\begin{equation}\label{convdist}
		\lim_{t\to\infty}\mathbb{E}^x[g(Y^\lambda_t)\, F(G_{t})] = 
		\int_{\mcal{X}} g(y)\,\phi_\lambda(y)^2\,\nu(\mathrm{d}y) \,
		\int_{-\infty}^\infty F(v)\,
		\frac{e^{-\frac{v^2}{2\Lambda''(\lambda)}}}{\sqrt{2\pi \Lambda''(\lambda)}} \,\mathrm{d}v.
	\end{equation}
	\item With $\rho_T(\cdot):=\frac{1}{T}\,\int_0^T \mathbf{1}_\cdot(Y^\lambda_t)\,\mathrm{d}t$, $\rho_T$ converges weakly to $\pi_\lambda$ as $T\to\infty$, almost surely. 
	\end{enumerate}
\end{prop}

\begin{proof}
	\textit{(a).} Define the unitary transformation $R_\lambda:\, L^2 \to L^2(\pi_\lambda)$, $f\mapsto\phi_\lambda^{-1}f$. We have the representation $Q_t^\lambda = R_\lambda\, e^{t(\mcal{L}_\lambda - \Lambda(\lambda))} R_\lambda^{-1}$. Since $\mcal{L}_\lambda$ generates a strongly continuous semigroup on $L^2$, $(e^{t(\mcal{L}_\lambda-\Lambda(\lambda))})_{t\ge0}=(e^{-t\Lambda(\lambda)}e^{t\mcal{L}_\lambda})_{t\ge0}$ is a strongly continuous semigroup on~$L^2$. Hence, $(Q_t^\lambda)_{t\ge0}$ is a strongly continuous semigroup on $L^2(\pi_\lambda)$ (cf.~the first paragraph in the proof of Theorems \ref{DriftGeneration} and \ref{DriftSpectrum} for $\mu>0$). Moreover, $Q_t^\lambda f$ can be chosen as a continuous function because $e^{t\mcal{L}_\lambda}(\phi_\lambda f)$ and $\phi_\lambda^{-1}$ are continuous.  
	
	\medskip
	
	\noindent\textit{(b).} Using the Markov property of $X$ and the identity in \eqref{FK eigen}, we find that $(M_t)_{t\ge0}$ is a martingale. Defining a probability measure via the density $M_t$ with respect to the law of~$X_{[0,t]}$, one can verify using Bayes' rule for conditional expectations, the Markov property of $X$, and the Feynman-Kac formula \eqref{appFK} that any conditional expectation $\mathbb{E}[f(X_{s_2})|\mcal{F}_{s_1}]$ is given by $Q_{s_2-s_1}^\lambda f(X_{s_1})$.
 	
 	\medskip
	
	\noindent\textit{(c).} The proof of Proposition \ref{SDEfacts}(c) can be repeated mutatis mutandis. 
	
	\medskip
	
	\noindent\textit{(d).} Consider any bounded measurable $f \!:\mcal{X}\to\mathbb{R}$, and define $g = f - \int_{\mcal{X}} f \,\mathrm{d}\pi_\lambda$ and $A_T = \frac{1}{T}\int_0^T g(Y^\lambda_t)\,\mathrm{d}t$. Assume first that the law of $Y^\lambda_0$ is $\pi_\lambda$ and note that for any $h\in L^2(\pi_\lambda)$, 
	\begin{equation}
	\int_{\mcal{X}} Q_t^\lambda h\,\mathrm{d}\pi_\lambda 
	= e^{-\Lambda(\lambda)t}\,\int_{\mcal{X}} \phi_\lambda\, e^{t\mcal{L}_\lambda}(\phi_\lambda h)\,\mathrm{d}\nu 
	= e^{-\Lambda(\lambda)t}\,\int_{\mcal{X}} (e^{t\mcal{L}_\lambda}\phi_\lambda)\,\phi_\lambda h\,\mathrm{d}\nu 
	=\int_{\mcal{X}} h\,\mathrm{d}\pi_\lambda. 
	\end{equation}
	Moreover, $\|Q_t^\lambda g\|_{L^2(\pi_\lambda)}\le e^{-ct}\|g\|_{L^2(\pi_\lambda)}$ for some $c >0$ due to $Q_t^\lambda = R_\lambda\, e^{t(\mcal{L}_\lambda - \Lambda(\lambda))} R_\lambda^{-1}$ and the corresponding property of $e^{t(\mcal{L}_\lambda - \Lambda(\lambda))}$ and $R_\lambda^{-1}g=\phi_\lambda g \perp \phi_\lambda$. We conclude that
	\begin{equation}
	\mathbb{E}^{\pi_\lambda}[A_T^2]
	=\frac{1}{T^2}\int_0^T \!\!\int_0^T \!\langle g, Q_{|t-s|}^\lambda g \rangle_{L^2(\pi_\lambda)} \,\mathrm{d}s \,\mathrm{d}t 
	= \frac{2}{T^2}\int_0^T \!(T-r)\,\langle g, Q_{r}^\lambda g \rangle_{L^2(\pi_\lambda)}\,\mathrm{d}r 
	\le\frac{2\|g\|^2_{L^2(\pi_\lambda)}}{cT}.
	\end{equation}
	Together with Chebyshev's inequality and the Borel-Cantelli lemma, this yields $A_{n^2}\to0$ as $n\to\infty$, almost surely.  If $n^2 \le T< (n+1)^2$, we use $A_T = \frac{n^2}{T}A_{n^2} + \frac{1}{T}\int_{n^2}^T g(Y_t^\lambda)\,\mathrm{d}t$ and
	\begin{equation}
	|A_T| \leq \frac{n^2}{T}\,|A_{n^2}| + \norm{g}_{\infty}\frac{(n+1)^2-n^2}{n^2}.
	\end{equation}
	Thus, $A_T\to0$ as $T\to\infty$, almost surely.
	
	\medskip
	
	Now, we let $Y^\lambda_0\!=\! x \!\in\!\mcal{X}$, and set $E_f = \{\lim_{T \to \infty} A_T \!=\! 0\}$ and $\mcal{X}_f = \{x \!\in\!\mcal{X}:\,\mathbb{Q}^\lambda_x(E_f)\!=\! 1\}$. We have shown that $\pi_\lambda(\mcal{X}_f)=1$. Moreover, part (b) and Assumption \ref{appass}(d) imply that for all $(t,x)\in(0,\infty)\times\mcal{X}$, the law of $Y_t^\lambda$ given $Y_0^\lambda=x$ is absolutely continuous with respect to~$\pi_\lambda$. Thus, 
	$Q_t^\lambda(\mathbf{1}_{\mcal{X}_f})(x)=1$, which by the Markov property of $Y^\lambda$ and the shift invariance of $E_f$ yields $\mathbb{P}(E_f|Y_0^\lambda=x)=1$. Finally, since $\mcal{X}$ is a separable metric space by Assumption~\ref{appass}(a), we can pick a countable collection $\mcal{U}$ of finite unions of open sets from a countable base. Then, given $Y_0^\lambda=x$, it holds with probability $1$ that $\rho_T(U)\to\pi_\lambda(U)$, $U\in\mcal{U}$. For any open $U\subset\mcal{X}$, there are $\mcal{U}\ni U_n\uparrow U$, and so on the same event of probability~$1$,
	\begin{equation}
	\liminf_{T \to \infty} \rho_T(U)\ge\liminf_{T \to \infty} {\rho_T(U_n)} = \pi_\lambda(U_n),\quad n\in\mathbb{N}.
	\end{equation}
	Taking $n\to\infty$, we arrive at the desired $\liminf_{T \to \infty} \rho_T(U)\ge\pi_\lambda(U)$. 
\end{proof}

The proofs in Sections 3--6 can be repeated with minor modifications, as outlined below. The only exception is the proof of convergence in Theorem \ref{MainTheorem}(b) for $\alpha>\frac{\beta}{\beta_c}$, which relies on the independent increments property of Brownian motion with drift and lacks a direct Markov process counterpart. The first minor modification in Sections~3--6 concerns the bound on the derivative of
\begin{equation}
G_\infty^x(\xi)=(0-\mathcal{L}_z)^{-1}V(x)
=\int_0^1 e^{s\mathcal{L}_z}V(x)\,\mathrm{d}s+e^{\mathcal{L}_z}(0-\mathcal{L}_z)^{-1}V(x)
\end{equation}
with respect to $\xi$ (the imaginary part of $z$) that is uniform over $0<\frac{\delta}{2}\le|\xi|\le R$ and $x\in K$ (see the end of Step 2 in the proof of Lemma \ref{CriticalAbounds}). For the first term on the right-hand side, we apply the Feynman-Kac formula \eqref{appFK}, differentiate with respect to $\xi$ inside the expectation, and use the boundedness of $V$. For the second term on the right-hand side, we recall that $(0-\mcal{L}_z)^{-1}$ is bounded-holomorphic on a neighborhood of the relevant segment. Thus, we may differentiate with respect to $\xi$ inside the expectation upon applying the Feynman-Kac formula \eqref{appFK}, then use the boundedness of $V$ and the local boundedness of $x\mapsto\norm{p_1(x,\cdot)}_{L^2}$ (see Assumption \ref{appass}(d)).

\medskip

Similarly, at the end of Step 3 in the proof of Lemma \ref{CriticalAbounds}, we need a bound on the derivative of 
\begin{equation}
	F_\infty^x(\xi)
	=\int_0^\infty e^{s\mcal{L}_z}Q_z^0V(x)\,\mathrm{d}s
	=\int_0^1 e^{s\mcal{L}_z}V(x)\,\mathrm{d}s
	-\int_0^1 e^{s\Lambda(z)}P_z^0V(x)\,\mathrm{d}s
	+e^{\mcal{L}_z}F_\infty^\cdot(\xi)(x)
\end{equation}
with respect to $\xi$ that is uniform over $|\xi|\leq\delta$ and $x\in K$. The first term on the right-hand side can be handled as in the case of $G^x_\infty(\xi)$. For the second term, we appeal to the holomorphicity of $z\mapsto\Lambda(z)$ and of $z\mapsto P^0_z$ from $L^2$ to $C(K)$ established after the proof of Theorem \ref{Generalspectrum}. Thus, we obtain the desired bound on the second term by differentiating with respect to $\xi$ inside the integral. Finally, the analogues of \eqref{F repr} and \eqref{Smallcontour}, together with the holomorphicity of $z\mapsto(\zeta-\mcal{L}_z)^{-1}$ and $z\mapsto Q^0_z$ on $L^2$, yield the holomorphicity of $F^\cdot_\infty$ on $L^2$. We conclude by again differentiating with respect to $\xi$ inside the expectation upon applying the Feynman-Kac formula \eqref{appFK}, then using the boundedness of $V$ and the local boundedness of $x\mapsto\norm{p_1(x,\cdot)}_{L^2}$.

\medskip

Lastly, at the end of Section 5, we need to show that for any compact $K_1,K_2\subset\mcal{X}$, 
\begin{equation}
\sup_{x \in K_1}\,\mathbb{E}_x\bigg[\int_0^\infty \mathbf{1}_{K_2}(X_t)\,\mathrm{d}t\bigg]<\infty.
\end{equation}
This follows from $P_t\mathbf{1}_{K_2}(x)\le\norm{p_1(x,\cdot)}_{L^2}\norm{P_{t-1}\mathbf{1}_{K_2}}_{L^2}$, $x\in K_1$, $t\ge1$, the local boundedness of $x\mapsto\norm{p_1(x,\cdot)}_{L^2}$ (see Assumption \ref{appass}(d)), and $\norm{P_{t-1}\mathbf{1}_{K_2}}_{L^2}\le e^{-\varepsilon(t-1)}\norm{\mathbf{1}_{K_2}}_{L^2}$. 

\medskip

Whether a semigroup $(e^{t\mcal{L}})_{t\ge0}$ is recurrent, or satisfies  $\sigma(\mcal{L})\subset(-\infty,-\varepsilon]$ for some $\varepsilon>0$, must be checked on a case-by-case basis. Here, it is worth noting that $\sigma(\mcal{L})\subset(-\infty,-\varepsilon]$ if and only if $-\langle \mcal{L} f , f \rangle_{L^2} \ge \varepsilon\norm{f}_{L^2}^2$ for all $f\in\mcal{D}(\mcal{L})$, which is usually easier to verify. We conclude by presenting two collections of processes that satisfy Assumption \ref{appass}.

\begin{lemma}
	\begin{enumerate}[(a)]
		\item Let $a \in C^2_b(\mathbb{R}^d;\mathbb{R}^{d\times d})$ be symmetric and uniformly elliptic (i.e.,  there exists some $c>0$ such that $y^\top a(x) y \ge c|y|^2$ for all $x, y \in \mathbb{R}^d$), and let $h \in C^1(\mathbb{R}^d)$. Set $\nu(\mathrm{d}x)=e^{2h(x)}\,\mathrm{d}x$, and assume that $\nu(\mathbb{R}^d)=\infty$. Suppose that the function $\frac{1}{2}\nabla \cdot a + a\nabla h$ is globally Lipschitz. Then, writing $\sigma$ for the positive definite square root of $a$, the SDE
		\begin{equation}\label{appsde}
			\mathrm{d}X_t 
			= \Big(\frac{1}{2}\nabla \cdot a + a\nabla h\Big)(X_t)\,\mathrm{d}t + \sigma(X_t)\,\mathrm{d}W_t
		\end{equation}
		is well-posed, and its solution $X$ satisfies Assumption \ref{appass}.
		\item Consider a countable state space $\mcal{X}$ with the discrete topology. Let $\nu\!:\mcal{X} \to (0,\infty)$ satisfy $\sum_{x \in \mcal{X}} \nu(x)=\infty$. Take $q(x,y)\ge0$ to satisfy the detailed balance condition $\nu(x)\, q(x,y) = \nu(y)\, q(y,x)$. Assume that the graph on $\mcal{X}$ with edges $\{(x,y):\, q(x,y)>0\}$ is connected and that 
		\begin{equation}
		N:= \sup_{x\in\mcal{X}}\,\sum_{y \neq x}q(x,y)<\infty.
		\end{equation}
		Then, the continuous time Markov jump process with the (bounded) generator
		\begin{equation}\label{jumpgen}
			\mcal{L}f(x) := \sum_{y \neq x} q(x,y)\,(f(y)-f(x))
		\end{equation}
		satisfies Assumption \ref{appass}.
	\end{enumerate}
\end{lemma}

\begin{proof}
	\textit{(a).} We first construct $X$ by means of Dirichlet forms. Start with a symmetric form on $C_c^\infty(\mathbb{R}^d)$ (dense in $L^2 = L^2(\mathbb{R}^d,\nu)$):
	\begin{equation}\label{predirichlet}
		\mcal{E}_0(f,g) := \frac{1}{2}\int_{\mathbb{R}^d} \nabla f^\top a \nabla g\,\mathrm{d}\nu.
	\end{equation}
	We claim that $\mcal{E}_0$ is closable. Indeed, take a sequence $(f_n)_{n\in\mathbb{N}}$ in $C_c^\infty(\mathbb{R}^d)$ such that $\|f_n\|_{L^2}\to0$ and $\mcal{E}_0(f_n-f_m,f_n-f_m)\to 0$. Then, the uniform ellipticity of $a$ implies that $(\nabla f_n)_{n\in\mathbb{N}}$ is Cauchy in $L^2$, and so $\nabla f_n\to G$ for some $G\in L^2$. 
	For any $\varphi \in C_c^\infty(\mathbb{R}^d)$, 
	\begin{equation}
	\int_{\mathbb{R}^d} G_i \,\varphi\,\mathrm{d}\nu 
	= \lim_{n\to\infty} \int_{\mathbb{R}^d}  \partial_{i}f_n \,\varphi\,\mathrm{d}\nu
	= - \lim_{n\to \infty} \int_{\mathbb{R}^d} f_n \,\partial_i (\varphi e^{2h})\, e^{-2h}\,\mathrm{d}\nu=0,\quad i=1,\,2,\,\ldots,\, d. 
	\end{equation}
	Thus, $G\equiv 0$, and we obtain $\mcal{E}_0(f_n,f_n)\to 0$ from $\nabla f_n^\top a \nabla f_n \le C|\nabla f_n|^2$. To show that $\mcal{E}_0$ is Markovian, we take $\delta>0$ and pick $\varphi_\delta\in C^\infty(\mathbb{R};(-\delta,1+\delta))$ such that $\varphi_\delta(x)= x$, $x \in [0,1]$ and $0\le \varphi_\delta' \le 1$. One checks directly that for any $f \in C_c^\infty(\mathbb{R}^d)$, $g_\delta := \varphi_\delta \circ f \in C_c^\infty(\mathbb{R}^d)$ and $\mcal{E}_0(g_\delta,g_\delta) \le \mcal{E}_0(f,f)$. By \cite[Theorem 3.1.2]{Fukushima1994-eo}, the smallest closed extension $\mcal{E}$ of $\mcal{E}_0$ is a regular Dirichlet form, and $C_c^\infty(\mathbb{R}^d)$ is a core for $\mcal{D}(\mcal{E})$. With $\mcal{E}_\alpha(f,g):=\mcal{E}(f,g)+\alpha\langle f,g\rangle_{L^2}$, we have the inequality $\norm{f}^2_{H^1(\mathbb{R}^d,\nu)}\le C\mcal{E}_1(f,f)$ for all $f \in C_c^\infty(\mathbb{R}^d)$. Since $C_c^\infty(\mathbb{R}^d)$ is a core for $\mcal{D}(\mcal{E})$, we have that $\mcal{D}(\mcal{E})\subset H^1(\mathbb{R}^d,\nu)$, and $\mcal{E}(f,g)$ is given by the right-hand side of \eqref{predirichlet} for all $f,g \in \mcal{D}(\mcal{E})$. In particular, $\mcal{E}$ is strongly local in the sense of \cite[p.~6]{Fukushima1994-eo}.
	
	\medskip
	
	Let $\mcal{L}$ be the generator corresponding to $\mcal{E}$ (see \cite[Theorem 1.3.1]{Fukushima1994-eo}). We seek to identify the action of $\mcal{L}$ on $C_c^\infty(\mathbb{R}^d)$. Define the operator $L$ on $C_c^\infty(\mathbb{R}^d)$ by $Lf = \frac{1}{2}e^{-2h}\nabla \cdot (e^{2h}a\nabla f) = \frac{1}{2}a : \nabla^2 f + \big(\frac{1}{2}\nabla \cdot a + a\nabla h\big) \cdot \nabla f$. Fix $f \in C_c^\infty(\mathbb{R}^d)$ and $g \in \mcal{D}(\mcal{E})$. Since $C_c^\infty(\mathbb{R}^d)$ is a core for $\mcal{D}(\mcal{E})$, we can take $g_n \in C_c^\infty(\mathbb{R}^d)$ such that $g_n \to g$ in $L^2$ and $\mcal{E}(g-g_n,g-g_n)\to 0$. Integration by parts yields
	\begin{equation}
	\mcal{E}(f,g) = \frac{1}{2}\lim_{n \rightarrow \infty}\langle a \nabla f,  \nabla g_n \rangle_{L^2} = - \lim_{n \rightarrow \infty} \langle L f, g_n \rangle_{L^2} = -\langle Lf , g \rangle_{L^2}. 
	\end{equation}
	It follows that $f \in \mcal{D}(\mcal{L})$ and $\mcal{L}f = Lf$. 
	
	Furthermore, by \cite[Theorems 7.2.2 and 4.5.3]{Fukushima1994-eo}, there exists a $\nu$-symmetric continuous Hunt process $(\Omega,(\mcal{F}_t)_{t \geq 0},(Y_t)_{t \geq 0},\{\mathbb{P}_x\}_{x \in {\mathbb R}^d})$ whose Dirichlet form is $\mcal{E}$. The process has no killing inside $\mathbb{R}^d$, but may explode a priori. 
	
	\medskip
	
	Next, we infer from $a \in C^2_b(\mathbb{R}^d;\mathbb{R}^{d\times d})$ that $\sigma$ is bounded and globally Lipschitz (see \cite[Section 27]{Varadhan1981-yf}). Since $\frac{1}{2}\nabla \cdot a + a\nabla h$ is assumed to be globally Lipschitz, well-posedness and non-explosivity of the SDE \eqref{appsde} are standard (see, e.g., \cite[Theorem~5.2.5, Theorem~5.2.9, and Proposition~5.3.20]{Karatzas2004-ky}). In view of \cite[Corollaries 5.4.8 and 5.4.9]{Karatzas2004-ky}, the corresponding local martingale problem is well-posed. Moreover, by \cite[Theorem 1.2]{Menozzi}, the solution of \eqref{appsde}  admits positive transition densities $(q_t(x,y))_{t > 0}$ which are continuous in $(x,y)$. We define $p_t(x,y) = e^{-2h(y)}q_t(x,y)$, which is also continuous in $(x,y)$ and positive.
	
	\medskip
	
	We aim to show that the Hunt process $Y$ has the same law as the solution of \eqref{appsde}. Let $g_N\in C_c^\infty(\mathbb{R}^d)$ satisfy $g_N(x)=1$ for $|x|\le N$ and $g_N(x)=0$ for $|x|\ge 2N$, let $\tau_N := \inf \{t\ge0\!:|Y_t| \ge N\}$, and let 
	\begin{equation}
	\mcal{D}_0 := \{x_i \, g_N(x):\, 1\le i\le d,\, N \in \mathbb{N}\}
	\cup\{x_ix_j \, g_N(x):\, 1 \le i,j\le d,\, N \in \mathbb{N}\}\subset \mcal{D}(\mcal{L}). 
	\end{equation}
	Then, \cite[proof of Theorem 5.2.2, display (5.2.26)]{Fukushima1994-eo} implies that for $\nu$-almost every $x\in\mathbb{R}^d$, the process $f(Y_{\cdot\wedge \tau_N}) - f(x) - \int_0^{\cdot\wedge \tau_N} \mcal{L}f(Y_s)\,\mathrm{d}s$ is a martingale for all $(f,N)\in \mcal{D}_0\times\mathbb{N}$. (Strictly speaking, the display yields the martingality outside of a quasi-exceptional set, but this shows the martingality for $\nu$-almost every $x$ by the paragraph below \cite[p.~66, display (2.1.6)]{Fukushima1994-eo}. The fact that a martingale AF, as defined on \cite[p.~242]{Fukushima1994-eo}, is a martingale in the usual sense follows from \cite[p.~222, property (A.2)]{Fukushima1994-eo}.) For such $x$, we conclude that the law of $Y_{\cdot \wedge \tau_N}$ under $\mathbb{P}_x$ agrees with the law of the solution to the SDE \eqref{appsde} started from~$x$ and stopped at $\inf \{t\ge0\!: |X_t| \ge N\}$ (cf.~\cite[proof of Proposition 5.4.6]{Karatzas2004-ky}). Thanks to the non-explosivity of \eqref{appsde}, the laws of $Y$ and $X$ then agree for $\nu$-almost every $x$.
	
	\medskip
	
	Let $(P_t)_{t\ge0}$ denote the (strongly continuous, self-adjoint) semigroup of $Y$. For $\nu$-almost every $x$, we have
	\begin{equation}
	P_tf(x) =  \mathbb{E}_x[f(X_t)] = \int_{\mathbb{R}^d} q_t(x,y)\, f(y)\,\mathrm{d}y 
	= \int_{\mathbb{R}^d} p_t(x,y)\, f(y)\,\nu(\mathrm{d}y), \quad (t,f)\in(0,\infty)\times L^2.
	\end{equation}
	The self-adjointness of $P_t$ implies that $p_t(x,y)=p_t(y,x)$ for $\nu\times\nu$-almost every $(x,y)\in\mathbb{R}^{2d}$. Due to the continuity of both sides, $p_t(x,y)=p_t(y,x)$ holds for all $(x,y)\in\mathbb{R}^{2d}$. We claim that $(P_t)_{t\ge0}$ is irreducible. Suppose for contradiction that there is an invariant $E\subset\mathbb{R}^d$ with $\nu(E),\nu(E^c)>0$. Pick any $E_1\subset E$, $E_2\subset E^c$ with $0<\nu(E_1),\nu(E_2)<\infty$. Since $E$ is invariant, $P_t \mathbf{1}_{E_1}\equiv 0$ on $E^c$, and so $\langle P_t \mathbf{1}_{E_1}, \mathbf{1}_{E_2} \rangle_{L^2}=0$. This contradicts the positivity of $p_t(x,y)$, yielding the irreducibility of $(P_t)_{t\ge0}$.
	
	\medskip
	
	It remains to verify parts (d) and (e) of Assumption \ref{appass}. By $p_t(x,y)=p_t(y,x)$ and the Chapman-Kolmogorov equation, $\langle p_t(x,\cdot),p_t(y,\cdot)\rangle_{L^2} = p_{2t}(x,y)$, $x,y\in\mathbb{R}^d$. Thus, the continuity of $p_{2t}$ implies that for any compact $K\subset\mathbb{R}^d$, we have $\sup_{x \in K} \norm{p_t(x,\cdot)}_{L^2}<\infty$. Moreover, $\norm{p_t(x,\cdot)-p_t(y,\cdot)}_{L^2}^2=p_{2t}(x,x)+p_{2t}(y,y) - 2p_{2t}(x,y)$, 
	which together with the continuity of $p_{2t}$ reveals that $\mathbb{R}^d\to L^2$, $x \mapsto p_t(x,\cdot)$ is continuous, and so $x\mapsto P_tf(x)$ is continuous for all $(t,f)\in(0,\infty)\times L^2$. 
	
	For part (e), we fix an $f\in L^\infty(\mathbb{R}^d)$ with compact support, take a non-empty, bounded, open set $U$ with $C^1$-boundary containing the support of $f$, and let $\mcal{G} \subset \mcal{D}(\mcal{L})$ be bounded in the graph norm. Then, $\mcal{E}(g,g)=-\langle \mcal{L}g,g\rangle_{L^2}\le \norm{\mcal{L}g}_{L^2}\norm{g}_{L^2}$ shows that $\mcal{G}$ is bounded in the norm associated with $\mcal{E}_1$. The uniform ellipticity of $a$ and the uniform positivity of $e^{2h}$ on $U$ imply that $\sup_{g \in \mcal{G}} \norm{g}_{H^1(U,\mathrm{d}x)}^2\le C\sup_{g \in \mcal{G}}\mcal{E}_1(g,g)<\infty$. The Rellich-Kondrachov Theorem (see \cite[Section 5.7, Theorem 1 and the remark following its proof]{evans2010partial}) reveals that for any sequence in $\mcal{G}$, there is a subsequence $(g_n)_{n\in\mathbb{N}}$ converging in $L^2(U,\mathrm{d}x)$ to some $g$. Since $f\in L^\infty(\mathbb{R}^d)$ with support contained in $U$, we have that $f g_n$ converges in $L^2(\mathbb{R}^d,\mathrm{d}x)$ to $fg$. This convergence also holds in $L^2(\mathbb{R}^d,\nu)$ because $h$ is locally bounded.
	
	\medskip
	
	\noindent\textit{(b).} Define a symmetric form on $L^2=L^2(\mcal{X},\nu)$ by
	\begin{equation}
	\mcal{E}(f,g)=\frac{1}{2}\sum_{x,y\in\mcal{X}} \nu(x)\, q(x,y)\,(f(y)-f(x))\,(g(y)-g(x)).
	\end{equation}
	Since $|\mcal{E}(f,g)|\le 2N\norm{f}_{L^2}\!\norm{g}_{L^2}$ by the Cauchy-Schwarz inequality and the detailed balance condition, $\mcal{E}$ is well-defined. Verifying directly that $\mcal{E}$ is a regular Dirichlet form on~$L^2$, we apply \cite[Theorem 7.2.1]{Fukushima1994-eo} to obtain a Hunt process $(\Omega,(\mcal{F}_t)_{t \geq 0},(Y_t)_{t \geq 0},\{\mathbb{P}_x\}_{x \in \mcal{X}})$. An elementary calculation shows that 
	\begin{equation}
	\mcal{E}(f,g) = - \langle \mcal{L}f, g \rangle_{L^2},\quad f,g\in L^2,
	\end{equation}
	where $\mcal{L}$ is defined by \eqref{jumpgen}, so that $\mcal{L}$ is the corresponding generator. 
	
	\medskip
	
	Define a transition matrix by $K(x,y) = q(x,y)/N$, $x \neq y$ and $K(x,x)=1-\sum_{y\neq x}K(x,y)$. Since $K = I + \mcal{L}/N$, $K$ is self-adjoint and the semigroup $(P_t)_{t \geq 0}$ associated with $\mcal{L}$ reads
	\begin{equation}
	P_t = e^{t\mcal{L}}=e^{-tN}e^{tNK}=e^{-tN} \sum_{n=0}^\infty \frac{(tN)^n}{n!}K^n,\quad t>0.
	\end{equation}
	Together with the connectivity of the graph on $\mcal{X}$ with edges $\{(x,y)\!: q(x,y)>0\}$, this implies $P_t(x,y)>0$, $x,y\in\mcal{X}$. The irreducibility of $(P_t)_{t \geq 0}$ now follows by the same argument as in part (a). Furthermore, Assumption \ref{appass}(d) holds due to the detailed balance condition and the Chapman-Kolmogorov equation, which lead to $\|p_t(x,\cdot)\|^2_{L^2}=p_{2t}(x,x)$, $x\in\mcal{X}$. Finally, to see Assumption \ref{appass}(e), note that for any $f\in L^\infty(\mcal{X})$ with non-empty compact support and any bounded $\mcal{G}\subset L^2$, the set $\{(f(x)g(x))_{x \in\mathrm{supp}(f)}\!: g \in \mcal{G}\}$ is bounded in $\mathbb{R}^{\mathrm{supp}(f)}$, and so is relatively compact.
\end{proof}

\bigskip

\bibliographystyle{plain}
\bibliography{bibliography.bib}   

\medskip\quad\medskip

\end{document}